\documentclass{amsart} \usepackage{amsmath,amssymb,bbm}
\usepackage[dvips]{epsfig}

\input xy

\xyoption{all}
\newcommand{\opname}[1]{\operatorname{#1}}
\newtheorem{definition}{Definition}
\newcommand{\bldsym}[1]{\boldsymbol{#1}}
\newcommand{\catname}[1]{\boldsymbol{\opname{{#1}}}}
\newcommand{\FSN}{\catname{FinSet_0}}
\newcommand{\MFSN}{\catname{MFinSet_0}}

\newcommand{\Set}{\catname{Set}}
\newcommand{\SetZ}{\catname{Set [ Z ]}}
\newcommand{\UR}{\mathbbm{R}^+ \otimes U(1)}
\newcommand{\MRp}{\mathbbm{R}^+ \otimes M}
\newcommand{\MRz}{(\mathbbm{R}^+ \otimes M)[\![z]\!]}
\newcommand{\URz}{(\mathbbm{R}^+ \otimes U(1))[\![z]\!]}
\newcommand{\Cz}{\mathbbm{C}[\![z]\!]}
\newcommand{\inprod}[1]{\left\langle{#1}\right\rangle}
\newtheorem{theorem}{Theorem}
\newtheorem{varexample}{Example}
\newenvironment{example}{\begin{varexample}\em}{\em\end{varexample}}
\newcommand{\mathd}{\mathrm{d}}
\newcommand{\mathi}{\mathrm{i}}
\newcommand{\mathe}{\mathrm{e}}
\newcommand{\Longrightarrowlim}{\mathop{\longrightarrow}\limits}
\newtheorem{varremark}{Remark}
\newenvironment{remark}{\begin{varremark}\em}{\em\end{varremark}}
\newcommand{\rightarrowlim}{\mathop{\rightarrow}\limits}
\newcommand{\longrightarrowlim}{\mathop{\longrightarrow}\limits}
\newcommand{\leftarrowlim}{\mathop{\leftarrow}\limits}
\newcommand{\longleftarrowlim}{\mathop{\longleftarrow}\limits}

\begin{document}
\raggedbottom

\title{Categorified Algebra and Quantum Mechanics} \author{Jeffrey Morton}
\address{University of California, Riverside} \email{morton@math.ucr.edu}

\begin{abstract}
  Interest in combinatorial interpretations of mathematical entities
  stems from the convenience of the concrete models they
  provide. Finding a bijective proof of a seemingly obscure identity
  can reveal unsuspected significance to it.  Finding a combinatorial
  model for some mathematical entity is a particular instance of the
  process called ``categorification''.  Examples include the
  interpretation of $\mathbbm{N}$ as the Burnside rig of the category
  of finite sets with product and coproduct, and the interpretation of
  $\mathbbm{N} [ x ]$ as the category of combinatorial species.  This
  has interesting applications to quantum mechanics, and in particular
  the quantum harmonic oscillator, via Joyal's ``species'', a new
  generalization called ``stuff types'', and operators between these,
  which can be represented as rudimentary Feynman diagrams for the
  oscillator.  In quantum mechanics, we want to represent states in an
  algebra over the complex numbers, and also want our Feynman diagrams
  to carry more structure than these ``stuff operators'' can do, and
  these turn out to be closely related.  We will show how to construct
  a combinatorial model for the quantum harmonic oscillator in which
  the group of phases, $U(1)$, plays a special role.  We describe a
  general notion of ``$M$-Stuff Types'' for any monoid $M$, and see
  that the case $M=U(1)$ provides an interpretation of time evolution
  in the combinatorial setting, as well as other quantum mechanical
  features of the harmonic oscillator.
\end{abstract}
\maketitle

\pagebreak

\tableofcontents

\pagebreak

\section{Introduction}

One reason for the success of category theory has been its ability to
describe relations between seemingly separate areas of mathematics.
Here, we will describe how category theory can be used to see
relationships between enumerative combinatorics and quantum
mechanics. Specifically, we examine an approach to understanding the
quantum harmonic oscillator by ``categorifying'' the Fock
representation of the Weyl algebra - which is the algebra of operators
on states of the oscillator.  This builds on work described by Baez
and Dolan\cite{finfeyn}.  Since we do not expect readers to be expert
in both quantum mechanics and category theory, we have tried to make
this paper as self-contained as possible.  Many definitions and
explanations will be well-known to the experts in each field, and are
arranged by section as much as possible to allow readers to skip
familiar material.

Categorification is best understood as the reverse of
``decategorification''.  This is a process which begins with some
category, and produces a structure for which isomorphisms in the
original category appear as equations between objects instead.
Categorification is the reverse process, replacing equations in some
mathematical setting with isomorphisms in some category in a
consistent - but possibly non-unique - way.  One example is the way we
can treat the category of finite sets as a categorification of the
natural numbers $\mathbbm{N}$.  The set $\mathbbm{N}$ can be seen as a
set of cardinalities indexing isomorphism classes of finite sets, and
which get their addition and multiplication from the categorical
coproduct and product on the category of finite sets.  We will see
further examples of a connection between decategorification and
cardinality.

Joyal {\cite{joyal2}} has described the category of {\textit{structure
types}}, which can be seen as the categorification of a certain ring
of power series.  These structure types play an important role in
enumerative combinatorics, in which ``generating functions'' of given
types of structures can be used in a purely algebraic way to count the
number of such structures of various sizes. These generating functions
are ``decategorified'' versions of structure types - or, equivalently,
cardinalities of them.  In section \ref{struct-type-sec}, we describe
how this works in more detail, and give some examples.

This leads, in section \ref{harmoscsec}, to the application to quantum
mechanics.  In that section, we describe briefly the quantum harmonic
oscillator.  It has a Hilbert space of states, and the Weyl algebra is
the algebra of operators on this space.  The Weyl algebra has a
representation as operators on Fock space - the space of formal power
series in one variable with a certain inner product.  Here, it is
generated by two operators - the creation and annihilation operators.
We show how this algebra can be categorified using the category of
structure types as a replacement for Fock space, and a certain class
of functors as the operators.

In section \ref{stuff-types-sec}, we find that structure types do not
have a rich enough structure to capture all properties of power
series.  In particular, they do not provide a natural way to treat
power series as functions which can be evaluated or composed in a way
which is compatible with the idea of cardinality for structure types.
To properly categorify these ideas, one can extend the notion of
structure type to a so-called ``stuff type''.  This makes use of the
fact that there are two ways of seeing structure types as functors.
One is as a functor taking each finite set to the set of structures of
a certain type which can be placed on them - the functor gives
``coefficients'' associated to finite sets.  The other point of view
treats structure types as ``bundles'' over the category of finite
sets, whose projections take structured sets to their underlying sets.
The latter point of view allows a larger class of ``total spaces'' for
the bundle - in fact, it can be any groupoid.  We describe a
classification of functors, and show how dropping the requirement of
faithfulness on the projection functor for the bundle leads to stuff
types.  In appendix \ref{sec:slicecategory} we provide more details,
showing how stuff types form a category of ``groupoids over finite
sets''.

Section \ref{stuff-quantum-sec} examines operations on stuff types
which are useful to the program of categorifying quantum mechanics.
The category of stuff types naturally gets an ``inner product'' on
objects by means of a pullback construction.  We then show this is a
categorification of the usual inner product on Fock space.  Then we
describe the equivalent of linear operators on Fock space - ``stuff
operators'', and show how they can be seen directly as categorified
matrices.  These can act on stuff types by a construction similar to
the inner product, as one might expect.  We then develop some
particular examples of stuff operators, namely the equvalent of the
creation and annihilation operators.  These provide a connection to
Feynman diagrams.

Then, in \ref{mstuff}, we introduce the idea of $M$-sets, labelled
with elements of a general monoid $M$, extend this from sets to
general groupoids, and describe an idea of cardinality for these.  Of
special interest for quantum mechanics is the group $U(1)$, the group
of phases.  We see that it is possible, by ``colouring'' sets
(interpreted as sets of quanta) with these phases, to recover more
features of the quantum harmonic oscillator.  In particular, we
describe an $M$-stuff operator which corresponds to time evolution of
a state without interactions.  We also demonstrate a connection
between stuff operators and Feynman diagrams.

Finally, we summarize the results, and suggest directions in which
this work could be extended.

\section{Structure Types\label{struct-type-sec}}

\subsection{Categorification of $\mathbbm{N}$ and $\mathbbm{N} [\![z]\!]$}

Before we can study structure types, we need to see how a category
with products and coproducts gives rise to a {\textbf{rig}}, which is
to say an object like a {\textit{ring}}, possibly without
{\textit{n}}egatives (see Appendix A for a more precise definition).
A simple example is the {\textit{free rig on no generators}} - that
is, $\mathbbm{N}$ (generated by the nonzero element $1$ under
addition, this has no extra generators or relations, so it has a
natural homomorphism into any rig).  Natural numbers are called by
this name because they arise naturally as {\textit{counting}} numbers
- namely, numbers we use to give the cardinality of some finite set of
things.  Bijections between finite sets are what make counting
possible (for example, bijections of fingers and sheep), so these
cardinalities are actually equivalence classes of finite sets under an
equivalence relation given by bijections.

This suggests looking the category of finite sets, with a cardinality
map given by taking sets to their isomorphism classes.  The
cardinality map turns this category into a rig.  Since its
decategorification is a rig, we say this category is an example of a
{\textit{2-rig}} - a category with a monoidal operation like
multiplication, and a coproduct structure giving addition (for a more
precise definition, see Appendix A).  The cartesian product of sets
gives the multiplication in $\mathbbm{N}$, and disjoint union gives
addition.  An analogous process makes sense for any 2-rig, but a
reverse process, starting with any rig, is more difficult.

For a more involved example, consider the problem of categorifying the
{\textit{free rig on one generator}}, $\mathbbm{N} [\![z]\!]$ (the rig
of polynomials in $z$ with natural number coefficients).  We can think
of this as a rig of functions from $\mathbbm{N}$ to $\mathbbm{N}$ in
at least two different ways.  One treats an element $f \in \mathbbm{N}
[\![z]\!]$ as the map taking $n$ to the natural number $f ( n )$.  The
second takes $n$ to the coefficient of the $n^{\opname{th}}$ power of
$z$ in $f ( z )$, denoted $f_n$.  It is interesting to compare how
multiplication of these formal power series is represented in each
representation.  In the first case, we have pointwise multiplication:
\[ f g ( n ) = f ( n ) g ( n ) \]
In the other, mulitplication looks like convolution:
\begin{equation}
  ( f g )_n = \sum_{k = 0}^n f_k g_{n - k} . \label{seriesmult}
\end{equation}

We'll return to the first way of looking at them later when we study
``stuff types''.  The second representation will be better for now,
because it treats power series as purely formal, rather as functions.
To categorify the rig $\mathbbm{N} [ z ]$ seen as the rig of functions
$f : \mathbbm{N} \rightarrow \mathbbm{N}$, we naturally expect to look
at a corresponding 2-rig of functors (see appendix
\ref{sec:catalgebra} for more on 2-rigs):

\begin{definition}
  A {\textbf{structure type}} is a functor from the category
  $\FSN$ whose objects are finite sets, and whose
  morphisms are the bijections{\footnote{$\FSN$ is
  Joyal's {\textbf{$B$}} (for ``bijection'').  It makes no difference from
  which universe we take these sets - a skeletal version of
  $\FSN$ consisting of only pure sets, one per
  cardinality, will do as well as any other for our purposes.}}, into the
  category $\Set$:
  \[ F : \FSN \rightarrow \Set \]
  These functors naturally form a category whose morphisms are the natural
  transformations $\alpha : F \rightarrow F'$.  We denote this category by
  $\SetZ$.
\end{definition}

We say that image of a set $S$ is the set of all ``structures of type
$F$'' which can be placed on $S$.  Now recall that structure types can
be seen as categorified formal power series.  The ``coefficients'' are
in the category of sets and maps than the groupoid of finite sets and
bijections (though, since all maps in $\FSN$ are bijections, so are
their images).  Since we allow the possibility of infinite sets as
coefficients, these are more general than power series.  There is no
loss in adopting this approach, except when it comes to taking
cardinalities - an issue we will consider when we discuss stuff types
(in section \ref{stufftypes}).

An example of a structure type is the type of ``graphs on finite sets
of vertices''.  So then the image of a given finite set $S$ is just
the set of all such graphs on $S$.  The morphisms in $\FSN$, $f : S
\rightarrow S'$ give $F ( f ) : F ( S ) \rightarrow F ( S' )$, which
are maps of the structures (graphs) on $S$ to those on $S'$.  These
maps are compatible with the given bijection of underlying elements:
they amount to consistent relabellings of all the vertices in all the
graphs according to the bijection $F$.  In particular, permutations of
$S$ give automorphisms of $F(S)$.

To take the cardinality of a structure type, let $F_n$ be the set of
$F$-structures on the finite set $n$ (we will elide the difference of
notation between a set and its cardinality).  Then the cardinality of
the structure type $F$ is the formal power series
\begin{equation}\label{eq:strtypecard}
 |F| = \sum_{n = 0}^{\infty} \frac{|F_n |z^n}{n!}
\end{equation} where $|F_n|$ is just the usual set cardinality of
$F_n$.  We will see later how this formula for the cardinality is a
manifestation of ``groupoid cardinality'', but note for now that the
formula for $|F|$ given above is known as the ``generating function''
for $F$-structures.  This is a well known and useful idea in
combinatorics (and generalizes considerably beyond this example - see,
e.g. \cite{wilf} and \cite{BLL} for more on the whole subject).

\begin{example} The simplest example is for the type $Z$, which we
call ``being a one-element set''.  This structure can be put on a
one-element set in just one way, and in no ways on any other.  The set
of all $Z$-structures on $S$ contains just $S$ itself if $S$ has
exactly one element, and is empty otherwise.  The cardinality of the
type $Z$ is easily seen to be just $z$.

Similarly, we have the type ``being an $n$-element set'', denoted by
$\frac{Z^n}{n!}$, since it has cardinality $\frac{z^n}{n!}$.
\end{example} 

Generating functions (cardinalities of structure types) can be used to
find cardinalities of types defined in terms of simpler types.  We
make this more precise with the following definition and theorem:

\begin{definition}
  Given two structure types $F$ and $G$, there are sum and product structure
  types $F + G$ and $F \cdot G$, defined as follows.  Putting an $F +
  G$-structure on a set $S$ consists of making a choice of $F$ or $G$, and
  putting a structure of that type on $S$.  Putting an $F \cdot G$-structure
  on $S$ consists of splitting $S$ into an ordered pair of disjoint subsets,
  then putting an $F$-structure on the first part and a $G$-structure on the
  second part.
\end{definition}

Conceptually, we associate addition with ``or'' (an $F$-structure or a
$G$-structure), and we associate multiplication with ``and'' (a
splitting into an $F$-structure and a $G$-structure).  This is similar
to the categorified notions of addition and multiplication for
$\mathbbm{N}$ as disjoint unions and cartesian products in $\Set$.
This reappears when we look at functors from $\FSN$ to $\Set$, and
allows us to categorify the algebraic operations on the rig
$\mathbbm{N}[\![z]\!]$ as well.

\begin{theorem}\label{thm:productthm}
  If $F$ and $G$ are two structure types, then $|F + G| = |F| + |G|$ and $|F
  \cdot G| = |F| \cdot |G|$
  
  \begin{proof}
    To see that $|F + G| = |F| + |G|$, just note that the set of $F +
    G$-structures on the set $n$ consists of the disjoint union of the set of
    $F$-structures and the set of $G$-structures, since by definition such a
    structure consists of a choice of $F$ or $G$ together with a structure of
    the chosen type.  Thus, the set cardinalities satisfy $| ( F + G )_n | =
    |F_n | + |G_n |$, from which the result follows from the definition by
    linearity.
    
    Now, as for $|F \cdot G|$, we have that
    \begin{eqnarray}
         |F \cdot G| ( z ) & = & \sum_{n = 0}^{\infty} \frac{| ( F \cdot G )_n
        |z^n}{n!}\\
         \nonumber & = & \sum_{n = 0}^{\infty} \sum_{k = 0}^n \frac{z^n}{n!}
        \binom{n}{k} |F_k | \cdot |G_{n - k} |\\
         \nonumber & = & \sum_{n = 0}^{\infty} \sum_{k = 0}^n \frac{z^n}{n!}
        \frac{n!}{k! ( n - k ) !} |F_k | \cdot |G_{n - k} |\\
         \nonumber & = & \sum_{n = 0}^{\infty} \sum_{k = 0}^n \frac{|F_k |z^k}{k!} \cdot
        \frac{|G_{n - k} |z^{n - k}}{( n - k ) !}\\
         \nonumber & = & |F| ( z ) \cdot |G| ( z )
    \end{eqnarray}
    This follows directly from the fact that the number of $F \cdot
    G$-structures on an $n$-element set is a sum over all $k$ from $0$ to $n$
    of a choice of a $k$-element subset of $n$, multiplied by the number of
    $F$-structures on the chosen $k$-element subset and of $G$-structures on
    the remaining $( n - k )$-element set.
  \end{proof}
\end{theorem}

These facts suffice to prove that the functor category $\SetZ$ is a
2-rig (in fact, it is the free 2-rig on one generator),

\begin{theorem}
  The category $\SetZ$ is a 2-rig whose monoidal
  operation $\otimes$ is the product $\cdot$ defined above.
  
  \begin{proof}
    See Appendix A.
  \end{proof}
\end{theorem}

The 2-rig structure of $\SetZ$ lets us find useful information about
types of structure defined in terms of simpler types using the sum and
product operations we have defined, such as many recursively-defined
structures.  In particular, it makes sense of many calculations done
with generating functions in combinatorics.  A simple example shows
that the factor of $\frac{1}{n!}$ in the cardinality formula is due to
the fact that we do not think of the set $n$ as ordered.

\begin{example} Define the type $Z^n$ using the product operation in
terms of the type $Z$, ``being a one-element set''.  In particular, we
say $Z^1=Z$, and recursively define $Z^n = Z \times Z^{n-1}$.  Then by
theorem \ref{thm:productthm}, we have $|Z^n| = z^n$.  Now we observe
that we can interpret $Z^n$ as the type of ``total orderings on an
$n$-element set''.  Since a structure in the product $Z \times
Z^{n-1}$ involves a choice of two distinguishable subsets, we can find
a unique total ordering on the set $n$ by assuming one to precede the
other, and defining the total order recursively.  So in fact a $Z^n$
structure can be seen as just a total order.  And, indeed, there are
$n!$ total orderings, so the type of the cardinality is $\frac{n! 
z^n}{n!} = z^n$.

Now note that the structure ``being a finite set'' is the sum over all
$n$ of the types ``being an $n$-element set'': $\sum_n
\frac{Z^n}{n!}$.  Since each coefficient is 1, the cardinality of this
type is $e^z$, so we denote the type by $E^Z$.

If $T$ is the structure ``being a totally ordered set'' (that is, an
$T$ structure on $S$ is a total ordering on $S$), then $|T_n| = n!$
so that:
\begin{equation}
  |T| = \sum_{n = 0}^{\infty} \frac{n!}{n!} z^n = \frac{1}{1 - z}
\end{equation}
\end{example}

In {\cite{BLL}}, it is shown that many tree-like structures can be
defined using the operations on structure types we have described.
Binary trees provide an elementary example.

\begin{example}
  \label{binarytrees}An example of a structure-type is the type
  \textit{Binary trees}, which we denote $B$.  To put a $B$-structure
  on a finite set $S$ is to make $S$ into the set of leaves of a
  binary tree.  This is a recursively-defined tree structure, which is
  either a bare node (leaf), or a node with two branches, where each
  branch is another binary tree.  That is: \[ B \cong Z + B^2 \] since
  the structure type $Z$ is the type ''being a one-element set'' (a
  leaf), and $B^2$ is the type which, put on some set of elements,
  divides them into two subsets and puts a $B$-structure on each one.
  Some typical binary trees are shown in figure \ref{fig:binarytrees}
  

\begin{figure}[h]
\begin{center}
\includegraphics{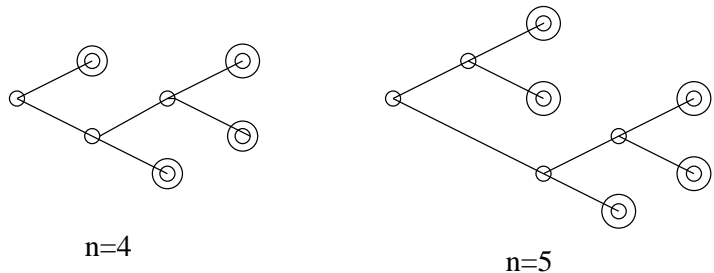}
\end{center}
\caption{\label{fig:binarytrees}Some Binary Trees with $n$ Leaves}
\end{figure}

  
  This highlights the close relationship between structure types and
  power series, since solving this recursive formula directly (for
  instance, by repeated substitution of the definition of $B$ into the
  $B^2$ term in the definition (or by solving the quadratic equation
  for $B$!) shows that $B$ is isomorphic to a structure type which we
  can write as the analog of a power series, beginning: \[ B \cong Z +
  Z^2 + 2 Z^3 + 5 Z^4 + \ldots \] where the coefficients are the
  Catalan numbers.  This enumerates binary trees of each size: $B$ is
  equivalent to a direct sum over all sizes $n$ of sets of some number
  of copies of the structure ``being an $n$-element set''.
  Specifically, this number is $n!$ times the $n^{\opname{th}}$
  Catalan number: the number of labellings of the leaves of an
  $n$-leaf binary tree by the elements of some given $n$-element set $S$.
\end{example}

\section{Structure Types and the Harmonic Oscillator} \label{harmoscsec}

We have defined structure types as functors of a certain kind (faithful
functors $\text{$F : \FSN \rightarrow
\Set$}$), and made the analogy between these and formal power
series.  Just as a single formal power series is really only of interest in
the setting of the space of all formal power series, so too a given structure
type acquires more meaning when we think of it in the setting of all such
functors.  In particular, one thing we are interested in for the purposes of
categorifying quantum mechanics are the categorified versions of algebras of
operators on this space.

To see why this is so, we will first describe the Weyl algebra - an algebra of
operators on the Hilbert space of states of the quantum harmonic oscillator -
then see how to use structure types to categorify it.  Readers who are familiar with quantum mechanics may wish to skip to section \ref{sec:strtypeweyl}

\subsection{The Quantum Harmonic Oscillator}

\subsubsection{The Weyl Algebra and its Representations}

States of the quantum harmonic oscillator can be represented as formal
power series with coefficients in $\mathbbm{C}$, where the state $z^n$
corresponds to the pure state with energy $n$, and a power series
represents a \textit{superposition} (linear combination) of these pure
states with given complex amplitudes.  We will study the harmonic
oscillator using the Weyl algebra, which consists of operators on the
space of states, and is generated by two operators satisfying certain
commutation relations.  There are two important representations the
Weyl algebra, which are easy to describe in terms of generators.
These are the Fock and Schr\"odinger representations.

\begin{definition}
  The {\textbf{Weyl algebra}} is the complex algebra generated by the
  {\textbf{ladder operators}}, namely the {\textbf{creation}} operator
  $a^{\ast}$ and the {\textbf{annihilation}} operator $a$.  These
  satisfy the relations $[a,a^{\ast} ] =a a^{\ast} -a^{\ast} a= 1$.
  The {\textbf{Fock representation}} of the Weyl algebra on the vector
  space of formal power series in $z$ is determined by the effect of
  these generators: \begin{equation}
     af ( z ) = \frac{\mathd f ( z )}{\mathd z}
  \end{equation}
  and
  \begin{equation}
     a^{\ast} f ( z ) = z f ( z )
  \end{equation}
\end{definition}

In other words, $a= \partial_z$, the derivative operator, and
$a^{\ast} = M_z$, the operator ``multiplication by $z$''.  It should
be clear that these satisfy the defining commutation relations, since
$[a,a^{\ast} ] ( z^k ) =a a^{\ast} ( z^k ) -a^{\ast} a( z^k ) =
\frac{\mathd}{\mathd z} ( z^{k + 1} ) - z ( k z^{k - 1} ) = z^k$.
These operators do not correspond to quantum-mechanical observables
(which must be self-adjoint); instead, these are operators which add
or remove a quantum of energy to a state.

We can build many operators from just $a$ and $a^{\ast}$.  One which
often appears is the operator
\begin{equation}\label{eq:fieldoperatordef}
  \phi = a + a^{\ast}
\end{equation} called the ``field operator'', which produces a
superposition of the states in which the system has lost one quantum
or gained one quantum of energy (in some interaction).

Another important operator is the {\textit{number operator}}, denoted
$N$, which is just
\begin{equation}
N=a^{\ast}a
\end{equation}
This is related to the energy of the system - which is the Hamiltonian
for the evolution of the oscillator.  It should be clear that the
eigenvalues of $N$ are just the natural numbers $\mathbbm{N}$, and an
eigenstate corresponding to the eigenvalue $n$ is the state $z^n$ (a
pure state with energy $n$).  These ladder operators give us the key
to seeing the Weyl algebra in the categorified setting, when we pass
from formal power series to structure types.

The other representation mentioned was the Schr\"odinger representation.
Here we think of the Weyl algebra as generated by a different pair of
generators:

\begin{definition}
  Two generators $p$ and $q$ of the Weyl algebra are
  given in terms of $a$ and $a^{\ast}$ by 
\begin{equation}
 \begin{array}{lll}
  q=\frac{a+a^{\ast}}{\sqrt{2}} &
  and &
  p=\frac{a-a^{\ast}}{\sqrt{2} \mathi}
 \end{array}
\end{equation}
or equivalently,
\begin{equation}
 \begin{array}{lll}
  a= \frac{q+ \mathi p}{\sqrt{2}} &
  and &
  a^{\ast} = \frac{q- \mathi p}{\sqrt{2}}
 \end{array}
\end{equation}
\end{definition}

Physically, $p$ is the momentum operator, whose eigenstates are
pure states with definite momenta which are the eigenvalues of $p$,
and $q$ is the position operator, whose eigenstates have position
given by the corresponding eigenvalues.  We can take these as the defining
generators of the Weyl algebra, since these generate everything in it (given
either pair of generators, the other is uniquely defined).  Moreover, they
satisfy the relations:
\begin{equation}
 [p,q] =p q - q p = -\mathi
\end{equation}
We could have taken this as the definition of the Weyl algebra, instead.

\begin{definition}
  A Schr\"odinger representation of the Weyl algebra is a
  representation on a space of functions $\psi : \mathbbm{R}
  \rightarrow \mathbbm{C}$, with the {\textbf{position}} operator $q$
  and {\textbf{momentum}} operator $p$ represented as \begin{equation}
   \begin{array}{lll}
    p \psi ( x ) = - i \psi' ( x ) &
    and &
    q \psi ( x ) = x \cdot \psi ( x )
   \end{array}
  \end{equation}
\end{definition}

The space of functions on which the $p$ and $q$ operators act is
commonly taken to be the Schwartz functions.  These are smooth
functions all of whose derivatives (including the functions
themselves) decay to zero faster than the reciprocal of any polynomial
(so in particular they are $L^2$ functions).

We can note that these $p$ and $q$ satisfy the commutation relations
above, by exactly the same argument as used for the ladder operators
in the Fock representation.  It is interesting, but potentially
confusing, that in both representations the generators can be
represented as multiplication and differentiation.  In fact, every
representation of the Weyl algebra has such a form, but the variables
in which power series are expanded will have different
interpretations.  The variable $x$ of functions in the Schr\"odinger
representation is literally the position variable for the oscillator,
whereas $z$ in the Fock representation is simply a marker, whose
exponent represents the energy of a state.

In fact, we will focus entirely on the Fock representation, and will
see that it has a natural combinatorial interpretation, which we will
describe in terms of the generators $a$ and $a^{\ast}$, and involving
structure types.  We can think of this interpretation as a
categorification of the Fock representation of the Weyl algebra.
Developed further, it will show close connections to the theory of
Feynman Diagrams, as we shall see in section \ref{sec:feynman}.

\subsubsection{The Inner Product on Fock Space} \label{fockinnerprod}

Now, we should remark here that in representing the quantum harmonic
oscillator, like any quantum system, the two essential formal entities
we need are a Hilbert space and an algebra of operators on this space,
including self-adjoint operators corresponding to the physical
observables of the system.  We have described the Weyl algebra as its
acts on the Hilbert space of power series, but we haven't completely
described this space, since to make it a Hilbert space, we need an
inner product.  Clearly, as a vector space, $H$ is spanned by $\{ z^k
|k \in \mathbbm{N} \}$, but there are many ways to put an inner
product on this - and from each one, we get a Hilbert space from the
subset $\mathbbm{C}[z]$ consisting of the elements with finite norm.
So we need to choose the physically significant inner product in order
to specify the Hilbert space of states of the oscillator.

The inner product represents the complex amplitude, whose squared norm
is a probability, for some combination of a state-preparation process
and a measurement process (a ``costate'').  The inner product of two
states $\inprod{\psi,\phi}$ is a ``transition amplitude'' the
amplitude for finding a system set up in state $\phi$ to be in state
$\psi$ on measurement.  If the system undergoes a change of state
between set-up and measurement, there will be an operator applied to
one side.  Self-adjoint operators correspond to \textit{observable
quantities}, whose eigenvalues are the possible values which can be
observed.

For physical reasons, we say that the inner product of two states
having different energy should be zero: $\inprod{z^n,z^m}=0 \text{ if
} n \neq m$.  That is, the probability of setting up a state with
energy $m$ and observing the state to have a different energy $n$
should be zero by conservation of energy, if there are no intermediate
interactions.

Since position and momentum are observable quantities, we want the
corresponding operators $p$ and $q$ we described along with the
Schr\"odinger representation - to be self-adjoint.  But then we must
have $a^{\ast} =a^{\dagger}$: i.e. the ladder operators must be
adjoints.  A straightforward calculation then reveals:
\begin{equation}
  \begin{array}{lll}
    \inprod{z^n,z^n} & = & \inprod{a^{\ast} z^{n
    - 1},z^n} \\
     & = & \inprod{z^{n - 1},az^n}\\
     & = & \inprod{z^{n - 1}, n \cdot z^{n - 1}} \\
     &  = & n \cdot \inprod{z^{n - 1}, z^{n - 1}}
  \end{array} \label{innerproductuncat}
\end{equation} Normalizing so that $\inprod{1, 1} = 1$ (here, the
vector $1$ represents the vacuum, or ground, state, where there are no
quanta of energy present), we get that the physically meaningful inner
product on power series is the following:

\begin{definition}
  The inner product on $\mathbbm{C}[z]$ is defined by its operation on
  the basis $\{z^n\}$ by $\inprod{z^n, z^m} = \delta_{n, m} n!$.  The
  space of states of the harmonic oscillator consists of all power
  series with finite norm according to this inner product.
\end{definition}

The particular form of this inner product turns out to be closely
involved with the connection between this quantum system and structure
types, and in fact this inner product will turn out to have a natural
categorification in the setting of ``stuff types'', which we see in
section \ref{stufftypes}.

A standard question one may ask about the quantum harmonic oscillator
is to find an ``expectation value'' like $\inprod{z^n,\phi^k z^m}$,
where.  Recall that the field operator $\phi$ takes a state $\psi$ and
gives a superposition of states in which it has gained and lost one
quantum of energy.  Thus, the expectation value may be nonzero, so
long as $k$ is at least as large as the difference between $n$ and
$m$.  This is the kind of value which can be calculated by means of
Feynman diagrams.

What we mean to show now is that the categorified expectation value
above has a direct interpretation in terms of a groupoid whose objects
just look like Feynman diagrams.  To see this, we first start with a
description of a categorified Weyl algebra.

\subsection{Structure Types and the Weyl Algebra\label{sec:strtypeweyl}}

We have said already that structure types are a categorified version
of formal power series (with natural number coefficients).  We want to
use them to help us categorify the Weyl algebra, which is generated by
$a=\frac{\partial}{\partial z}$ and $a^{\ast} = M_z$ in the algebra of
operators on such series.  To do this, we must develop operators on
structure types (that is, natural transformations between functors)
which correspond to these in the decategorified form.

This amounts to finding a \textit{combinatorial interpretation} of the
operators $a$ and $a^{\ast}$, since a structure type is a
combinatorial entity: it identifies a kind of structure which can be
put on a finite set, and decategorifying it gives a ``generating
series'' for those sorts of structures.  The coefficients of this
series count the number of such structures, which is the cardinality -
and thus decategorification - of the \textit{set} of such structures.
The combinatorial flavour of structure types is made clearer, for
example, when we can define them recursively, or otherwise show some
relationships between the structures on different sets. These
properties can sometimes be expressed as algebraic or differential
equations involving the generating series.

A pure state with $n$ quanta of energy, in the representation of the
Weyl algebra described above, corresponds to the state $z^n$.  The
categorified version of this state, $Z^n$, in the structure-type
setting, is the structure ``being a finite set with $n$ elements''
(that is, this is the structure which can be put in exactly one way on
an $n$ element set and no ways on any other set).  It seems natural to
identify the elements of the sets on which we put our structures with
\textit{quanta of energy} of the quantum harmonic oscillator, and this
is what we will do.  A categorified state $\Phi$ is a structure type -
a \textit{type} of structure which can be put on some set of quanta of
energy, which we can express in the basis $\{ Z^n \}$ of the 2-Hilbert
space of states ($\SetZ$).  It is characterized by the set of
$\Phi$-structures on each size of finite set, and the ways those
structures transform as we relabel the underlying set of quanta.

We want to define two operators on structure types, $A$ and
$A^{\ast}$, which correspond to differentiation with respect to $z$
and multiplication by $z$ in the Fock representation of the Weyl
algebra.  In fact, these are just the ``insertion'' and ``removal''
operators, familiar in combinatorics:

\begin{definition}\label{structladder}
  The structure-type operator $A$ acts on a structure type $F$ to give a
  structure type $A F$, for which putting an $A F$-structure on a set $S$ is
  to adjoin a new element, which we denote $\star$, to $S$ and then put an
  $F$-structure on $S \cup \{ \star \}$.  The adjoint operator, $A^{\ast}$ is
  the one which acts on a structure type $F$ by giving a structure type
  $A^{\ast} F$ for which putting an $A^{\ast} F$-structure on a set $S$ is the
  same as removing an element from $S$ and putting an $F$-structure on the
  resulting set.
\end{definition}

It should be clear that $A$ acts like differentiation, by seeing how
it acts on $Z^n$, since $A ( Z^n )$ is the structure which, to put it
on a set, means putting the structure of ``being a totally ordered
$n$-element set'' on $S \cup \{ \star \}$.  There are $n$ ways to do
this, provided $S$ is an $n - 1$ element set (one for each position
$\star$ might take in the total ordering), so $A ( Z^n )$ is
equivalent to putting the structure of ``being a totally ordered $( n
- 1 )$-element set on $S$ and also choosing one of $n$ positions:
\begin{equation}
  A ( Z^n ) \cong n \cdot Z^{n - 1}
\end{equation}
And in general, $A$ acts like differentiation on structure types (we can
extend the above property linearly).

Now to see how $A^{\ast}$ acts, note that since an $A^{\ast}
F$-structure on S is equivalent to a way of splitting S into two
parts, putting an $F$ structure on one, and the structure of being a
1-element set on the other, we can see even more directly that
$A^{\ast}$ acts as (categorified) multiplication by $Z$.  The
commutation relation $a a^{\ast} -a^{\ast} a = 1$ now can be seen as a
decategorification of the corresponding property of the operators $A$
and $A^{\ast}$.  To summarize:

\begin{theorem}
  The structure type operators $A$ and $A^{\ast}$ satisfy
  \begin{equation}
   A \circ A^{\ast} = A^{\ast} \circ A + 1
  \end{equation}
and
  \begin{equation}
   |A F| = \frac{\mathd}{\mathd z} |F|
  \end{equation}
and
  \begin{equation}
   |A^{\ast} F| = z|F|
  \end{equation}
\end{theorem}

To see these ideas more concretely, consider the following examples:

\begin{example}
  A categorified state is just a structure type.  One such type is
  ``being a finite set'', which we can denote $E^Z$ ($E$ from the
  French ``ensemble'', but also appropriate in light of the generating
  series for this type).  Then putting an $A ( E^Z )$-structure on a
  finite set $S$ is the same as putting an $E^Z$ structure on $S \cup
  \star$, and there is exactly one way to put an $E^Z$-structure on
  ANY finite set, hence exactly one way to put an $A ( E^Z
  )$-structure on $S$ as well.  That is, $A ( E^Z ) = E^Z$.  This
  makes sense, since $|E^Z | = \mathe^z$, and so this equation becomes
  $\frac{\partial}{\partial z} \mathe^z = \mathe^z$ when we take its
  cardinality.
\end{example}

\begin{example}
  Consider the structure type $O$, ``being a totally ordered finite set'' (or
  ``of total orderings on a finite set'').  There are $n!$ such structures on
  a finite set $n$, so we can see this type as isomorphic to
  \begin{equation}
    \frac{1}{1 - Z} = \sum_{n = 0}^{\infty} Z^n
  \end{equation} (since $Z^n$ is the type of totally ordered
  $n$-element sets). Here, the sum is to be understood as a coproduct.
  Then to put an $A ( O )$-structure on a set $S$, one puts a total
  order on $S \cup \{ \star \}$: this amounts to splitting $S$ into
  two (ordered) parts and putting a total order on each part.  That
  is, an $A ( O )$-structure is the same as an $O^2$-structure, which
  is a combinatorial interpretation of the algebraic fact that
  $\partial_z \frac{1}{1 - z} = \frac{1}{( 1 - z )^2}$.  Extending
  this further, an $A^2 ( O )$-structure on $S$ consists of putting a
  total order on $S$ with two extra elements adjoined - thereby
  dividing $S$ into three parts (in order) and totally ordering these.
  There are two ways to build such a structure with $O$- and $A ( O
  )$-structures: to divide $S$ in two and put an $A ( O )$-structure
  (equivalently, $O^2$-structure) on the first part and an
  $O$-structure on the second, or to do this in the reverse order.
  Thus, $A^2 ( O ) \cong A ( O ) \cdot O + O \cdot A ( O ) \cong O^3 +
  O^3$ - a combinatorial interpretation of the fact that $\partial_z^2
  ( \frac{1}{1 - z} ) = \frac{2}{( 1 - z )^3}$.  This pattern
  continues for general $A^n ( O )$.
\end{example}

This representation of the Weyl algebra in terms of these operators on
structure types gives us a model in which the ladder operators have immediate
meaning in terms of adding and removing elements to sets, and the states of
the system are given as kinds of structures which can be put on those sets. 
The sets in question are sets of {\textit{quanta of energy}} in the system.  A
categorified state consists of a type of structure which can be put on these
quanta.  The number of such structures for each number of quanta is related
to the amplitude for the state to have that energy by means of the inner
product on the Hilbert space given in equation \ref{innerproductuncat}.  We
shall see a combinatorial interpretation for this inner product in a later
section, and see that it has a direct relationship to Feynman diagrams for
energy quanta in the harmonic oscillator.  This suggests what we must do
next: we have recovered the basic structure of the Weyl algebra in a
categorified form, but so far we are missing some concepts which, though quite
natural in the decategorified setting, are difficult to see in this
combinatorial picture.  These are evaluation of power series at particular
points, and the inner product of the Hilbert space (which will become a
2-Hilbert space in the categorified setting, once we have described the inner
product). To describe these adequately, we shall examine a generalization of
structure types, which we call ``stuff types''.

\section{Stuff Types} \label{stuff-types-sec}

We have described a structure type as a functor $F \in \SetZ$, where
the image of each finite set $S$ is the set of structures ``of type
$F$'' which can be put on $S$, but it should be obvious that the
category $\Set$ is rather larger than we need, since for most sorts of
structures we can think of, almost all sets will not appear in the
image at all - most are not sets of $F$ structures on an $n$-element
set for any $F$ or $n$.  If we think of a category $\catname{X}$ whose objects
are precisely the structures of type $F$, and whose morphisms are
those maps which arise from bijections of the underlying sets, we have
a category better suited to $F$. This category is, in fact, a groupoid
(i.e. a category in which all morphisms are iso), since $\FSN$ is.

If we do this, however, we now have not sets of $F$ structures as objects, but
$F$-structures themselves, and there are many such structures corresponding to
each $n$.  So now it is more natural to think of $F$ as a functor:
\[ F : \catname{X} \rightarrow \FSN \]
where each object of $\catname{X}$ is taken to its underlying set, and each morphism to
the underlying bijection of sets.  So every morphism in the image,
$\text{$\FSN$}$, then comes from morphisms in $\catname{X}$ under
this $F$.  A functor with this property is called {\textbf{faithful}}.

Why should we make this changed to the definition of a structure type?  There
are at least two good reasons to take this approach.  One is that, while the
previous definition fit well with the view of formal power series in which we
are interested in finding the coefficient of the $n^{\opname{th}}$ power of $z$,
this definition allows us to think of $F$ as corresponding to a power series
which we evaluate at various (positive) real numbers to get other (positive)
real numbers.  To see how this works, we first remark that it will make sense
to think of a structure type $F$ being evaluated at a {\textit{groupoid}}, and so
we need some useful facts about these.

\subsection{Groupoids}

To begin with, we recall what kind of category we are dealing with here:

\begin{definition}
  A {\textbf{groupoid}} is a category in which all morphisms are
  invertible.
\end{definition}

A group is a special example of a groupoid, with only one object -
then the elemnts of the group correspond to morphisms of the groupoid.
Another special case of a groupoid is a groupoid which has only one
(identity) morphism per object - such a groupoid is just equivalent to
the \textit{set} of its objects.  General groupoids are different from
either extreme case, since they can have many objects and many
morphisms.  However, the idea that a set may have a cardinality leads
us to try to extend this idea to more general groupoids.

\subsubsection{Groupoid Cardinality}

There indeed is a notion of cardinality for any groupoid, which in
general can give any positive real number (though it may also be
divergent). The notion of cardinality is closely related to the idea
of ``decategorification''.  This is a process which takes a category
and gives the set of isomorphism classes of objects.  Similarly, there
is a notion of cardinality which takes a set and gives a number, one
which takes a structure type and gives a formal power series, and one
which takes a monoidal category and gives a monoid (a category with
extra structure producing a set with extra structure).

\begin{definition}
  The {\textbf{groupoid cardinality}} of a groupoid $\mathcal{G}$ is
  \begin{equation}
   |\mathcal{G}| = \sum_{[ x ] \in \underline{\mathcal{G}}} \frac{1}{|
     \opname{Aut} ( x ) |}
  \end{equation}
  where $\underline{\mathcal{G}}$ is the set of isomorphism classes of objects of
  $\mathcal{G}$.  We call a groupoid {\textbf{tame}} if this sum converges.
\end{definition}

That is, each isomorphism class of objects of $\mathcal{G}$
contributes a term inversely proportional to the size of the
automorphism group of a typical element.  Note that groupoid
cardinalities of finite groupoids are just positive rational numbers -
a finite sum of reciprocals of the sizes of finite groups.  However,
since a general groupoid may be infinite, its cardinality may be an
infinite sum. Thus, groupoid cardinalities can be any nonnegative real
number (including infinity).  It is worth noting that, as with sets,
this idea of cardinality agrees with two natural operations we can
perform on groupoids.  These are the disjoint union (sum) and product,
in the following sense:

\begin{theorem}
  If $\mathcal{G}$ and $\mathcal{G'}$ are tame groupoids, then so are
  $\mathcal{G}+\mathcal{G}'$ and $\mathcal{G} \times \mathcal{G}'$,
  and we have $| \mathcal{G} + \mathcal{G}' | = | \mathcal{G}
  | + | \mathcal{G}' |$ and $| \mathcal{G} \times \mathcal{G}' | = |
  \mathcal{G} | \times | \mathcal{G}' |$.  If $\mathcal{G}$ and
  $\mathcal{G'}$ are equivalent, $| \mathcal{G} | = | \mathcal{G'} |$.
\end{theorem}

\begin{proof}
  The groupoid $\mathcal{G} + \mathcal{G'}$ is the category whose set of
  objects is the disjoint union of the sets of objects of $\mathcal{G}$ and
  $\mathcal{G'}$, and since all morphisms are internal to these groupoids, so
  is the set of isomorphism classes of objects.  So the fact that $|
  \mathcal{G} + \mathcal{G}' | = | \mathcal{G} | + | \mathcal{G}' |$ follows
  directly from the definition.
  
  The groupoid $\mathcal{G} \times \mathcal{G}'$ has objects which are ordered
  pairs of objects from $\mathcal{G}$ and $\mathcal{G'}$, and morphisms
  likewise ordered pairs of morphisms, which are iso precisely when both
  elements of the pair are iso.  So the isomorphism classes of objects are
  again ordered pairs of isomorphism classes of objects from $\mathcal{G}$ and
  $\mathcal{G'}$.  The automorphism group of any object $( g, g' ) \in
  \mathcal{G} \times \mathcal{G}'$ is the direct product $\opname{Aut} ( g )
  \times \opname{Aut} ( g' )$, so
    \begin{eqnarray}
        |\mathcal{G} \times \mathcal{G}' | & = & \sum_{[ ( g, g' ) ] \in
       \mathcal{\underline{G \times \mathcal{G}'}}} \frac{1}{| \opname{Aut} ( g,
       g' ) |}\\
        \nonumber & = & \sum_{[ g ] \in \underline{\mathcal{G}}} \sum_{[ g' ] \in
       \underline{\mathcal{G'}}} \frac{1}{| \opname{Aut} ( g ) \times \opname{Aut}
       ( g' ) |}\\
        \nonumber & = & \left( \sum_{[ g ] \in \underline{\mathcal{G}}} \frac{1}{|
       \opname{Aut} ( g ) |} \right) \times \left( \sum_{[ g' ] \in
       \mathcal{\underline{G'}}} \frac{1}{| \opname{Aut} ( g' ) |} \right)\\
       \nonumber & = & | \mathcal{G} | \times | \mathcal{G}' |
    \end{eqnarray}
  To see that cardinality is preserved under equivalence, just note that there
  is a 1-1 correspondence between isomorphism classes of their objects, and
  equivalent objects have isomorphic automorphism groups, since an equivalence
  is a full, faithful, and essentially surjective functor.

\end{proof}

This allows us to begin to elide the distinction between groupoids and their
cardinalities, passing back and forth as convenient (part of the intent of
categorification, just as when we conflate $\mathbbm{N}$ and
$\FSN$).

We can make an analogy with sets here: in a set seen as a category
with only identity morphisms, an element of the set is exactly the
same as an isomorphism class of objects, so we might think of these
classes as ``elements'' of a groupoid $\mathcal{G}$.  In this case,
the cardinality function (decategorification) which we have described
gives us potentially fractional values for each ``element'' of a
groupoid $\mathcal{G}$, so we can think of a groupoid as a way of
getting a ``fractional'' set - at least from the point of view of this
cardinality function.  So we can think of a groupoid as an entity
whose cardinality is a nonnegative extended real number (given by a
possibly infinite sum of positive rationals), and so if a structure
type is an entity whose cardinality is a formal power series with
natural number coefficients, it should not be too surprising that we
can view this structure type as a function which takes a groupoid and
produces another groupoid, just as a power series can be evaluated at
a real number and yields another real number.  Moreover, although
groupoid cardinalities can diverge, creating a problem of
well-definedness, by passing to the categorified setting, we can
eliminate this problem, and deal directly with the groupoids instead.
The way we do this is to define
\begin{equation}
  F ( Z_0 ) = \sum_{n \in \mathbbm{N}} (F_n \times Z_0^n)/\,\!\!/S_n
\end{equation} where the sum is interpreted as a coproduct, $F_n$ is
the groupoid whose object set is the $n^{\opname{th}}$ coefficient of
$F$ (as a set of structures) and whose morphisms are as above (the
groupoid of $F$-structured finite sets is a direct sum of groupoids
$F_n$, of structures whose underlying sets have size $n$).  In this
definition, $Z_0^n$ is a product of $n$ copies of $Z_0$, and $S_n$ is
the permutation group on $n$ elements.  The quotient which appears
inside the sum (coproduct) is a weak quotient of the groupoid $F_n
\times Z_0^n$ by the group $S_n$.  This requires some explanation.
First we will describe groupoid-coloured sets, then we will describe
the construction of a weak quotient.

\subsubsection{Groupoid-Coloured Sets}\label{colouredsets}

For any groupoid $Z_0$, we can speak of ``$Z_0$-coloured'' sets.  It
is easier to understand what these are if we think of sets as
groupoids themselves.  In particular, given a set $S$, we can think of
it as a groupoid whose objects are the elements of $S$, and the only
morphisms present are the identity morphisms (which must exist by
definition).  This is clearly a groupoid.  Then a map from a set into a
groupoid is just a functor.  So:

\begin{definition}
  A \textbf{$Z_0$-coloured set} is a set $S$ equipped with a
  \textbf{colouring} map $c:S \rightarrow Z_0$.  \textbf{Maps of
  $Z_0$-coloured sets} in $\hom((S,c),(S',c'))$ are bijections
  $\sigma: S \rightarrow S'$ together with, for each $x \in S$, a
  morphism $f_x \in \hom(c(x),c'(\sigma(x))$.  That is,
  \begin{equation}
   \xymatrix{
    S \ar^{\sigma}[r] \ar_c[d] & S' \ar^{c'}[d] \\
    Z_0 \ar@2{->}[r]_{\{f_x\}} & Z_0
   }
  \end{equation}
\end{definition}

\begin{figure}[h]
\begin{center}
\includegraphics{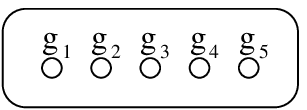}
\end{center}
\caption{\label{fig:colouredset}A $Z_0$-Coloured Set}
\end{figure}

We can think of these as sets having each element $x \in S$
``coloured'' by an object of $Z_0$, namely its image under the map
$f$, which is a functor between two groupoids.  The general form of
such a thing is shown in (\ref{eq:colouredset}), and an illustration
of a coloured set is shown in figure \ref{fig:colouredset}.

\begin{equation}\label{eq:colouredset}
 \xymatrix{
   Z_0 & *+[F-,]{\bullet \dots \bullet} \ar^f[l]
 }
\end{equation}

\begin{figure}[h]
\begin{center}
\includegraphics{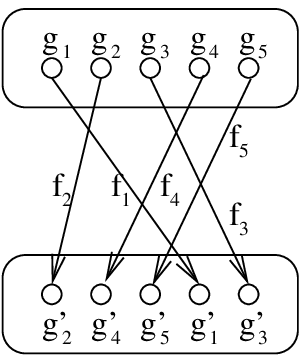}
\end{center}
\caption{\label{fig:colouredsetmap}A Morphism of $\mathcal{G}$-Coloured Sets}
\end{figure}

Morphisms of $Z_0$-coloured sets can be seen as bijections
$\alpha$ of the underlying sets where ``strands'' of the bijections
are labelled by morphisms of $Z_0$ between the $Z_0$-objects colouring
the elements of $S$ and $S'$ they connect, as shown in figure
\ref{fig:colouredsetmap}, where all the $g_i$ and $g'_i$ are in
$\mathcal{G}$, and $f_i \in \hom(g_i,g'_i)$.

\begin{equation}
\label{eq:colouredsetmap}
  \xymatrix{
    Z_0 \ar@2{->}_{\{f_i\}}[d] & *+[F-,]{\bullet \dots \bullet} \ar^{c}[l] \ar^{\alpha}[d] \\
    Z_0 & *+[F-,]{\bullet \dots \bullet} \ar^{c'}[l] \\
 }
\end{equation}

\begin{theorem} The collection of $Z_0$-coloured finite sets forms the
object set of a category $Z_0$-$\Set$ whose morphisms are as
described.  Moreover, $Z_0$-$\Set$ is a groupoid.
\begin{proof}
The morphisms described can be composed in the obvious way - composing
bijections of the sets, and labelling strands of the result with the
composites of the $Z_0$-morphisms labelling each strand:
\begin{equation}
  \xymatrix{
    S\ar^{\sigma}[r] \ar_{c}[d] \ar@/^1pc/^{\sigma ; \sigma'}[rr] & S'\ar^{\sigma'}[r] \ar_{c'}[d] & S'' \ar_{c''}[d] \\
    Z_0 \ar@2{->}[r]^{\{f_x\}} \ar@2{->}@/_1pc/_{\{f_x;f'_{x'}\}}[rr] & Z_0 \ar@2{->}[r]^{\{f'_x\}} & Z_0 \\
  }
\end{equation}
This notion of composition is well defined since if the strand
$\sigma(x) = x'$ is labelled by $f \in \hom(x,x')$ and $\sigma '(x') =
x''$ by $f' \in
\hom(x',x'')$, the strand in the composite $\sigma '\sigma(x)=x''$ can
be labelled by $f'f$ since these are composable.  This composition
rule inherits all the usual properties (e.g. associativity) from
bijections in $\Set$ and morphisms in $Z_0$.

The identity morphism from a $Z_0$-coloured set to itself is clearly
the morphism with identity bijection whose strands are labelled by
identity morphisms on the labels - again, properties of the identity
are inherited from $\Set$ and $Z_0$.  So in fact
$Z_0$-$\Set$ is a category.

Moreover all morphisms of this kind are invertible, since both
bijections $\sigma$ and all morphisms from $Z_0$ labelling strands are
invertible (i.e. $Z_0$-$\Set$ inherits the property of being a
groupoid from the fact that both $Z_0$ and $\Set$ are.
\end{proof}
\end{theorem}

In short, $Z_0$-$\Set$ is a groupoid of sets labelled by
objects of $Z_0$ in a way compatible with the groupoid
structure of $Z_0$.

\begin{remark}
  In the special case where $Z_0$ is a trivial groupoid with only
  identity morphisms, (which can be seen as a set) the definition of a
  morphism reduces to bijections compatible with the colourings, that
  is $\sigma : S \rightarrow S'$ gives a morphism between the
  $Z_0$-coloured sets $c : S \rightarrow Z_0$ and $c' : S' \rightarrow
  Z_0$ provided $\phi \circ c' = c$), i.e. that:
  \begin{equation}
    \xymatrix{    
      S\ar_{c}[d]\ar^{\sigma}[r] & S'\ar^{c'}[dl] \\
      Z_0 & \\
    }
  \end{equation}
  commutes.

  This is an example of an ``over category'', also known as a
  ``slice'' category.  See appendix \ref{sec:slicecategory} for more
  comments on this.
\end{remark}

Having constructed a groupoid $Z_0$-$\Set$ from $Z_0$, it makes sense
to ask about its cardinality.  However, this is a special case of what
we really wish to do: given a stuff type $\Phi$, find the cardinality
of $\Phi$-stuffed, $Z_0$-coloured finite sets.  In particular, if
$\Phi$ is the stuff type (in fact, structure type) ``being a finite
set'', then this is exactly the cardinality of $Z_0$-$\Set$.  To
describe the general case, we need to understand the weak quotient of
groupoids by groups, which will account for the effect of the
permutations $\sigma$ in the above construction.

\subsubsection{Weak Quotients of Groupoids by Groups}

We want to define the weak quotient of any groupoid $\mathcal{G}$ by a group
$G$ which acts on it, giving a groupoid $\mathcal{G} / G$.  This will be
particularly nice in the special case where $\mathcal{G}$ is just a trivial
groupoid - i.e. a set $S$, seen as a category - and in this case we can also
speak of the quotient of a set by a group, which will be a groupoid $S / G$.

\begin{definition}
  A {\textbf{strict action}} of a group $G$ on a category $C$ is a map $A$
  which for every $g \in G$ gives a functor $A ( g ) : C \rightarrow C$ such
  that $A ( g h ) = A ( g ) A ( h )$ and $A ( 1 ) = \opname{Id}_C$.  If there is
  such a strict action, then the {\textbf{strict quotient}} of $C$ by $G$ is
  a category $C / G$ together with a quotient functor $j : C \rightarrow C /
  G$ such that $j \circ A ( g ) = j$ for all $g \in G$.
\end{definition}

Clearly, in the special case where $\mathcal{G}$ is a groupoid, so is
$\mathcal{G}/\,\!\!/G$, since all morphisms are generated, by
composition, from invertible morphisms.  Moreover, when $C$ is a
trivial groupoid - i.e., a set - the definition of a strict action of
$G$ on $C$ is identical to the usual definition of a group action on a
set, but when there are nontrivial morphisms, it carries more
information because $A(g)$ must be functorial.  A strict quotient
agrees with the usual intuition of how a quotient should work on
objects, in that if there is a group element taking an object $x \in
C$ to $y \in C$, then $j(x)=j(y)$: the objects of $C / G$ are just
equivalence classes of objects in $C$ which are equivalent if they lie
in the same orbit under the action $A$.  Also, the fact that $A$ and
$j$ are functorial means that morphisms of $C$ are taken to morphisms
of $C / G$ compatibly with composition.

This strict action, and strict quotient, require too much to be very useful. 
More generally, it is not necessary that $j \circ A ( g )$ actually be equal
to $j$, so long as they are isomorphic in a reasonable way.

\begin{definition}\label{def:weakquotient}
  A {\textbf{weak quotient}} of a category $C$ by a group $G$, acting
  on $C$ by an action $A$ (as above), is a groupoid $C/\,\!\!/G$ whose
  objects are objects of $C$.  Its morphisms are generated by
  composition from morphisms in $C$ with morphisms of the form
  $A(h,g): g \rightarrow g'$, whenever $A(h)(g)=g'$, where $g, g' \in
  C$, $h \in G$.  Any relations which hold in $C$ hold in
  $C/\,\!\!/G$, together with relations: $A(h',g')\circ A(h,g)=
  A(h'h,g)$ whenever $A(h)g=g'$; and $A(h)f \circ A(h,g'') = A(h,g)
  \circ f$ for all $f: g'' \rightarrow g$ (the action $A$ is
  functorial).
\end{definition}

Notice that all the morphisms added by the group action are
invertible, since $h \in G$ is invertible.  Now we illustrate the weak
quotent $C/\,\!\!/G$ showing a few representative morphisms:
\begin{equation}
 \xymatrix{
  \bullet_{g}\ar@/_3pc/[rrrr]_{A(h,g)} & \bullet_{g''}\ar@/_1pc/[l]^f\ar@/^1pc/[rrr]^{A(h)f \circ A(h,g'') = A(h,g) \circ f} \ar@/_1pc/[r]_{A(h,g'')} & \bullet \ar@/_1pc/[rr]_{A(h)f} & \bullet & \bullet_{g'} \\
 }
\end{equation}

\begin{remark}
It is possible to describe a weak quotient by means of a (weak)
universal property.  It will be a groupoid $C /\,\!\!/ G$ together
with a quotient functor $j : C \rightarrow C /\,\!\!/ G$ such that
there is a natural isomorphism $\tau ( g ) : j \circ A ( g )
\Longrightarrowlim^{\sim} j$ for all $g \in G$.  We require that the
natural isomporphism satisfy the coherence condition $\tau ( g h ) =
\tau ( g ) \circ \tau ( h )$, and that the weak quotient should be
``weakly initial'' among all groupoids with these properties.  We will
not describe this in detail here, however, since we have a concrete
construction.
\end{remark}

Once we have this concept of the weak quotient of a groupoid by a
group, we naturally want to find its groupoid cardinality; this will
generally be smaller than $|\mathcal{G}|$ since we have added new
isomorphisms, hence potentially increased $|\opname{Aut}(x)|$ for some
objects $x$.  In fact, we have a better result:

\begin{theorem} The cardinality of the weak quotient of a groupoid
$\mathcal{G}$ by a group $G$ satisfies
\begin{equation}
  |\mathcal{G}/\,\!\!/G| = \frac{|\mathcal{G}|}{|G|}
\end{equation}
\begin{proof}
We have by definition that
\begin{equation}
  |\mathcal{G}/\,\!\!/G| = \left( \sum_{[ g ] \in \underline{\mathcal{G}/\,\!\!/G}} \frac{1}{|\opname{Aut}(g)|\cdot|\opname{Stab}(g)|} \right)
\end{equation} where $\opname{Stab}(g)$ is the stabilizer subgroup of
a representative $g$ in $\mathcal{G}$.  This is since the isomorphism
classes in $\underline{\mathcal{G}/\,\!\!/G}$ are given by considering
the isomorphism classes in $\mathcal{G}$ and identifying any which are
related by the action of $G$.  For each of these isomorphism classes,
the automorphism group consists of transformations taking one
equivalent object to another.  Any given object $[g]$ in such a class
will have as automorphism group the product of the automorphism of a
corresponding object in $\mathcal{G}$ with its stabilizer subgroup
$\opname{Stab}([g])$ in $G$.  So each isomorphism class contributes a
term $\frac{1}{|\opname{Aut}(g)|\cdot|\opname{Stab}(g)|}$.

On the other hand,
\begin{equation}
 \frac{|\mathcal{G}|}{|G|} = \frac{1}{|G|} \sum_{g \in \underline{\mathcal{G}}} \frac{1}{|\opname{Aut}(g)|}
\end{equation} Here the isomorphism classes are in $\mathcal{G}$: for
each isomorphism class in $\mathcal{G}/\,\!\!/G$, there will be
$|G|/|\opname{Stab}(g)|$ isomorphism classes in $\mathcal{G}$, since each
object $g \in \mathcal{G}$ is acted on by each element of $G$ and
taken to one of these classes.  So in fact this is the same as
$|\mathcal{G}/\,\!\!/G|$.
\end{proof}
\end{theorem}

\begin{remark} It is worth noting what happens in the special case
where $\mathcal{G}$ is just a discrete groupoid (i.e. a set, whose
groupoid cardinality is just its set cardinality).  If the group
action happens to be free, this result just says that the number of
orbits is the cardinality of the set divided by the size of the group
- but the result holds even when the action is not free, as in
this picture illustrating a $\mathbb{Z}_2$ action on a 3-element set,
giving a groupoid with cardinality $\frac{3}{2}$:

\vspace{0.2in}

\begin{equation}
 \xymatrix{
  \bullet\ar@/^2pc/@{<->}[rr]\ar@(dr,dl)[]^{id} & \bullet\ar@(dr,dl)[]^{id}\ar@{<->}
@(ul,ur)[] & \bullet\ar@(dr,dl)[]^{id} \\
 }
\end{equation}
\end{remark}

\subsection{Stuff Types as a Generalization of Structure Types}

We are now ready to describe stuff types and some of their properties.

\subsubsection{Stuff Types}\label{stufftypes}

We have already described \textit{structure types} as functors
\[ F : \FSN \rightarrow \Set \]
and also as special functors (in particular, faithful ones)
\[ \tilde{F} : \catname{X} \rightarrow \FSN \]
where $\catname{X}$ is the groupoid of $F$-structures on finite sets, and the functor
$\tilde{F}$ takes each to its underlying set. Faithfulness means that we do
not have two distinct morphisms in $\catname{X}$ with the same image in
$\FSN$ - that is, that a map between structures is
completely determined by its effects on the underlying set.  So we say that we
have ``forgotten structure'' - since there may be morphisms in
$\FSN$ which do not correspond to any in $\catname{X}$ (because
they do not ``preserve $F$-structure'').  In section \ref{sec:forgetful}, we
describe in more detail what is meant by saying that a functor forgets
\textit{properties}, \textit{structure}, and \textit{stuff}.  For the moment,
we extend the notion of structure types, and describe property types,
structure types, and stuff types, all of which are functors from some groupoid
$\catname{X}$ to $\FSN$ which forget the given sort of
information.

\begin{definition}
  A {\textbf{stuff type}} is a functor $\Phi : \catname{X} \rightarrow
  \FSN$, where $\catname{X}$ is a groupoid.  If $\Phi$ is faithful
  but perhaps not full or essentially surjective (forgets structure) we say it
  is a {\textbf{structure type}}; if it is full and faithful, but perhaps
  not essentially surjective (forgets properties) it is a {\textbf{property
  type}}; if it is an equivalence (forgets nothing), it is a
  {\textbf{vacuous property type}}.
\end{definition}

We think of these functors as giving the ``underlying sets'' of the stuff type
in question.  We can think of a stuff type $\text{$\Phi : \catname{X} \rightarrow
\FSN$}$ as a ``groupoid over $\FSN$'', so that any object of the groupoid $\catname{X}$ 
has an underlying set in $\FSN$, and by analogy with the terminology
``$F$-structured finite set'', we will describe it as a
``$\Phi$-stuffed finite set''.  This should suggest the idea that a
stuff type gives us a collection of objects which correspond to finite
sets $S$, but which possibly have extra information associated with
them (the ``stuff'' forgotten by $\Phi$).

Returning to our connection with the quantum harmonic oscillator, if
we view stuff types this way, and think of the dual entity,
$\Phi^{\ast} : \FSN \rightarrow \Set$, taking each finite set
$\bldsym{n}$ to the set of ways of putting ``$\Phi$-stuff'' on it - in
this context, this represents the set of ways for a state with energy
$n$ to occur in $X$, the groupoid associated with the state (which,
however, we just denote $\Phi$).

We can see how this connection develops in two steps.  First, we
replace a function $\psi :\mathbbm{N} \rightarrow \mathbbm{C}$
function $\psi : \mathbbm{N} \rightarrow \Set$: that is, a structure
type described as a map giving, for any finite set, the set of all
$\Psi$-structures on it.  The category of sets takes the role of the
complex numbers (in the first stage), which we understand in the
following way.  The complex values $\psi_n$, which in quantum
mechanics represent the ``amplitude'' for the particular pure state
with energy $n$ to occur are replaced in this first stage by sets of
$\Psi$-structures on a particular finite set $S$ of size $n$.  We can
think of these as ``the set of ways for possibility $n$ to occur''.
(In section \ref{mstuff} we will see how amplitudes arise).

\begin{figure}[h]
\begin{center}
\includegraphics{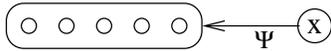}
\end{center}
\caption{\label{fig:typeobject}Example Object of a Stuff Type $\Psi$}
\end{figure}

The net step is to replace the natural numbers $\mathbbm{N}$ by the
groupoid $\FSN$, so that stuff types get a richer structure - in
particular, we can describe them as functors, since $\FSN$ is a
groupoid rather than merely a set.  Now we can think of these in
several ways, including as {\textit{bundles}} over the groupoid
$\FSN$.  In this setting, a ``point'' in the type $\Phi$ is an object
of $\catname{X}$ together with its underlying set - the base point in
the bundle.  We we can also think of the ``point'' as a finite set
with extra ``$\Psi$-stuff'', which we depict as a label $x$ attached
to $\Psi(x)$.  This is illustrated in figure \ref{fig:typeobject}.

If $\Psi$ is a structure type, this label can be seen as some
structure put on the set $\Psi(x)$, but in general this will not be
the case - the label could, for example, be other sets in a tuple of
which $\Psi(x)$ is one part (this is example \ref{st-ex:tuple}).  To
be general, we will not specify what this label contains, and simply
say that it contains ``stuff''.

\begin{figure}[h]
\begin{center}
\includegraphics{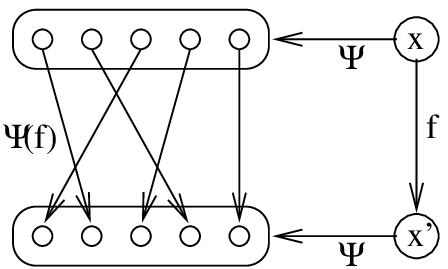}
\end{center}
\caption{\label{fig:typemorphism}Morphism In the Groupoid of Stuff Type $\Psi$}
\end{figure}

A morphism $f$ in $\catname{X}$ is iso, since $\catname{X}$ is a
groupoid: this gives a bijection $\Psi(f)$ of underlying sets in
$\FSN$.  Together, this data is a morphism in $\Psi$, show generally
in (\ref{eq:typemorphism}) and illustrated by the example in figure
\ref{fig:typemorphism}.
\begin{equation}\label{eq:typemorphism}
 \xymatrix{
  *+[F-,]{\bullet \dots \bullet}\ar_-{\Psi(f)}^{\wr}[d] & *++[o][F-]{x}\ar_-{f}^{\wr}[d]\ar^-{\Psi}[l] \\
  *+[F-,]{\bullet \dots \bullet} & *++[o][F-]{x'}\ar^-{\Psi}[l] \\
 }
\end{equation}

\begin{example}
  We've seen that the structure type ``being a finite set'' can be
  represented as $E^Z$, since there is exactly one way to put this
  structure on any finite set $S$ (in fact, it is the ``vacuous
  property type'').  We might ask if there is a structure type
  ``$E^{E^Z}$'' - shorthand for ``being a finite set of finite sets''.
  But, in fact, this is impossible.  Certainly, this will not have
  cardinality $\mathe^{\mathe^z}$, since there are an infinite number
  of ways of putting this structure on a given finite set $S$ - that
  is, taking a finite set of (disjoint) finite sets such that their
  union is $S$. In particular, any number of copies of the empty set
  may be in an ``$E^{E^Z}$-structure'' on $S$.  A worse problem is
  that there are distinct morphisms of such structure, differing only
  in their effect on these empty sets, which correspond to the same
  morphism on the underlying set $S$.  So these empty sets constitute
  extra ``stuff'', and this $E^{E^Z}$ must be a stuff type which is
  not a structure type.
\end{example}

This example, $E^{E^Z}$, highlights two points which we should address
about stuff types.  First, it is a special case of something we would
like to be able to do more generally, namely given two types $F$ and
$G$, to find a type $F \circ G$ which is in some reasonable sense the
composite of $F$ and $G$.  In this example, both $F$ and $G$ are the
vacuuous property type $E^Z$, yet even in this elementary case, their
composite already needs to be a stuff type.  The second issue is
cardinality: our concept of composition should include the fact that
when we take cardinalities, we should have that $|F \circ G| = |F|
\circ |G|$.  This does not work with the cardinality operation we have
developed for structure types - at least when we try to generalize it
naively.  The stuff type we are optimistically calling $F \circ G =
E^{E^Z}$ is well defined, for any finite set $S$, but the set of these
structures has infinite cardinality, so the corresponding generating
function is not $\mathe^{\mathe^z}$.  So we need a new concept of
cardinality for stuff types, which should be consistent with the old
definition in the special case where they are structure types.  We
will deal with this issue first.

\subsubsection{Cardinality of Stuff Types}

We have said that a stuff type is a groupoid over
$\FSN$.  By analogy with the situation for structure
types, then, we can define a notion of cardinality for stuff types.  For
structure types, we had a notion of cardinality for which 
\begin{equation}
 |F| = \sum_{n \in \underline{\FSN}} \frac{|F_n |}{n!} z^n
\end{equation}
Here, as before, $\underline{\FSN}$ is the set of
isomorphism classes of $\FSN$ (that is, $\mathbbm{N}$),
and $|F_n |$ is the cardinality of the set of $F$-structures on $n$ - where
this is now the usual cardinality function on sets.  The use of a common
notation for both cardinality operators emphasizes that any notion of
cardinality $| \cdot |$ is in some sense a decategorification operation.  For
simplicity, we will replace $\underline{\FSN}$ with
$\mathbbm{N}$, but we should remember that it derives from the original
groupoid $\FSN$.

In the same way, we can define the cardinality of a stuff type:

\begin{definition}\label{def:stuffcard}
  Given a stuff type $\text{$\Phi : \catname{X} \rightarrow
  \FSN$}$, we define its cardinality by
  \begin{equation}
    | \Phi | = \sum_{n \in \underline{\mathbbm{N}}} | \Phi_n |z^n
  \end{equation}
  where $|\Phi_n|$ is the {\textit{groupoid cardinality}} of the preimage of $n$ under
  $\Phi$.
\end{definition}

\begin{remark}
  Note that this really only makes sense if we are working with a
  ``skeletal'' version of $\FSN$, so that the preimage of a set of
  size $n$ is well defined.  This is a category which is equivalent as
  a category to the standard version of $\FSN$ from set theory - that
  is, there are functors between them whose composites are naturally
  isomorphic to the identities.  However, the skeletal version of
  $\FSN$ has only one object per isomorphism class - that is, per
  finite cardinal number.  We will assume $\FSN$ is skeletal whenever
  convenient.
\end{remark}

\begin{remark}
  Definition (\ref{def:stuffcard}) is consistent with the definition
  of structure types, since if we think of a structure type as a
  groupoid over $\FSN$, its morphisms all arise from permutations of
  the underlying sets, and so the groupoid cardinality of the preimage
  of $n$ is precisely $| \Phi_n | / n!$, since the groupoid is
  equivalent to the weak quotient $\Phi_n /\,\!\!/ S_n$.  Just as
  structure types had cardinalities which were power series whose
  $n^{\opname{th}}$ coefficients were integers divided by $n!$, stuff
  types have cardinalities which are power series with nonnegative
  real coefficients.
\end{remark}

An example shows that stuff types significantly generalize structure types, at
least as concerns the kinds of generating functions (power series) which can
appear as their cardinalities:

\begin{example}
  The stuff type ``being the first of a $k$-tuple of equal-sized
  finite sets'', for any value $k > 1$, is a stuff type which is not a
  structure type.  This can be represented as a functor $\Phi : \catname{X}
  \rightarrow \FSN$, where $\catname{X}$ is the groupoid whose
  objects are $k$-tuples of finite sets of the same size, and whose
  morphisms are $k$-tuples of bijections.  The components of $\catname{X}$ are
  $X_n$, which are the groupoids of $k$-tuples of $n$-element sets.  The
  cardinality of $\Phi$ is
  \begin{eqnarray}
     | \Phi | & = & \sum_{n \in \mathbbm{N}} |X_n |z^n\\
     \nonumber & = & \sum_{n \in \mathbbm{N}} | 1 /\,\!\!/ S_n |^k z^n\\
     \nonumber & = & \sum_{n \in \mathbbm{N}} \left( \frac{1}{n!} \right)^k z^n
  \end{eqnarray}
  This cardinality cannot appear as that of any structure type, since the
  coefficients of $z^n$ in such power series must be of the form
  $\frac{a_n}{n!}$ for $a_n \in \mathbbm{N}$, which is not the case here
  unless $k = 1$.
\end{example}

\subsection{Stuff Types As Power Series}

Once we have the mechanics of groupoids, stuff types, and their
cardinalities, it becomes possible to extend the analogy with formal
power series which began structure types.  Although for the purposes
of quantum mechanics, we are mainly interested in the \textit{Hilbert
space} structure of formal power series, it is worth pointing out that
more of the algebraic and arithmetical properties of power series can
also be recovered.  In particular, we will next see that the
cardinality maps for stuff types and groupoids lets us find analogs of
the evaluation of a power series at a real number - and, by extension,
the composition of two power series.  These could not, in general, be
done with just structure types.

In section \ref{stuff-quantum-sec} we return to features specially
relevant to quantum mechanics.

\subsubsection{Evaluation of Stuff Types}

One useful fact about stuff types is that we can define a sensible
notion of evaluation, which brings us back to the fact that we had two
ways of looking at power series as functions: either by the evaluation
map $z \mapsto f ( z )$ or by the map picking out coefficients in the
power series expansion, $n \mapsto f_n$.  Structure types have given
us a good categorified way of looking at the latter, but with stuff
types we will be able to do the former as well.  This also lets us
talk sensibly about the composition of types, just as we might talk
about composition of functions.  Since this is one motivation for
passing from structure types to stuff types, let's consider an example
of an algebraic operation with power series which can't be extended to
the setting of structure types in a way which is compatible with the
correspondence between structure types and power series:

This example of a stuff type suggests the way we speak of evaluating stuff
types at groupoids - by ``colouring'' elements of a set with objects of the
groupoid.  That is:

\begin{definition}
  Given a groupoid $Z_0$ and a stuff type $\Phi : \catname{X}
  \rightarrow \FSN$, the {\textbf{evaluation of $\Phi$ at $Z_0$}}.  is
  the groupoid $\Phi ( Z_0 )$ of $\Phi$-stuffed, $Z_0$-coloured finite
  sets, whose objects are pairs $( \phi, z_0 ) \in \Phi_n \times
  Z_0^n$, where $\phi$ is a way of equipping an $n$-element set with
  $\Phi$-stuff, and $z_0$ is a map $f : n \rightarrowlim Z_0$
  equipping each element of $n$ with an object of $Z_0$.  The
  morphisms of $\Phi ( Z_0 )$ are bijections of the underlying
  $n$-element sets with strands labelled by morphisms in $Z_0$ between
  the objects labelling corresponding elements.
\end{definition}

Notice that there is an action of $S_n$ on the objects of $\Phi ( Z_0
)_n$ - the ways of putting $\Phi$-stuff on a set of size $n$, together
with a $Z_0$-colouring of this set - which comes immediately from the
action of $S_n$ on the underlying set.

In the case where $\Phi$ is a structure type $F$, the groupoid
$F(Z_0)_n$ is just a weak quotient $(F_n \times Z_0^n) /\,\!\!/ S_n$.
That is, its objects consist of pairs: $F$-structures which can be put
on an $n$-element set, together with colourings of an $n$-element set
by $Z_0$ objects.  These have an action of the permutation group on
the underlying set (two such objects are isomorphic by reindexing the
sets in in the $F$-structure and the colouring in the same way).  So
then the groupoid cardinality of $F(Z_0)$ is
\begin{equation}
| \sum_{n = 0}^\infty (F_n \times Z_0^n) /\,\!\!/ S_n | = \sum_{n=0}^\infty \frac{|F_n| |Z_0|^n}{n!} = |F|(|Z_0|)
\end{equation}
where the first ``sum'' is a categorical coproduct.  This formula is
consistent with the formula for the generating function of a structure
type which we've seen before in equation (\ref{eq:strtypecard}).  The
analogous fact is true for any stuff type, though recall that in that
case we will not explicitly refer to the action of $S_n$ on the
groupoid.  However, just as above, we have
\begin{equation}
| \Phi(Z_0) | = |\Phi|(|Z_0|)
\end{equation}

\begin{example}\label{ex:compcolour}
  Take a groupoid $\mathcal{G}$, and say $\catname{X}$ is the groupoid
  of $\mathcal{G}$-coloured finite sets, and consider this groupoid as
  the stuff type \[ \Phi : \catname{X} \rightarrow \FSN \] where
  $\Phi$ is the forgetful functor which takes a $\mathcal{G}$-coloured
  finite set to its underlying set.  An object in $X_n$ (the preimage
  of $n$ in $\catname{X}$) consists of $n$ objects from $\mathcal{G}$,
  and morphisms of $X_n$ are $n$-tuples of morphisms in $\mathcal{G}$
  composed with permutations of the $n$ elements.  Thus, $X_n \simeq
  G_0^n / S_n$, the weak quotient by the action we have
  described. Taking its cardinality, we find that
  \begin{eqnarray}\label{eq:compcolour}
     |\Phi| & = & \sum_{n \in \mathbbm{N}} |X_n |z^n\\
     \nonumber & = & \sum_{n \in \mathbbm{N}} | {\mathcal{G}^n}/\,\!\!/{S_n}
     | z^n\\
     \nonumber & = & \sum_{n \in \mathbbm{N}} \frac{|\mathcal{G}|^n}{n!} z^n\\
     \nonumber & = & \mathe^{| \mathcal{G} |z}
  \end{eqnarray}
  This follows since groupoid cardinalities are compatible with both
  powers and weak quotients.

  Both the parallel with the generating function $\mathe^{|
  \mathcal{G} |z}$ and the notion that $F$ is the type of
  ``$\mathcal{G}$-coloured finite sets'' suggests that it makes sense
  $F$ should be seen as $E^{\mathcal{G} \times Z}$ - that is,
  composition of the type $E^{\Box}$ (``finite sets of'') with
  $\mathcal{G} \times Z$ (the product type whose objects are
  one-element sets together with objects of $\mathcal{G}$ -
  i.e. $\mathcal{G}$-coloured one-element sets).
  
  Some interesting special cases of this example: when $\mathcal{G}$
  is a groupoid which is just a set $k$ with only identity morphisms,
  we have a {\textit{structure type}} of ``$k$-coloured finite sets'';
  $\mathcal{G}$ is a group $G$ seen as a one-object groupoid, we have
  a notion of $(1/\,\!\!/G)$-coloured sets.

  Furthermore, the calculation above gives us the cardinality of the
  groupoid $\catname{X}$ of $\mathcal{G}$-colored finite sets itself
  (i.e. not as a stuff type) to be simply $\mathe^{|\mathcal{G}|}$
  (since no powers of $z$ appear, but the calculation is otherwise the
  same).  So we can see this groupoid as the evaluation of $E^{\Box}$
  at $\mathcal{G}$.
\end{example}

\subsubsection{Composition of Stuff Types}

At last we can return to the question of how to compose stuff types.
We have seen how to evaluate stuff types at groupoids - given a
groupoid $Z_0$, evaluating the stuff type $\Phi$ at it gives
$\Phi(Z_0)$, the groupoid of ``$\Phi$-stuffed $Z_0$-coloured finite
sets''.  Since evaluating a stuff type at a groupoid (whose
cardinality is a real number) yields another groupoid (whose
cardinality is again a real number), we should be able to repeat this
process as many times as we like.  In principle, for instance, we
should be able to describe $\Psi \circ \Phi ( Z_0 )$ as
``$\Psi$-stuffed, $\Phi ( Z_0 )$-coloured finite sets'' - a set with
$G$-stuff, whose elements are labelled by finite sets with $F$-stuff
and elements labelled with objects of $Z_0$.

Since a stuff type $\Phi$ is itself a groupoid over finite sets, $\Phi
: \catname{X} \rightarrow \FSN$, and we have a way of evaluating a stuff type at
a groupoid, we get a notion of composition for stuff-types. We have
seen in section \ref{colouredsets} that there is a groupoid of
$Z_0$-coloured finite sets, whose morphisms are bijections of sets
with strands labelled by morphisms in $Z_0$.  We saw in
(\ref{ex:compcolour}) that its cardinality is $\mathe^{|Z_0|}$.

We should think of this as an illustration of the above where $F$ is
the structure type, in which we had the stuff type ``being a finite
set'' composed with the type $Z_0 \times Z$ - $Z_0$-coloured
one-element sets, but we can generalize this to any stuff types $F$
and $G$, to obtain stuff types $F \circ G = F ( G )$ and in such a way
that $|F \circ G| = |F| \circ |G|$.

\begin{figure}[h]
\begin{center}
\includegraphics{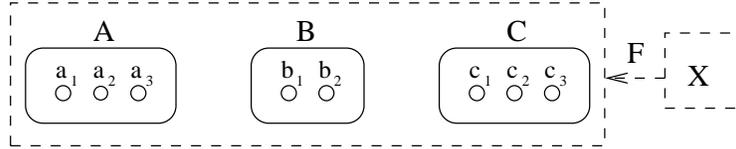}
\end{center}
\caption{\label{fig:compositetypeobject}An Object in a Composite Stuff Type}
\end{figure}

We can describe what we get here as a type which, evaluated at $Z_0$,
gives ``$F$-stuffed $G(Z_0)$-labelled finite sets''.  This has objects
(as shown in figure \ref{fig:compositetypeobject}) which consist of
finite sets equipped with $F$-stuff (say $F: \catname{X} \rightarrow S$).  The
elements of $F$ are labelled by objects of $G(Z_0)$: that is, the
labels themselves consist of finite sets labelled with $G$-stuff,
denoted by capital letters in the figure, whose elements are
themselves labelled in turn by objects from $Z_0$, denoted by
lower-case letters.

\begin{figure}[h]
\begin{center}
\includegraphics{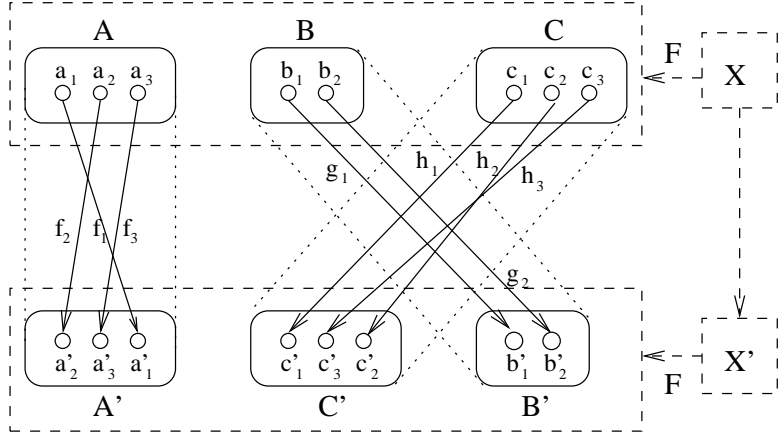}
\end{center}
\caption{\label{fig:compositetypemorphism}A Morphism in a Composite Stuff Type}
\end{figure}

Morphisms in the groupoid of $F$-stuffed $G(Z_0)$-coloured finite sets
(illustrated in figure \ref{fig:compositetypemorphism}) consist of
maps between the objects of $F$'s groupoid of stuff, which project
down to morphisms of the underlying $G(Z_0)$-coloured sets.  This has
its top level as a bijection of the underlying sets.  The strands of
the top-level bijection are labelled by morphisms of the groupoid of
$G(Z_0)$-coloured finite sets.  These include morphisms of the objects
associated to the $G$-stuffed $Z_0$-coloured sets (the dotted lines),
which are associated to bijections of the $Z_0$-coloured sets - with
strands labelled by morphisms of $Z_0$.

Clearly, we could in principle continue this sort of construction
recursively, to define the composite of any number of stuff types.

\subsection{Forgetful Functors: Properties, Structure, and Stuff\label{sec:forgetful}}

Here we take a detour through some ideas from category theory which
apply directly to stuff types.  The first is a classification of
functors by degrees of forgetfulness, which is the source of the term
``stuff'' in ``stuff type''.  In fact, we want to explain the terms
``structure'' and ``stuff'' which appear before the word ``type'' in
our terminology.  We began by trying to categorify a certain rigs,
replacing it by a 2-rig of functors - and now we look at a
classification of functors to see how we generalized this process.  In
particular, we are interested in functors which are, respectively:
faithful; full and faithful; or full and faithful and essentially
surjective:

\begin{definition}
  A functor $F : \mathcal{C} \rightarrow \mathcal{D}$ is
  {\textbf{{\textit{essentially surjective}}}} if the images of objects in
  $\mathcal{C}$ cover all objects of $\mathcal{D}$ in the sense that
  for any $d \in \opname{Ob} ( \mathcal{D} )$, there is some $c \in
  \mathcal{C}$ with $F(c) \cong d$.  It is {\textbf{{\textit{full}}}} if
  for any $c, c' \in \mathcal{C}$, the map between the sets of
  morphisms, $F : \opname{Hom}_{\mathcal{C}} ( c, c' ) \rightarrow
  \opname{Hom}_{\mathcal{D}} (F(c), F(c'))$, is surjective (in the
  set-theoretic sense). The functor $F$ is {\textit{{\textbf{faithful}}}}
  if for morphisms in its image, $F(f), F(f') \in
  \opname{Hom}_{\mathcal{D}} (F(c), F(c'))$, we have $F(f) = F(f')
  \Rightarrow f = f' \opname{in} \opname{Hom}_{\mathcal{C}} (c,c')$.
\end{definition}

Each of these can be seen as a form of surjectivity.  The notion of essential
surjectivity is a version of {\textit{onto}} for categories at the level of
objects. Fullness means that a functor is ``{\textit{onto}} for morphisms'';
faithfulness is ``{\textit{onto}} for equations between morphisms'' in the sense
that any equation between morphisms in $\opname{Hom}_{\mathcal{D}} ( F ( c ), F
( c' ) )$ comes from some equation between the preimage morphisms in
$\opname{Hom}_{\mathcal{C}} ( c, c' )$.

We have a classification of functors, then, by which of these it satisfies
(which we will attempt to explain next):

\begin{center}
  \begin{tabular}{|c|c|}
    \hline
    {\textbf{Functor}} & {\textbf{Forgets}}\\
    \hline
    essentially surjective, full, faithful & ``nothing''\\
    \hline
    full, faithful & ``properties''\\
    \hline
    faithful & ``structure''\\
    \hline
    all & ``stuff''\\
    \hline
  \end{tabular}
\end{center}

When we move to the setting of $n$-categories, we have not only objects and
morphisms, and the possibility that morphisms may be equal, but also
2-morphisms between morphisms (so that they may be 2-isomorphic, rather than
merely equal, or indeed might have non-iso 2-morphisms between them), and
3-morphisms between 2-morphisms, and so on.  Important to notice is that
essential surjectivity involves a weakening (a functor is essentially
surjective if it is surjective onto isomorphism classes, but not necessarily
objects).  This is because we should not distinguish between equivalent
categories, and since every category is equivalent to a skeletal category with
only one object in each equivalence class.  We will, in fact, want similarly
weakened versions of {\textit{full}}, {\textit{faithful}}, and their
higher-dimensional counterparts, when we use $n$-categories.

For now, though, we will explain the second column of this table, and see how
it applies to structure types.  The intuition begins with the commonplace fact
that a map between sets is an isomorphism if it is both injective and
surjective.  For functors between categories, essential surjectivity is the
natural analog of surjectivity, but \textit{full} functors are the natural
analogs of injective maps.  A map is injective if no two distinct objects have
the same image - that is, any equation of objects in the image comes from an
equation in the domain.  A functor is full if every morphism between objects
in the image comes from a morphism in the domain, so injectivity of set maps
is a special case, considering a set as a category with only identity
morphisms. In this trivial sort of category, every functor is faithful.

So a bijection between sets is a full, faithful, essentially surjective
functor, and if such a functor exists, we should treat the sets as ``the
same'' - the functor has lost no important information (of which there is very
little, for a set), which is reflected in the fact that it is invertible.  If
we have more general categories, with nontrivial morphisms, then in much the
same way, we have that a functor $F : \mathcal{C \rightarrow \text{D}}$, is an
equivalence of $\mathcal{C}$ and $\mathcal{D}$ (i.e. there is a functor $F^{-
1}$ which is an inverse to $F$ up to natural isomorphism) if $F$ is full,
faithful, and essentially surjective.  Similar results are true for
higher-dimensional categories ($n$-categories, for any $n$), while functors
which fail to have these properties are in various senses ``forgetful'' - not
being equivalences, they must forget information about the source category. 
This gives a ``periodic table'' of grades of ``stuff'' forgotten by functors
which fail to be onto for objects, morphisms, 2-morphisms (between morphisms),
3-morphisms (between 2-morphisms), and so on.  We restrict our attention to
the case $n = 1$, but note that the pattern we will develop continues for
higher $n$.

We see that these classes of functor can be used to talk about classes
of ``type'': structure types and stuff types are functors from
groupoids into the groupoid $\FSN$ - in fact, they are functors which
forget ``structure'' and ``stuff'' respectively.

\subsubsection{Examples of Forgetful Functors}

Some examples illustrate the classes of functors we have described.  So for
instance, if a functor $F : \mathcal{C} \rightarrow \mathcal{D}$ is not
essentially surjective, but is full and faithful, we have a subcategory of
$\mathcal{C}$ in $\mathcal{D}$, namely the image of $F$, which is equivalent
to $\mathcal{C}$, but which does not exhaust the isomorphism classes of
$\mathcal{D}$:

\begin{example}
  The functor \[ I : \catname{AbGrp} \rightarrow \catname{Grp} \]
  embeds the category of abelian groups and their homomorphisms into
  the category of all groups and homomorphisms.  This functor has
  ``forgotten a property'', namely the property of being Abelian.  The
  category $\catname{Grp}$ does not discriminate between objects with
  and without this property, whereas $\catname{AbGrp}$ is
  distinguished by the fact that it does.
\end{example}

If a functor fails to be essentially surjective \textit{and} fails to be full,
we can have not only the sort of situation above, but the target category can
have morphisms which do not correspond to those in the source:

\begin{example}
  The functor \[ I : \catname{OrdFinSet_0} \rightarrow \FSN \] takes
  any element in the groupoid $\catname{OrdFinSet_0}$ of totally
  ordered finite sets (whose morphisms are order-preserving
  bijections) into $\FSN$ - taking ordered sets and order-preserving
  maps to the underlying sets and set-maps.  This fails to be full
  since there can be bijections between the underlying sets of two
  ordered sets which fail to preserve order.  This functor ``forgets
  structure'' - namely, that structure which must be preserved by
  morphisms in $\catname{OrdFinSet_0}$, the total ordering on its
  objects.
\end{example}

Both of these examples are faithful functors, in the sense that each morphism
in the source category is sent to a distinct morphism in the target.  There
are functors which lack this property as well:

\begin{example}
  Consider the functor: \[ P_1 : \catname{Vect^2} \rightarrow
  \catname{Vect} \] The objects and morphisms in $\catname{Vect^2}$
  are ordered pairs of those in $\catname{Vect}$, namely vector spaces
  and linear transformations between them, while the functor $P_1$ is
  just projection onto the first component of these pairs.  Clearly,
  this is not faithful, since there are many pairs morphisms in
  $\catname{Vect^2}$ with the same first component.  What this functor
  has ``forgotten'' is, for each object, an entire vector space, and
  all the information about morphisms associated to these.  More than
  simply forgetting about properties shared by all objects, or
  structure which must be preserved by morphisms, we say this functor
  forgets ``stuff'' - parts of objects, in this case.
\end{example}

\begin{example}
  \label{st-ex:tuple}Our previous example of a completely forgetful
  functor involved the category of vector spaces, but similarly, there
  is an obvious class of stuff types associated to the groupoid
  $\FSN^n$, whose objects and morphisms are $n$-tuples of finite sets
  and bijections:
  \[ 
     P_j : \FSN^n \rightarrow \FSN
  \]
  where $P_j$ is the projection onto the $j^{\opname{th}}$ coordinate.
  The set in this entry is the ``underlying set'', and the ``stuff''
  being forgotten consists of all the other sets in the tuple.  This
  stuff type is a completely forgetful functor.
\end{example}

\section{Stuff Types And Quantum Mechanics} \label{stuff-quantum-sec}

Now we will return to our original motivation - seeing stuff types as
a categorification of states of the quantum harmonic oscillator.  We
have already seen how stuff types can take the role of formal power
series, at least up to the level of linear structure.  The space of
formal power series should be the Hilbert space of states of the
quantum oscillator: so composition of stuff types, which we have
explored, will not enter into this picture, though the linear
structure will.  If stuff types are categorified power series, and the
2-category $\catname{StuffTypes}$ is a categorified Fock space some of
the basic structure for a Hilbert space{\footnote{A treatment of the
construction of Hilbert spaces from structure types appears in
{\cite{maassen}}.}}.  One difficulty is that we do not have additive
inverses, and another is that we lack an inner product (which any
Hilbert space must have, and which we need to calculate probabilities
in a quantum system).

Interestingly, while the additive inverses are rather tricky to
define, and will have to wait until section \ref{mstuff} when we
define $M$-stuff types, the inner product does not even need to be
defined as a special construction, or imposed as extra structure: a
canonical one arises directly from the categorified framework as a (in
part because by describing stuff types as functors into $\Set$, we are
in effect choosing an ordered basis for this categorified equivalent
of a Hilbert space). What is more, a Hilbert space should have an
algebra of linear operators - endofunctions - which acts on it.  We
need to see what the equivalents of operators on stuff types are before
we can use them to categorify the oscillator.  We will see that these
are directly related to the inner product.

\subsection{Inner Product of Stuff Types} \label{stuffinnerprod}

The first feature of a Hilbert space we need to recover is its inner
product.

\subsubsection{Inner Product as Pullback}

When we were discussing the Hilbert space of states of the quantum
harmonic oscillator, we described the inner product on this Hilbert
space, which gave us, on the basis of pure states $z^k$, the form
$\inprod{ z^n, z^m } = \delta_{n, m} n!$.  We would like to give this
a combinatorial, or categorified, interpretation in terms of our
categorified states - stuff types.  The inner product for stuff types
(now thought of as categorified states) will turn out to satisfy the
property that $| \inprod{ \Phi, \Psi } | = \inprod{ | \Phi |, | \Psi |
}$ just as we saw for composition of types.

Previously we described the categorification of formal power series in
two steps - first replacing the complex numbers by some groupoid
$\catname{X}$, then replacing the natural numbers by $\FSN$.  We
follow the same two-step process to describe the categorified inner
product.

Now, in a quantum mechanical setting, the space of states of a system is a
space of $L^2$ functions over some configuration space.  In the case
of the harmonic oscillator, the configuration space is just a set of
energy levels, equivalent to $\mathbbm{N}$, and so the space of states
can be seen as $\ell^2$, the space of square-summable sequences (of
complex numbers).  

In the categorified space of states, a natural way to get an inner
product is to extend the definition $\inprod{ \psi, \phi } = \sum_{n
\in \mathbbm{N}} \psi_n \phi_n$ in the case of complex numbers to
become
\begin{equation}
  \inprod{ \psi,\phi } = \sum_{n \in \mathbbm{N}} \psi_n \times \phi_n
\end{equation} where the sum is now a coproduct in $\Set$ (i.e. a
disjoint union) and the multiplication is the categorical (Cartesian)
product.  This way of looking at the inner product gives a set, and
the equivalent condition to square-summability is the finiteness of
this set.  The cardinality operator gives us a notion of
square-summability in the decategorified setting: we can dispense with
this condition in the categorified setting, just as we can treat
structure types with infinite sets as coefficients.  In the case where
the set we get is finite, it should be clear that when we take
cardinalities, we just get that the cardinality $|\inprod{\psi,\phi}|$
is just $\inprod{|\psi|,|\phi|}$.  So this is naturally seen as an
inner product in the sense that when we take cardinalities, we get a
number which will be the inner product of the vectors in the space of
power series which are the cardinalities of the two types.  In the
finite case these are just polynomials: this is as yet not very
interesting since we are dealing with sets, whose cardinalities are
just integers.  We will address this when we pass to the case of a
groupoid over $\FSN$.

To see how this will work, note that, just as structure types and
stuff types themselves, we can treat the inner product as a ``bundle
over $\mathbbm{N}$'', with projection maps taking the individual
elements of this disjoint union of products down to $\mathbbm{N}$ by
the obvious projection taking an element of $\psi_n \times \phi_n$ to
$n$ (which is well defined, since all elements of the inner product
are of this form).  So if we think of the two states $\psi,
\phi$ as corresponding to two bundles $F : S \rightarrow \mathbbm{N}$ and $G :
S' \rightarrow \mathbbm{N}$ for $S, S' \in \opname{Set}$ we have, as
we've seen:
\begin{equation}
  \inprod{ F, G } = \sum_{n \in \mathbbm{N}} \psi_n \times
  \phi_n = \sum_{n \in \mathbbm{N}} F^{- 1} ( n ) \times G^{- 1} ( n )
\end{equation}
This can be understood in categorical terms as the fibered product of
the two bundles, also written as $X \times_{\mathbbm{N}} Y$ where the
projection maps $F, G$ of the bundles are understood.  Another way to
say this is to describe it as a {\textit{pullback}} of the two projections
$X \rightarrowlim^F \mathbbm{N} \leftarrowlim^G Y$, which is to say an
object which is initial among objects of the form $X
\leftarrowlim^{P_x} O \rightarrowlim^{P_y} Y$ making the square
\begin{equation}
 \xymatrix{
  O \ar[r]^-{P_x}\ar[d]_-{P_y} & X\ar[d]^-{F} \\
  \catname{Y} \ar[r]_-{G} & {\mathbbm{N}}
 }
\end{equation} commute.  (That is, given any other such object $O'$
with maps into $\catname{X}$ and $\catname{Y}$, there is a unique map from $O$ to $O'$
making the combined diagram commute.)

\subsubsection{The Categorified Case: Inner Product As Weak Pullback}
The next level of categorification of this description will give us a
definition of the inner product of two stuff types as a groupoid - namely the
pullback of the two functors from groupoids $\catname{X}$ and $\catname{Y}$ into
$\FSN$.   The inner product on this space is the same
one we described when talking about the harmonic oscillator (equation
\ref{innerproductuncat}).  To get this more fully categorified inner product -
now an inner product of stuff types - we should replace $\mathbbm{N}$,
the ``base space'' by $\FSN$, which is the free symmetric monoidal
category on one generator, just as $\mathbbm{N}$ is the free
commutative monoid on one generator.

So suppose we have two stuff types, namely functors $\Psi : X
\rightarrow \text{$\FSN$}$ and $\Phi : \catname{Y} \rightarrow \text{$\FSN$}$,
for some groupoids $X, \catname{Y} \in \catname{Gpd}$.  We want to do the
equivalent of taking the pullback of these two functors:
\[ \catname{X} \rightarrowlim^{\Psi} \FSN \leftarrowlim^{\Phi} \catname{Y} \] Since these
are not just functions between sets, but functors between categories,
the pullback is in the 2-category $\catname{Cat}$, or indeed in
$\catname{Gpd}$.  So defining a pullback becomes slightly more
complicated - in fact, we can and should weaken the requirement that
the pullback square

\begin{equation}
 \xymatrix{
   & \inprod{ \Psi, \Phi } \ar[rd]^-{P_X}\ar[ld]_-{P_Y} & \\
  \catname{Y} \ar[rd]_-{\Phi} & & X\ar[dl]^-{\Psi} \\
   & \FSN & \\
 }
\end{equation} commutes exactly, and allow it to commute only up to a
2-isomorphism between the two composite projections, so that what we
want is the {\textit{weak pullback}} (an example of a ``pseudo-limit''
in a 2-category; see \cite{street}) :

\begin{equation}\label{eq:weakpullback}
  \xymatrix{ 
     & \inprod{ \Psi, \Phi } \ar[rd]^-{P_X}\ar[ld]_-{P_Y} & \\ 
    \catname{Y} \ar[rd]_-{\Phi} & & X\ar[dl]^-{\Psi}\ar@{=>}[ll]^{\alpha}_{\sim} \\
     & \FSN & \\ 
  }
\end{equation}

This is like the fibrewise product over $\mathbbm{N}$ which we
described above, and we can also denote this by $\catname{X}
\times_{\FSN} \catname{Y}$, emphasizing the groupoids rather than the
functors.  Let's understand this weak pullback better by seeing what
this groupoid actually looks like internally, and then seeing that the
groupoid cardinality of this inner product of stuff types corresponds
to the inner product of two states $| \Psi |$ and $| \Phi|$.

\begin{definition} The groupoid $\inprod{ \Psi, \Phi } = \catname{X}
\times_{\FSN} \catname{Y}$ has objects which are pairs $(x,y) \in
\catname{X} \times \catname{Y}$ equipped with an isomorphism
$\alpha_{(x,y)} : \Psi(x)\tilde{\rightarrow}\Phi(y)$.  A morphism in
$\inprod{\Psi,\Phi}$ is a morphism in $\catname{X} \times
\catname{Y}$, say $(f,g) : (x,y) \rightarrow (x',y')$, such that
\begin{equation}\label{inprodnaturality}
\xymatrix{   
\Psi(x) \ar^-{\Psi(f)}[r] \ar_-{\alpha_{x,y}}[d] & \Psi(x') \ar^-{\alpha_{x',y'}}[d] \\   
\Phi(y) \ar_-{\Phi(g)}[r] & \Phi(y') }   
\end{equation} 
commutes.  That is, $\alpha_{x',y'} \circ \Psi(f) = \Phi(g) \circ \alpha_{x,y}$.
\end{definition}

The isomorphism $\alpha$ is from the definition of weak pullback; in
a strict pullback, $\alpha_{( x, y )}$ would always be the identity -
that is, we would require that $\Psi(x)=\Phi(y)$.

\begin{figure}[h]
\begin{center}
\includegraphics{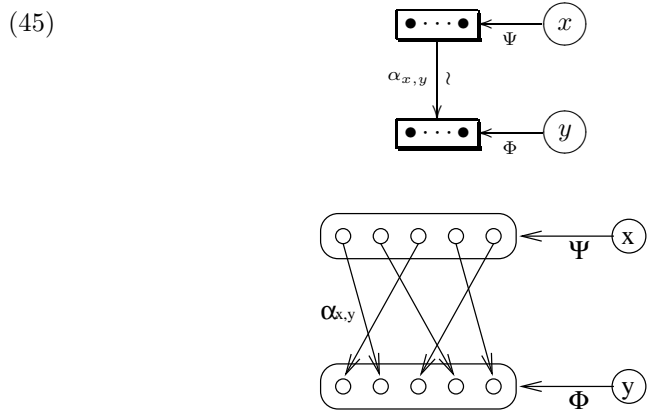}
\end{center}
\caption{\label{fig:pullbackobject}Object of an Inner Product (Pullback) Groupoid}
\end{figure}

Now, what does this all mean?  The isomorphism $\alpha$ is a bijection
of underlying sets - so an object of $\inprod{ \Psi, \Phi }$ is a pair
of objects $(x,y)$, ($\Psi$- and $\Phi$-stuff respectively on their
underlying sets), together with a bijection between the underlying
sets, as shown generall in (\ref{eq:pullbackobject}) and illustrated
in figure \ref{fig:pullbackobject}.

\begin{equation}\label{eq:pullbackobject}
 \xymatrix{
  *+[F-,]{\bullet \dots \bullet}\ar_-{\alpha_{x,y}}^{\wr}[d] & *++[o][F-]{x}\ar^-{\Psi}[l] \\
  *+[F-,]{\bullet \dots \bullet} & *++[o][F-]{y}\ar^-{\Phi}[l] \\
 }
\end{equation}

\begin{figure}[h]
\begin{center}
\includegraphics{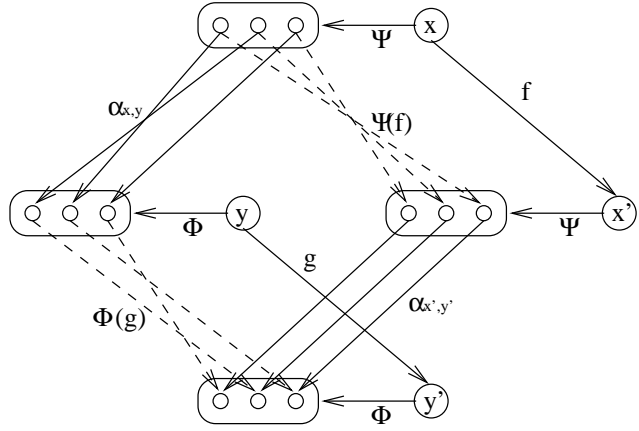}
\end{center}
\caption{\label{fig:pullbackmorphism}A Morphism in an Inner Product (Pullback) Groupoid}
\end{figure}

A morphism between two such objects includes morphisms of the form
(\ref{eq:typemorphism}) between the objects of $\Psi$ and $\Phi$; this
gives bijections between the underlying sets which must be compatible
with those which are part of the $\inprod{\Psi,\Phi}$ objects (of the
form (\ref{eq:pullbackobject})) themselves.  The general result is
illustrated in (\ref{eq:pullbackmorphism}) and an example appears in
figure \ref{fig:pullbackmorphism}.

\begin{equation}\label{eq:pullbackmorphism}
 \xymatrix{
   & *++[o][F-]{x}\ar_-{\Psi}[ld]\ar^-{f}_{\sim}[rr] & & *++[o][F-]{x'}\ar^-{\Psi}[ld]  \\
  *+[F-,]{\bullet \dots \bullet}\ar_-{\alpha_{x,y}}^{\wr}[dd]\ar^-{\Psi(f)}_{\sim}[rr] & & *+[F-,]{\bullet \dots \bullet}\ar^(.3){\alpha_{x',y'}}_(.3){\wr}[dd] & \\
   & *++[o][F-]{y}\ar_-{\Phi}[ld]\ar^(.3){g}_(.3){\sim}[rr] & & *++[o][F-]{y'}\ar^-{\Phi}[ld]  \\
  *+[F-,]{\bullet \dots \bullet}\ar^-{\Phi(g)}_{\sim}[rr] & & *+[F-,]{\bullet \dots \bullet} & \\
 }
\end{equation}

So the naturality square (\ref{inprodnaturality}) means that the inner
square of set bijections in (\ref{eq:pullbackmorphism}) commutes: in
other words, that we have given a way to identify all these underlying
sets.  We can see this in the figure as the fact that following
strands of the set bijections (of 3-element sets) around the square
reveals that there are exactly three complete squares of strands.

\subsubsection{Stuff Type and Fock Space Inner Products Related}

The key example of the inner product which makes clear the connection
with the inner product on Fock space developed in section
\ref{fockinnerprod} is the inner product of stuff types $Z^n$ and
$Z^m$, which also happen to be structure types.  The stuff type (in
fact, structure type) $Z^n$ is the type of {\textit{total orders on an
$n$-element set}}, or ``being a totally ordered $n$-element set''.  An
object of the inner product groupoid $\inprod{ Z^n, Z^m }$ is thus a
pair $( x, y )$ which are $n$- and $m$-element totally-ordered sets
respectively, equipped with an isomorphism between the underlying
sets.  We may think of these as the sets $n = \{ 0, \ldots, n - 1 \}$
and $m = \{ 0, \ldots, m - 1\}$, and the isomorphism is just any
bijection $\alpha$ between these.  If $n \neq m$, there are no such
objects since there are no such bijections.  If $n = m$, then there
are $n!$ such bijections, given by the permutations of $\{ 1,
\ldots, n \}$.  The morphisms of $\inprod{ Z^n, Z^m }$ are just the
identity morphisms on these objects, since the only morphisms in the
groupoid of totally ordered sets are order-preserving bijections -
that is, objects have no nontrivial automorphisms.  Thus, we get:
\begin{equation}
 \inprod{ Z^n, Z^m } = n! \delta_{n, m}
\end{equation} where $n!$ is the groupoid with $n!$ objects (one per
permutation of $n$) and only identity morphisms, and $\delta_{n, m}$
is analogous to the usual Kronecker delta, being a groupoid with no
objects if $n \neq m$ and with one object and one morphism if $n = m$.
Thus,
\begin{equation}
 | \inprod{ Z^n, Z^m } | = \inprod{ |Z^n |, |Z^m | }
\end{equation} so that the inner product we found on Fock space is
natural in this setting. This also illustrates the reason for the
factor $n!$ which shows up in the $n^{\opname{th}}$ term in the
expansion of a structure type, or the power series which is its
generating function.  This $n!$ is the cardinality of the group $S_n$,
the group of permutations of the underlying $n$-element set.

\subsection{Stuff Operators}

So far we have described stuff types, and implied that they, as
extensions of structure types, are a useful way of categorifying
functions - in the setting where these functions are seen as states of
a certain quantum system.  The inner product defined in section
\ref{stuffinnerprod} gives these some of the structure of a Hilbert
space, and also makes a connection to Feynman diagrams for energy
quanta of a harmonic oscillator.  We want to describe more of the
structure of the 2-rig of stuff types - in particular its linear
structure, demanding a definition for the equivalent of a linear
operator.  We call such a thing a {\textit{stuff operator}}.  There will
be a category of these, called $\catname{StuffOps}$, with
higher-dimensional algebraic structure similar to that of the algebra
of linear operators on a Hilbert space, with an action on the category
$\catname{StuffTypes}$.

We need to describe this action: a stuff operator $T$, given a stuff
type $\Psi:X\rightarrow\FSN$, ought to produce another stuff type
$T(\Psi):T(X)\rightarrow\FSN$, using natural category-theoretic
operations, in a way that reflects the fact that $T$ is the
categorified equivalent of a linear operator on a Hilbert space of
states.  One way to motivate our approach to constructing this is to
remember that for a Hilbert space $H$, a linear operator $T$ can be
thought of as an element of the tensor product $H \otimes H^{\ast}$ of
$H$ with its dual, and given a basis of $H$, $T$ can be represented as
a matrix (a two-index tensor).  To put this into the same framework as
the stuff type $\Psi$, recall that this can be seen as a vector in a
space of states{\footnote{Strictly speaking, the ``ground field'' here
is not $\mathbbm{C}$ but rather the setting for cardinalities of
groupoids, $\mathbbm{R}^+$, so we do not yet have a vector space.
When we introduce $M$-stuff and get cardinalities in something like
$\mathbbm{C}$, we will really have a vectorspace.}}, and the
equivalent of resolving it in a basis arises by taking the preimages
of elements of $\FSN$ as the components in the basis.  In this
setting, applying $T$ to a vector $\bldsym{v}$ (using the $H^{\ast}$
in the description in terms of $H \otimes H^{\ast}$) amounts to
applying covectors in any decomposition of $T$ to $\bldsym{v}$, which
(since $H \cong H^{\ast}$) amounts to taking the inner product of
$\bldsym{v}$ with a vector in $H$.  Since we have already seen a
natural definition of inner product for stuff types, we will use a
similar construction.

A stuff type $\Psi: \catname{X} \rightarrow \FSN$ can be variously
seen as the projection map of a groupoid-bundle on finite sets, and
also as a way of picking out an index for some component of a vector
in a categorified equivalent of a Hilbert space.  A stuff
{\textit{operator}} should have {\textit{two}} such maps, since it is
to correspond to a linear operator.  We should think of this as
``resolving $T$ in the same basis''.

The analogy with matrices suggests that $T$ acts on stuff types in the
same way as the inner product, just as matrices act on vectors
(resolved in a basis) by way of the inner product in each index.  So
indeed, $T \Psi$ is defined as a weak pullback:

\begin{definition} A \textbf{stuff operator} is a groupoid $T$ with
two projection maps into $\FSN$:
\begin{equation}
 \xymatrix{
   & T\ar_-{p_1}[ld]\ar^-{p_2}[rd] & \\
  \FSN & & \FSN\\
 }
\end{equation}
\end{definition}

\begin{figure}[h]
\begin{center}
\includegraphics{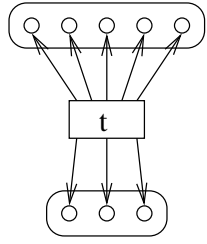}
\end{center}
\caption{\label{fig:stuffopobject}Object in the Groupoid of Stuff Operator $T$}
\end{figure}

We have seen the internal picture of objects of $\Psi$ in figure
\ref{fig:typeobject}, so we next look at a similar description for $T$.
 An object $t \in T$ is somewhat similar to an object of $\Psi$, but
this object has not just one underlying set, but two possibly distinct
ones, which we call $p_1(t)$ and $p_2(t)$, which in general need not
have the same cardinality.  We can think of this as two sets sharing a
common label, which we may think of as a ``process'' connecting
$p_1(t)$ with $p_2(t)$, where the label $t$ contains {\textit{stuff}}
associated to this transition, shown generally in
(\ref{eq:stuffopobject}) and illustrated in figure
\ref{fig:stuffopobject}.

\begin{equation}\label{eq:stuffopobject}
 \xymatrix{
  *+[F-,]{\bullet \dots \bullet} &
  *++[o][F-]{t}\ar^-{p_1}[l]\ar_-{p_2}[r] &
  *+[F-,]{\circ \dots \circ} \\
 }
\end{equation}

We are not yet at the point of recovering Feynman diagrams, but we can
already see an object which contains some sort of label for a process
connecting an $n$-element set of ``quanta'' (in the convention we have
already applied to the harmonic oscillator) to an $m$-element set of
``quanta''.  First we should see how these act like linear operators.

We can define an algebraic structure of stuff operators, when we
recall how the the categorical product and coproduct work for
groupoids $\catname{X}$ and $\catname{Y}$ in $\catname{Gpd}$.  The
coproduct $\catname{X} + \catname{Y}$ is just the direct sum of
groupoids, whose objects and morphisms are all objects and morphisms
of either $\catname{X}$ or $\catname{Y}$ (in a way which distinguishes
where they came from).  Their categorical product of has objects and
morphisms both consisting of ordered pairs of those from $\catname{X}$
and $\catname{Y}$.  We have seen that these operations are compatible
with groupoid cardinalities.

\begin{definition} Given two stuff operators $T$, $T'$ with projection
operators $p_1$, $p_2$ and $p'_1$, $p'_2$, respectively, the \textbf{sum} has
groupoid $T + T'$ whose projection functors $p_i + p'_i$ act like $p_i$
or $p'_i$ as appropriate.  The \textbf{product of $T$ by a groupoid} $\catname{G}$,
$\catname{G}T$ has objects is the product groupoid $\catname{G} \times
T$ (with projection operators acting on the $T$ component).
\end{definition}

These naturally have the properties that $(T + T')\Psi \cong T\Psi +
T'\Psi$ and $(\catname{G}T)\Psi \cong \catname{G}(T\Psi)$ in the sense
of the sum and ``scalar product'' of stuff types, and thus the
corresponding facts hold ``on the nose'' (i.e. as equations) for
cardinalities.  On the other hand, if we think of the stuff operators
as the categorified equivalent of infinite matrices (with projection
operators the equivalent of indexing), we can think of these as the
sum and scalar product for matrices, defining the linear structure of
the algebra of operators.  We leave the proof of this to the
interested reader.

More interesting is its internal multiplication, and the action on our
categorified Hilbert space, the 2-category $\catname{StuffTypes}$.

\begin{definition}\label{def:stuffopdefns} There is also an
\textbf{action of $T$ on a stuff type} $\Psi$ giving a
stuff type $T\Psi$ given by a weak pullback:
\begin{equation}
 \xymatrix{
   & & T\Psi\ar_-{P_T}[ld]\ar^-{P_X}[rd]\ar@/_2pc/_-{P_T;p_1}[ddll] & \\
   & T\ar^-{p_1}[ld]\ar_-{p_2}[rd]\ar@{<=}[rr]_{\alpha}^{\sim} & & X\ar^-{\Psi}[ld] \\
   \FSN & & \FSN &\\
 }
\end{equation}
and a \textbf{composite} of $T$ and $T'$ given similarly:
\begin{equation}
 \xymatrix{
   & & TT'\ar_-{P_T}[ld]\ar^-{P_T'}[rd]\ar@/_2pc/_-{P_T;p_1}[ddll]\ar@/^2pc/^-{P_T';p'_2}[ddrr] & & \\
   & T\ar^-{p_1}[ld]\ar_-{p_2}[rd]\ar@{<=}[rr]_{\alpha}^{\sim} & & T'\ar^-{p'_1}[ld]\ar_-{p'_2}[rd] & \\
   \FSN & & \FSN & & \FSN\\
 }
\end{equation} 
\end{definition}

Here, we note that $(TT')\Psi \cong T(T'\Psi)$, which can be seen by
considering the isomorphisms $\alpha$ of the weak pullbacks in the
diagrams for these two constructs.

Now, the groupoid $T \Psi$ resulting from the action of $T$ naturally
becomes a stuff type (a groupoid over $\FSN$) by composition of the
projections $P_T;p_2 = p_2 \circ P_T$, where $P_T$ is the projection
onto $T$ from $T \Psi$, which is the pullback of the functor $\Psi$
from $\catname{X}$ onto $\FSN$ along the projection $P$ from $T$ to
the same copy of $\FSN$.  To understand this construction better, we
should see what the objects and morphisms of the groupoid $T \Psi$ are
like internally.

\begin{figure}[h]
\begin{center}
\includegraphics{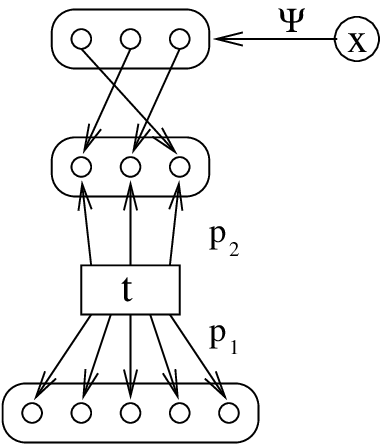}
\end{center}
\caption{\label{fig:stuffopresult}Object of Stuff Type $T\Psi$}
\end{figure}

The stuff type $T \Psi$ is a weak pullback of $T$ and $\Psi$, over the
copy of $\FSN$ which is the target of $\Psi : \catname{X} \rightarrow
\FSN$, and also of $p_2 : T \rightarrow \FSN$.  Its objects will be
pairs of objects $x \in \catname{X}$ and $t \in T$ together with
isomorphisms $\alpha_{x,t}:\Psi(x) \rightarrow p_2(t)$.  These are
isomorphisms of the underlying sets, so in particular they only exist
if these sets have the same cardinality.  Thus, an object of $T\Psi$
looks like an object of $T$ connected by a bijection of underlying
sets to an object of $\Psi$ (using the ``right-hand'' underlying set
of $T$).  The general form is shown in (\ref{eq:stuffopresult}) and an
example is illustrated in figure \ref{fig:stuffopresult}.

\begin{equation}\label{eq:stuffopresult}
 \xymatrix{
  *+[F-,]{\bullet \dots \bullet} & 
  *++[o][F-]{t}\ar^-{p_1}[l]\ar_-{p_2}[r] & 
  *+[F-,]{\circ \dots \circ} &   & \\
   &   &  *+[F-,]{\circ \dots \circ}\ar_-{\alpha_{x,t}}[u] &
  *++[o][F-]{x}\ar^-{\Psi}[l] \\
 }
\end{equation}

Since this is to be an object of $T\Psi$, we should see it as a stuff
type it its own right, over the copy of $\FSN$ mapped to under the
projection $p_1$ from $T$.  That is, everything in this picture - the
$t$ object, $x$ object, and the specific bijection $\alpha_{x,t}$
between the appropriate underlying sets - can be regarded as
``$T\Psi$-stuff'' attached to the underlying set $p_1(t)$.

\begin{figure}[h]
\begin{center}
\includegraphics{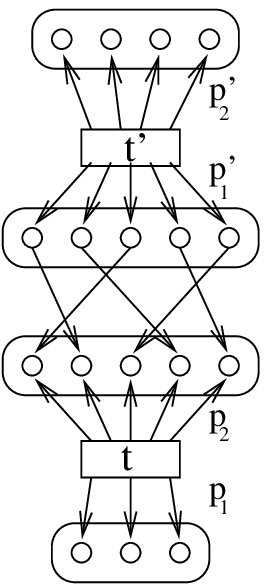}
\end{center}
\caption{\label{fig:stuffopcomposition}Object of Stuff Type $TT'$}
\end{figure}

Now as for the composite of $T$ and $T'$, similar reasoning holds
except that we have another stuff operator $T'$ in place of $\Psi$, so
the general form of an object of $TT'$ is as shown in
(\ref{eq:stuffopcomposition}), and an example is illustrated in figure
\ref{fig:stuffopcomposition}.

\begin{equation}\label{eq:stuffopcomposition}
 \xymatrix{
  *+[F-,]{\bullet \dots \bullet} & 
  *++[o][F-]{t}\ar^-{p_1}[l]\ar_-{p_2}[r] & 
  *+[F-,]{\circ \dots \circ} & & & \\
   &   &  *+[F-,]{\circ \dots \circ}\ar_-{\alpha_{t,t'}}[u] &
  *++[o][F-]{t'}\ar^-{p'_1}[l]\ar^-{p'_2}[r] & *+[F-,]{\diamond \dots \diamond}\\
 }
\end{equation}

The construction we have described for stuff operators is an example
of a ``span'': in particular, as morphisms from $\FSN$ to itself.  The
composition we have described above, as well as being analogous to
matrix multiplication of linear operators, satisfies the axioms for
composition of morphisms.  But in fact we have seen that these
operators also have an action on the category $\catname{StuffTypes}$,
derived from the fact that it is a category over $\FSN$.  In fact, we
can interpret stuff operators as endofunctors of
$\catname{StuffTypes}$, just as linear operators are endofunctions of
a vector space.  We describe this in more detail in appendix
\ref{sec:2cat}.

\subsection{Feynman Diagrams and Stuff Operators}\label{sec:feynman}

\begin{figure}[h]
\begin{center}
\includegraphics{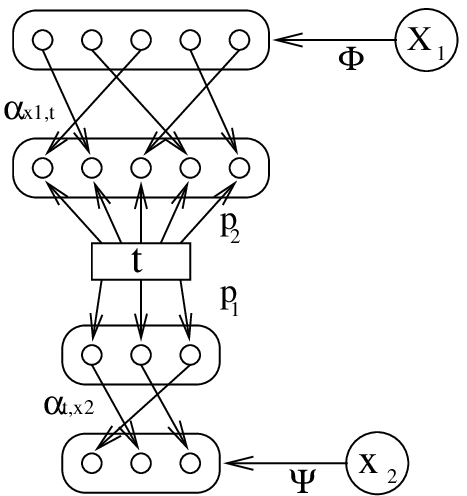}
\end{center}
\caption{\label{fig:stuffopinprodobject}Object of Groupoid $\inprod{\Psi,T\Phi}$}
\end{figure}

From a quantum mechanical point of view, we are often interested in
finding inner products such as $\inprod{ \Phi, T \Psi }$, and finding
these inner products can be done by means of Feynman diagrams.  That
is, in QM, the ``transition amplitude'' between states $\Phi$ and $T
\Psi$ is a sum of amplitudes associated to Feynman diagrams, each
showing one possible way of getting from state $\Phi$ to state $\Psi$
by process $T$.  We will show how this idea can be recovered in the
categorified setting, and in fact is given by exactly the algebraic
ideas we have defined.  In general, the groupoid $\inprod{ \Phi, T \Psi
}$ has objects as shown in figure \ref{fig:stuffopinprodobject}

We will see that we can think of these as Feynman diagrams, and
finding the sum over their amplitudes is exactly the process of taking
the groupoid cardinality of this inner product of stuff types - that
is, a sum over Feynman diagrams of some ``amplitudes''.  At the
moment, these amplitudes are all positive reals, rather than complex
numbers as is usual in quantum mechanics.  When we discuss $M$-stuff
types in section \ref{mstuff}, and in particular the case $M = U(1)$
in section \ref{mstuffu1}, we will see how this can be resolved by
introducing the quantum mechanical notion of {\textit{phase}}.

Since we are motivated here by the use of the algebra of stuff types
as a categorification of the Weyl algebra, we examine the stuff
operators $A$ and $A^{\ast}$, the annihilation and creation operators.
The annihilation operator $A$ can be realized in this form, with $T =
\FSN$ with two projection functors to $\FSN$, one of which is the
identity, the other of which is the functor whose action on objects is
to take a set $S$ and produce $S + \{ \star \}$:
\begin{equation}
  \FSN \longleftarrowlim^{\bldsym{1}}
  \FSN \longrightarrowlim^{+ \{ \star \}}
  \FSN
\end{equation}

To see that this reduces to our previous definition for $A$
(definition \ref{structladder}) on stuff types which happen to be structure
types, first recall that it said an $A F$-structure on a set $S$ is an
$F$-structure on the set $S \cup \{\star\}$.  If our stuff type $\Psi$
happens to be a structure type $F$ whose groupoid is just a set of
$F$-structured finite sets, then we have:
\begin{equation}
 \xymatrix{
   & & A S \ar_-{P_1}[ld]\ar^-{P_2}[rd]\ar@/_2pc/_-{AF=P_1}[ddll] & \\
   & \FSN\ar^-{\bldsym{1}}[ld]\ar_-{+\{\star\}}[rd]\ar@{<=>}[rr]_-{\alpha} & & S\ar^-{F}[ld] \\
   \FSN & & \FSN &\\
 }
\end{equation}

Tracing the map in the new type $A F$ from $A S$ to $\FSN$, we note
that we can pass through $S$ so that $A F = P_2;\alpha$, in which case
we see that since $\alpha$ must make the lower triangle commmute so
that $\alpha;+\{\star\}=F$ we get that the preimage of a given finite
set $S$ under $A F$ must correspond to the preimage of $S+\{\star\}$
under $F$.  So indeed, putting an $A F$-structure on $S$ amounts to
putting an $F$-structure on $S+\{\star\}$.

Similar reasoning shows that $A^{\ast}$, the adjoint of $A$, can be
realized in the same way, with groupoid $T^{\ast} = \FSN$ but with
the projections reversed:
\begin{equation}
  \FSN \longleftarrowlim^{+ \{ \star \}}
  \FSN \longrightarrowlim^{\bldsym{1}}
  \FSN
\end{equation}
Moreover, this acts like $A^{\ast}$, so that in the event that
$\Psi=F$ is a structure type, a $A^{\ast}F$-structure on a finite set
$S$ amounts to choosing an element of $S$ and putting an $F$ structure
on what remains.

\begin{figure}[h]
\begin{center}
\includegraphics{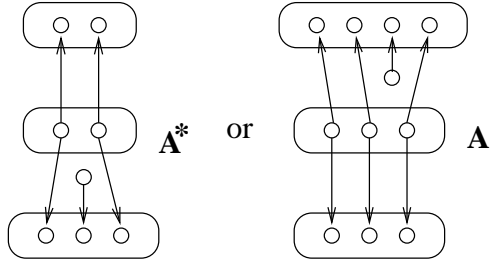}
\end{center}
\caption{\label{fig:fieldoperatorobject}Example Objects in the Categorified Field Operator}
\end{figure}

From the parallel with Fock space, one operator we should want to
define is the \textit{field operator} $A + A^{\ast}$.  As a stuff
operator, this behaves as one might expect: the groupoid in the stuff
operator is just the groupoid sum $T + T^{\ast}$ (i.e. two copies of
$\FSN$), and the projections just act as the projections on $T$ and
$T^{\ast}$ when applied to objects and morphisms from each of the two
copies.  So, the objects of the groupoid $\Phi$, following the pattern
we illustrated in (\ref{eq:stuffopobject}), look like either objects of
$A$ or of $A^{\ast}$, as shown in figure \ref{fig:fieldoperatorobject}.

\begin{figure}[h] 
\begin{center} 
\includegraphics{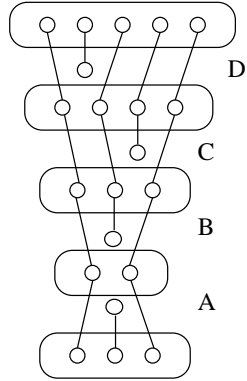} 
\end{center} 
\caption{\label{fig:phi-n}Example Object From $\Phi^4$} 
\end{figure}

In our categorified setting, this is written nearly the same way,
$\inprod{Z^n,\Phi^k Z^m}$, but we now have an interpretation of the
inner product as a groupoid over $\FSN$, obtained by a pullback
- and indeed, of the stuff operator $\Phi^k$ as a composite of stuff
operators, etc.  Now, this $\Phi^k$ has objects which are chains of
objects of the form in figure \ref{fig:fieldoperatorobject}, composed
as in figure \ref{fig:stuffopcomposition}.  We can draw these in
various ways (different drawing styles form objects in equivalent
categories), but for compactness, we will draw these in a style which
omits the ``internal'' bijections of the composite type, and also the
bijections from the $A$ and $A^{\ast}$ objects of $\Phi$ to $\FSN$.
Thus, each internal finite set would have previously been drawn three
times, with bijections between them.  This compact style is
illustrated in figure \ref{fig:phi-n}, which shows an object of
$\Phi^4$.  This one shows an object made from one annihilation and
three creation operators, in that order.  All other permutations of
four $\Phi$-objects are possible also.
 
\begin{figure}[h]  
\begin{center}  
\includegraphics{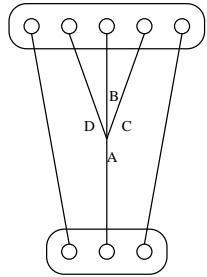}  
\end{center}  
\caption{\label{fig:field4oplabelled}Example Object From $\Phi^4$ (Alternate Style)}  
\end{figure}

For the sake of further compactness (and another drawing style
depicting an equivalent category), note that what we really have here
is an operation in which $k$ quanta of energy either appear or
disappear, and do so in a definite order.  We could really draw this
with just a single ``interaction'' vertex, incident with $k$ strands
(all other strands passing straight through the diagram, matching up a
quantum in the top with a quantum in the bottom).  These incidences
would have to be labelled with a total ordering, so that the object
shown in figure \ref{fig:phi-n} would be drawn as in figure
\ref{fig:field4oplabelled}.

There is an action of the permutation group $S_k$ on these objects,
changing the order in which we encounter objects from $A$ and
$A^{\ast}$ as we pass through $\Phi^k$.  That is, for every object we
get whose diagram has one order labelling the incidences, we will get
objects with all other possible orders exactly once through the action
of the permutation group $S_k$, also known as $k!$.  So if we want to
omit this labelling for clarity in the drawings, we can do so as long
as we remember that this means we are really drawing objects of the
weak quotient $\Phi^k /\,\!\!/ S_k$ (weak quotients are defined in
definition \ref{def:weakquotient}).  The objects of this weak quotient
are isomorphism classes of diagrams under permutations of labellings.
These permutations give the natural isomorphisms in the definition of
weak quotient by taking any labelled diagram to the same diagram with
permuted labels.  In this new category, the cardinality of the
groupoid is thus scaled by $\frac{1}{k!}$.  We also need to keep in
mind that an unlabelled diagram really stands in for possibly several
different inequivalent labellings of the incidences by distinct
orderings.

\begin{figure}[h]
\begin{center}
\includegraphics{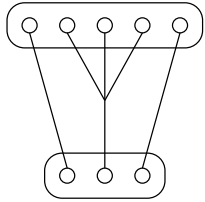}
\end{center}
\caption{\label{fig:field4op}Example Object in $\Phi^4 /\,\!\!/ S_4$}
\end{figure}

All this is really a notational convenience: really, to calculate
transition amplitude between states $\psi$ and $\phi$ for which we
have a description as stuff types $\Psi$ and $\Phi$, when we put the
system through some process $t$ which we describe as a stuff operator
$T$, we only need to find the groupoid cardinality
$|\inprod{\Psi,T\Phi}|$.  Simplifying diagrams and finding convenient
conventions for labelling them is really a way for getting a
calculational convenience out of diagrams like figure
\ref{fig:field4op}, which shows an object from the groupoid $\Phi^4
/\,\!\!/ S_4$.

\begin{figure}[h]
\begin{center}
\includegraphics{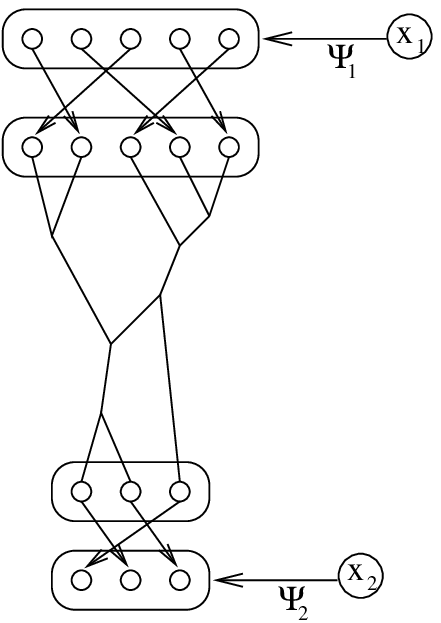}
\end{center}
\caption{\label{fig:feynmaninprod}An Object in $\inprod{\Psi_2,(\Phi^3 /\,\!\!/ S_3)^6 \Psi_1}$}
\end{figure}

Now, diagrams like this give vertices with $k$ incidences.  Taking
polynomials in the operators which give such diagrams gives operators
which can be interpreted in terms of diagrams having several such
vertices.  Such a diagram is shown in figure \ref{fig:feynmaninprod} -
note that here we continue the practice of omitting to draw the
internal finite sets in the composite stuff operator.

Given two stuff types $\Psi_1:\catname{X_1} \rightarrow \FSN$ and
$\Psi_2:\catname{X_2} \rightarrow \FSN$, we can take the inner product
$\inprod{\Psi_1,\Phi^{n}\Psi_2}$.  This applies this operator to
$\Psi_2$ to give compound objects involving objects of $\Phi^{n}$ and
of $\catname{X_2}$.  Taking the inner product with $\Psi_1$ gives
objects as illustrated in figure \ref{fig:feynmaninprod}.  The groupoid
cardinality of this inner product amounts to a sum over all such
diagrams, each with a weight related to the size of the symmetry group
of the diagram.

\begin{example}

We can use the above to show the categorical meaning of the usual
calculation of the expectation value of a power of the field operator.
In particular, suppose we want to calculate $\inprod{1,\frac{\phi^6}{6!} 1}$, the
vacuum expectation value of the 6th power of the (normalized) field operator.

To do this, we want to take a sum over objects which are equivalent to
ways of matching two empty sets with diagrams like figure
\ref{fig:field4op}, containing one vertex of valence 6.  We begin with
the case where incidences are labelled (as in figure
\ref{fig:field4oplabelled}).  Since the source and target sets are
empty, all edges must form loops touching the vertex at both ends.
The number of such diagrams is $(\binom{6}{2} \binom{4}{2}
\binom{2}{2})/3!=15$ (choosing the endpoints of three edges, without
order).  These give the objects of a groupoid one of which is shown in
figure \ref{fig:6vertex}.

\begin{figure}[h]
\begin{center}
\includegraphics{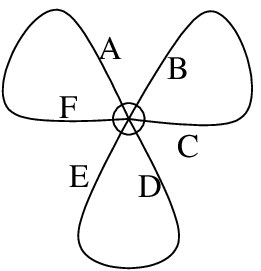}
\end{center}
\caption{\label{fig:6vertex}Example Object in $\inprod{1,\Phi^6 /\,\!\!/ S_6 1}$}
\end{figure}

The isomorphisms of this groupoid are given by permutations of the
labels.  Since this permutation group is $S_6$ with $720$ elements,
the groupoid cardinality should be $\frac{15}{720}=\frac{1}{48}$.  The
automorphism group of any such diagram is of size 48: there are 6 ways
to map the set of loops to themselves, each with the same orientation,
or reversed orientation.  Equivalently, we can think of the objects of
the groupoid as diagrams like this, but without labels.  In this case,
there is only one such diagram, with the automorphism group as just
described.  So we have $\inprod{1,\frac{\phi^6}{6!}1} =
|\inprod{1,\Phi^5 /\,\!\!/S_6 1}| = \frac{1}{48}$.

\end{example}

Now, the transition amplitude between two states in a quantum harmonic
oscillator which undergoes an interaction described by a given
operator in the Weyl algebra can be calculated, in part using a sum
over Feynman diagrams.  We have now seen how, in this categorified
setting, we can find a direct combinatorial interpretation for this
fact.

Unfortunately, some features of the diagrams used in quantum mechanics
are missing from this interpretation.  In particular, we do not have
any way to express the notion of ``phase'', or operators involving
propagators without interactions.  In quantum mechanics, states exist
in superpositions - linear combinations - having complex coefficients.
Non-interacting propagation in time involves the rotation of those
coefficients by a phase - that is, a unit complex number.  Thus, both
$\mathbbm{C}$ and the group $U(1)$ of phases are important.  We will
see in section \ref{mstuff} how to incorporate this into our
combinatorial picture.  As we shall see, this involves an explicit
decomposition of complex numbers into amplitudes and phases.

\section{$M$-Stuff Types And Quantum Mechanics} \label{mstuff}

We have described various entities under the heading \textit{types},
namely stuff, structure, and property types, and hinted at the
possibility that this sequence of classifications will continue as we
move into the setting of increasingly higher-dimensional categories.
For physical purposes, though, we are still missing some essential
properties which we would like a categorified version of quantum
mechanics to have.  In particular, all of our cardinalities, and hence
coefficients of our types, lie in \textit{rigs}, rather than rings -
they can be added and multiplied, but not subtracted and divided.  It
is possible to handle this in the abstract setting of structure types
by defining \textit{virtual structure types} as
equivalence classes of formal differences of structure types (and
similarly for stuff types) to make subtraction possible (see {\cite{BLL}}).  For quantum
mechanics, however, the coefficients of the power series which
represent states are complex numbers, and the physical significance of
\textit{phase} is of paramount importance, so we will still need
something more.  So for our purposes, it makes more sense to treat the
question in a different way.

\subsection{$M$-Stuff Types}

Now we will see an analog of a stuff type which can carry a phase - or
more generally, a weighting of some kind.  To reproduce some of the
features of quantum mechanics which don't appear in the picture of
stuff types as ``categorified states'', we should consider what is
missing.  First, states of a quantum system should form a Hilbert
space, and in particular a vector space.  Since we already have
something like an inner product, what is missing is the ability to
take linear combinations of states.  For this, they need to have a
notion of scalar product and of addition.  If we allow the
categorified states to carry a weight, this weight can play the role
of a scalar multiple, but these weights need to form a monoid, which
we think of as multiplicative.

This is the motivation for defining a notion of ``$M$-Stuff Types''
for some monoid $M$, and in particular $M = U(1)$, the group of
phases.  This is the case which is most interesting for quantum
mechanics.  We'll do this in general, since the construction does
not require $M = U(1)$.

\subsubsection{$M$-Sets}

Before we can talk about $M$-Stuff Types, we should start with a more basic
definition:

\begin{definition}
  If $M$ is a monoid, $\catname{MSet}$ is the category of ``sets over
  $M$''{\footnote{See appendix \ref{sec:slicecategory} for comments
  about such ``over categories''.  The usual definition applies here by
  treating $M$ as a set of elements, though we get some extra
  structure from the monoidal operation on $M$.}}, or ``$M$-sets''.
  Its objects are pairs $(S,f)$, for $S \in \Set$ and $f : S
  \rightarrow M$.  Morphisms between two $M$-sets $f_1 : S_1
  \rightarrow M$ and $f_2 : S_2 \rightarrow M$ are maps $g : S_1
  \rightarrow S_2$ in $\catname{Set}$ giving commuting triangles:
  \begin{equation}
   \xymatrix{
    S_1\ar_{f_1}[d]\ar^{g}[r] & S_2\ar^{f_2}[dl] \\
    M & \\
   }
  \end{equation}
\end{definition}

By abuse of notation, we will sometimes call the object just $S$ or
just $f$ if the meaning is clear by context.  A similar definition can
be made for $\catname{MFinSet}$ or $\MFSN$, where the sets $S$ lie in
$\catname{FinSet}$ or $\FSN$.

Note that for each set $S$ of cardinality $n$, the set of all $M$-sets
with ``overlying set'' $S$ is just $M^S$, equivalent to $M^n$.  The
morphisms which make this into an over category provide some extra
structure, however.

We also note that this definition is similar to that for a
$Z_0$-coloured set for a groupoid $Z_0$ - in fact, the image one
should have of an $M$-coloured set is just the same as figure
\ref{fig:colouredset}.  One difference is that in this case, the
picture we have for morphisms is different from that for
$Z_0$-coloured sets: for $M$-Sets, morphisms are just set maps which
are compatible with the labelling.  We could, of course, define a weak
over category of sets ``weakly over $M$'', as we did with
groupoid-coloured sets, for which strands of morphisms are also
labelled by elements of $M$, but as we shall see, this is not what we
want to do.  One result of this is that we lose some of the desirable
features of the category of sets, while retaining others.  For
example, we have the following:

\begin{theorem}
  $\catname{MSet}$ is a category with all colimits (in particular, it
  has coproducts).  
  
  \begin{proof}
    First, consider any diagram in $\catname{MSet}$, and take the
    underlying diagram in $\Set$.  Since $\Set$ is a cocomplete
    category, every diagram, and in particular this one, has a colimit
    $S$.  The diagram in $\catname{MSet}$ has a colimit provided we
    can construct a map from $S$ to $M$ which is compatible with the
    set-maps from the objects of the diagram in $\catname{MSet}$.  For
    any given element in $S$, every element taken to it by one of the
    maps in $\Set$ must have the same image in $M$ under the map for
    the corresponding object in $\catname{MSet}$, since if there is
    more than one, they must be taken to some common element by maps
    in the diagram.  Thus, we can consistently define $f(s)$ for any
    element in the colimit to be equal to the value for any preimage,
    and so all the maps in $\Set$ are compatible with the function
    into $M$, and the colimit in $\Set$ becomes a colimit in
    $\catname{MSet}$.
  \end{proof}
\end{theorem}

Coproducts in $\catname{MSet}$ can be interpreted as direct sums of
$M$-sets - and this makes it possible to define a cardinality for
$M$-sets.  We note that since $M$ has only one monoidal
operation, we could consider two interesting kinds of cardinality,
depending on whether we want this operation to look like addition or
multiplication.

For our purposes, it is better to think of $M$ as a
multiplicative monoid, since we will later want to take $M=U(1)$,
thought of as a subgroup of $\mathbbm{C}$.  So we would like to have a
notion of cardinality which gets along with multiplication in an
analogous way.  We should define a notion of cardinality
which reduces to set cardinality when we think of $M$ as
multiplicative.  Then we will find a ``tensor product'' compatible
with this notion of cardinality.

\begin{definition} The \textbf{cardinality} of an $M$-set $S
\rightarrowlim^{f} M$ is an element of $\mathbbm{N} \otimes M$ given
by
\begin{equation}
|S| = \sum_{s \in S} f(s)
\end{equation}
where the sum is taken in $\mathbbm{N}$.
\end{definition}

This cardinality operator is a kind of decategorification: it takes a
set $S$ labelled with values in $M$, and gives a formal sum of values
in $M$, each taken the number of times it appears in $S$.  Note that
this is again not compatible with the cartesian product in
$\catname{MSet}$, by the same argument as for the additive
cardinality.  Instead, we should take the following product:

\begin{definition} The \textbf{tensor product} of two $M$-sets $S
\rightarrowlim^{f} M$ and $S' \rightarrowlim^{f'} M$ is an $M$-set $S
\otimes S'$ has underlying set $\underline{S} \times \underline{S'}$
(the cartesian product of underlying sets in $\catname{Set}$).  The
map $(f \otimes f'): S \times S' \rightarrow M$ is given by $(f
\otimes f')(s,s') = f(s)\cdot f'(s')$.
\end{definition}

\begin{theorem} The tensor product of $M$-sets satisfies $|S \otimes
S'| = |S| \times |S'|$.  When $M$ is commutative, the tensor product
is symmetric.
  \begin{proof}
  \begin{eqnarray}
    |S \otimes S'| & = & \sum_{(s,s')\in \underline{S} \times \underline{S'}} f(s)f'(s') \\
   \nonumber & = & \bigl{(}\sum_{s \in \underline{S}} f(s) \bigr{)} \bigl{(} \sum_{s' \in \underline{S'}} f(s') \bigr{)}\\
   \nonumber & = & |S| \times |S'|
  \end{eqnarray} When $M$ is commutative, there is an obvious
  isomorphism between $S \times S'$ and $S' \times S$ taking $(s,s')$
  to $(s',s)$; the labelling is unchanged, since $f(s) f'(s') = f'(s')
  f(s)$.
  \end{proof}
\end{theorem}

These constructions for sets can be extended to groupoids, where
cardinalities start to look like complex numbers.

\subsubsection{$M$-Groupoids}

We would like to extend these results about $M$-sets to a notion of
$M$-groupoids and their cardinalities which is compatible with
cardinalities of $M$-sets and of ordinary groupoids in the suitable
special cases.  One way to see the correct approach to $M$-groupoids
is to take advantage our existing idea of $M$-sets and a connection we
already know between sets and groupoids.  This is the concept of a
groupoid-coloured set.  We define groupoid-coloured $M$-sets by
analogy with these.  Recall from \ref{colouredsets} that a set can be
seen as groupoid whose objects are the elements of the set, and with
only identity morphisms\footnote{Since we are already thinking of sets
as special kinds of categories here, this raises the question of what
happens if we categorify $M$.  We return to this in section
\ref{sec:categorifyingM}.}.  Then we have:

\begin{definition} 
  Given a groupoid $Z_0$, a \textbf{$Z_0$-coloured $M$-set} is an
  $M$-set $S$ equipped with a \textbf{colouring} map $c:S \rightarrow
  Z_0$.  \textbf{Maps of $Z_0$-coloured $M$-sets} are $M$-set
  bijections $\sigma: S \rightarrow S'$ together with, for each $x \in
  S$), a morphism $f_x \in \hom(c(x),c'(\sigma(x))$.  That is,
  \begin{equation}
   \xymatrix{
    S \ar^{\sigma}[r] \ar_c[d] & S' \ar^{c'}[d] \\
    Z_0 \ar@2{->}[r]_{\{f_x\}} & Z_0
   }
  \end{equation}
\end{definition}

This is essentially the same definition as appeared in
\ref{colouredsets}, but we note that now $\sigma$ is a bijection of
$M$-sets - that is, it is a set bijection \textit{which is compatible}
with the $M$-labelling.  But notice that groupoid-coloured $M$-sets
are just sets with two maps, one into a groupoid, and one into a
monoid:
\begin{equation}
\xymatrix{
   & S \ar_{c}[dl] \ar^{f}[dr] & \\
  Z_0 & & M \\
}
\end{equation}
Since the elements of this $Z_0$-coloured $M$-set are just elements of
$S$ labelled by both an object of $Z_0$ and an element of $M$, we
would like to be able to think of this as a set labelled by objects of
an ``$M$-groupoid'', which would look like objects of the groupoid
$Z_0$ labelled by element of $M$.  A consideration of what morphisms
of $Z_0$-coloured $M$-sets must be reveals how to define this:

\begin{definition} Given a monoid $M$, an \textbf{$M$-groupoid} is a
groupoid $\mathcal{G}$ with a functor $f$ from $\mathcal{G}$ into the
set $M$ regarded as a groupoid.  The \textbf{cardinality} of an
$M$-groupoid $\mathcal{G}_M$ is an element of $\MRp$, where
$\mathbbm{R}^+$ and $M$ are thought of as multiplicative monoids.  The
cardinality of $\mathcal{G}_M$ is given by the formal sum:
\begin{equation}
   |\mathcal{G}_M| = \sum_{[x] \in \underline{\mathcal{G}_M}} \frac{f(x)}{|\opname{Aut}(x)|}
\end{equation}
\end{definition}

\begin{remark} Note that $\MRp$ consists of all formal
$\mathbbm{R}^+$-linear sums of formal products $r \otimes m$ for $r
\in \mathbbm{R}^+$ and $m \in M$, subject to the distributive law
$(r+r')\otimes m = r \otimes m + r'\otimes m$.  It becomes a rig with
the obvious multiplication $(r\otimes m)(r' \otimes m') = (rr'
\otimes mm')$.

Since we are thinking of $M$ as a groupoid with only identity
morphisms, functoriality of $f$ means that for any $a$ and $b$ in
$\mathcal{G}$ and $g \in \opname{hom}(a,b)$, we have $f\circ g = f$.
That is, the following diagram commutes:
\begin{equation}
  \xymatrix{
    a \ar^{g}[r] \ar_{f}[d] & b \ar^{f}[dl] \\
    M & \\
  }
\end{equation} We see also that $f(x)$ is well defined for elements of
any given isomorphism, since any two objects with an isomorphism
between them will be sent under $f$ to the same element of $M$.  It
should be clear that in the case where the ``overlying'' groupoid of
an $M$-groupoid happens to be a set (i.e. groupoid with only trivial
morphisms), this reduces to the definition of an $M$-set and its
cardinality.  In the case where $M$ is the trivial groupoid, this
cardinality reduces to the usual groupoid cardinality.
\end{remark}

Given two $M$-groupoids, we define their product as with $M$-sets:

\begin{definition} The \textbf{tensor product} of two $M$-groupoids $X
\rightarrowlim^{f} M$ and $X' \rightarrowlim^{f'} M$ is an $M$-groupoid $X
\otimes X'$ which has underlying groupoid $\underline{X} \times \underline{X'}$
(the cartesian product of underlying groupoids in $\catname{Gpd}$).  The
map $(f \otimes f'): X \times X' \rightarrow M$ is given by $(f
\otimes f')(x,x') = f(x)\cdot f'(x')$.
\end{definition}

As with $M$-sets, this product gets along with $M$-groupoid
cardinalities.  The proof is essentially the same, except that
cardinalities involve factors of $|\opname{Aut}(x)|$.  This depends on
the fact that the automorphism group of an object $(x,x')$ in
$\underline{X} \times \underline{X'}$ is just the product of the
automorphism groups of $x$ and $x'$.

\subsubsection{$M$-Stuff Types and their Cardinalities}

We begin with a definition:

\begin{definition} An \textbf{$M$-stuff type} is an $M$-groupoid
$\catname{X} \rightarrowlim^{f} M$ equipped with a functor $\Psi :
\catname{X} \rightarrow \FSN$, where $\catname{X} \in \catname{Gpd}$.
\end{definition}

Typically, we will just think of $\catname{X}$ as an object of $\catname{MGpd}$
and blur the details, but this definition is what we always mean.  So
as with stuff types, we may think of $M$-stuff types as functors from
$M$-sets of ``$\Psi$-stuffed finite sets'' to their underlying finite
sets.  In the case where $\Psi$ is faithful we can say it is an
$M$-structure type.  Note that we are still thinking of $\catname{X}$
as lying over $\FSN$, not $\MFSN$ - we will return to this shortly.

Since stuff types (and $M$-stuff types) can be multiplied by
groupoids, whose cardinalities lie in $\mathbbm{R}^+$, this action by
$M$ gives another version of multiplication.  This will be
particularly interesting when we consider $M = U(1)$ in section
\ref{mstuffu1}, but first we should define the cardinality of an
$M$-stuff type:

\begin{definition}
  The \textbf{cardinality} of an $M$-stuff type $\Psi : \catname{X}
  \rightarrow \FSN$ is \begin{equation}
    |\Psi| = \sum_{n \in \underline{\FSN}} |\Psi_n|z^n
  \end{equation}
  where $|\Psi_n|$ is now the $M$-groupoid cardinality of the preimage of
  $n$ under $\Psi$.
\end{definition}

(As with a stuff type, this definition requires us to take $\FSN$ to
be skeletal to be well defined - or else to consider only the
\textit{essential preimage}.  We will do the former.)  This
cardinality is an element of $\MRz$: a formal power series in $z$
whose coefficients are formal combinations of pairs of groupoid
cardinalities and elements of $M$.

\begin{theorem} There are natural left and right actions of the monoid
$M$ on the $M$-stuff type $\Psi$.  If $M$ is abelian, these are the
same action, which satisfies
\begin{equation}
  |m\Psi|(z) = m|\Psi|(z)
\end{equation} 
 \begin{proof}
  We define the map $(m,\Psi) \mapsto m\Psi$, where $m\Psi :
  m\catname{X} \rightarrow \FSN$ acts as follows.  If $x\in
  \catname{X}$ is an object of $\catname{X}$ whose weight is $f(x)$,
  then the corresponding element $mx$ in $m\catname{X}$ has weight
  $m\cdot f(x)$. Then $m\Psi(mx)=\Psi(x)$. This is a left action on
  stuff types because it is a left action on $M$-groupoids together
  with a compatible map to $\FSN$.  The right action of $M$ is defined
  similarly.

  If $M$ is abelian, a left and right action are the same, and the
  result follows by direct calculation.
 \end{proof}
\end{theorem}

\subsubsection{$M$-Stuff Type Inner Product and $M$-Stuff Operators}

\begin{figure}[h]
\begin{center}
\includegraphics{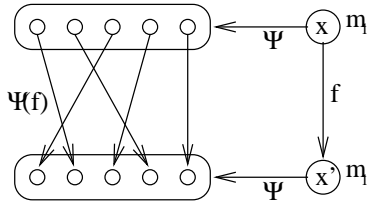}
\end{center}
\caption{\label{fig:mstuffmorphism}A Morphism in the Groupoid of an $M$-Stuff Type}
\end{figure}

Once we defined $M$-sets and hence $M$-groupoids, it was possible to
define $M$-stuff types simply by substituting these for groupoids in
the original definition of stuff types.  The only properties of $\FSN$
which were used in the original construction of stuff types and
operators was that it should be a groupoid: groupoids $\catname{X}$
with one or two functors into it were stuff types and operators
respectively.  Morphisms between the objects of a stuff type were
morphisms in $\catname{X}$ together with compatible bijections of sets
(recall figure \ref{fig:typemorphism}), but this depended only on the
fact that these were isomorphisms in the groupoid $\FSN$.  So, in the
same way, a morphism in the groupoid of an $M$-stuff type consists of a
morphism in its $M$-groupoid, together with a compatible bijection of
underlying sets, as illustrated in figure \ref{fig:mstuffmorphism}.
Notice that the objects $x$ and $x'$ are labelled by the same element,
$m_1 \in M$, since $f$ is an isomorphism in $\catname{X}$.

Now, the inner product of stuff types $\Psi:\catname{X} \rightarrow
\FSN$ and $\Phi:\catname{Y} \rightarrow \FSN$ was defined to be a
weak pullback, as described in (\ref{eq:weakpullback}).  The same
definition will apply if we let these be $M$-stuff types, allowing
$\catname{X}$ and $\catname{Y}$ to be $M$-groupoids.

So we have a weak pullback of $\Psi$ along $\Phi$, which gives a
groupoid $\inprod{\Psi,\Phi}$, and can define canonical projection
maps to $\catname{X}$ and $\catname{Y}$.  The groupoid is the
fibrewise product $\catname{X} \times_{\FSN} \catname{Y}$, where we
must use the tensor product of $M$-groupoids rather than the cartesian
product of groupoids to assign elements of $M$ to its objects.

\begin{definition} The given two $M$-stuff-types $\Psi: \catname{X}
\rightarrow \FSN$, and $\Phi: \catname{Y} \rightarrow \FSN$, the
$M$-groupoid $\inprod{ \Psi, \Phi } = \catname{X} \otimes_{\FSN}
\catname{Y}$ is the weak pullback of $\Psi$ and $\Phi$ over $\FSN$.
It has objects which are pairs $(x,y) \in \catname{X} \otimes
\catname{Y}$ equipped with an isomorphism $\alpha_{(x,y)} : \Psi(x)
\tilde{\rightarrow} \Phi(y)$.  A morphism in $\inprod{\Psi,\Phi}$ is a
morphism in $\catname{X} \otimes \catname{Y}$, say $(f,g) : (x,y)
\rightarrow (x',y')$, such that
\begin{equation}\label{eq:M-inprodnaturality}
\xymatrix{   
\Psi(x) \ar^-{\Psi(f)}[r] \ar_-{\alpha_{x,y}}[d] & \Psi(x') \ar^-{\alpha_{x',y'}}[d] \\   
\Phi(y) \ar_-{\Phi(g)}[r] & \Phi(y') }   
\end{equation} 
commutes.  That is, $\alpha_{x',y'} \circ \Psi(f) = \Phi(g) \circ \alpha_{x,y}$.
\end{definition}

\begin{figure}[h]
\begin{center}
\includegraphics{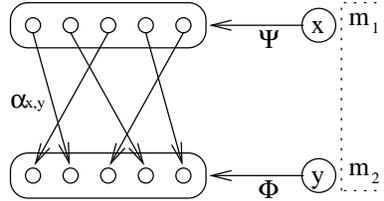}
\end{center}
\caption{\label{fig:mstuffinprodobject}An Object In the Inner Product of two $M$-Stuff Types}
\end{figure}

An object in the inner product groupoid looks like figure
\ref{fig:mstuffinprodobject}, where $m_i$ are elements of $M$.  This
figure is analogous to the previous inner product (figure
\ref{fig:pullbackobject}).  Note that, in contrast to the case in
figure \ref{fig:mstuffmorphism}, the objects $x \in \catname{X}$ and
$y \in \catname{Y}$ are not in the same groupoid, hence not related by
any morphism, so there is no requirement that $m_1$ and $m_2$ should
be equal.  The object illustrated is labelled be the element $m_1
\cdot m_2 \in M$ (as highlighted).

Similar changes apply to the other constructions defined using the
weak pullback, so that we have nearly identical categorical diagrams
definining morphisms in the inner product, as well as the action of an
$M$-stuff operator on an $M$-stuff type and the composition of two
$M$-stuff operators, as in definition \ref{def:stuffopdefns}.  The
sole change at this level is the replacement of groupoids with
$M$-groupoids, and thus in figures \ref{fig:stuffopresult} and
\ref{fig:stuffopcomposition}, we have labels in $M$ on every
groupoid object, preserved under every isomorphism.

\begin{definition} An \textbf{$M$-stuff operator} is an $M$-groupoid $T$ with
two functors from the underlying groupoid of $T$ into $\FSN$:
\begin{equation}
 \xymatrix{
   & T\ar_-{p_1}[ld]\ar^-{p_2}[rd] & \\
  \FSN & & \FSN\\
 }
\end{equation} It acts on an $M$-stuff type to give another $M$-stuff
type by weak pullback over one copy of $\FSN$.
\end{definition}

Just as with the former constructions, we have:

\begin{theorem} If $\Psi$ and $\Phi$ are two $M$-stuff types, then
$|\inprod{\Psi,\Phi}|=\inprod{|\Psi|,|\Phi|}$.
  \begin{proof}
    At the level of the underlying sets and groupoids, every product
    over a finite set in the (skeletal version of) $\FSN$ in the
    fibrewise product looks just the same as for regular stuff types.
    Each of these products is a product of $M$-groupoids, which are
    compatible with cardinality.  So the result holds.
  \end{proof}
\end{theorem}

\subsection{Quantum Mechanics:  $M=U(1)$}
\label{mstuffu1}

As remarked earlier, the notion of a {\textit{phase}} is crucial in
quantum mechanics.  Stuff types, and in particular stuff operators and
the inner product of stuff types, proved in the last section to have a
close connection to entities which resemble Feynman diagrams, but the
only notion of cardinality we had for these was groupoid cardinality,
which yields positive real values.  We would like to be able to do
more, since in quantum mechanics, these diagrams should have not a
real cardinality, but a complex amplitude, which has both a magnitude
and a phase.  This leads us to the idea of $U(1)$-stuff types,
since $U(1)$ is the group of phases, corresponding to the unit
circle in $\mathbbm{C}$.

\subsubsection{$U(1)$-Stuff Types}

From here on, we take $M=U(1)$ - an abelian monoid, and in fact an
abelian group - we know that $U(1)$-stuff types exist, and that they
have cardinalities in $\URz$, which has the obvious homomorphism into
$\Cz$.  As we have seen, the cardinalities of $M$-sets for an abelian
monoid $M$ lie in $M$, which we can think of ; cardinalities for
$M$-groupoids lie in $\MRp$; cardinalities for $M$-stuff types lie in
$\MRz$. When $M=U(1)$, this gives $\UR$, which has a homomorphism onto
$\mathbbm{C}$ 
\begin{equation}
  h : \UR \rightarrow \mathbbm{C}
\end{equation}

We should note that this description of $\mathbbm{C}$ in terms of
$\UR$ explicitly separates complex numbers into a magnitude and a
phase, and while it has a multiplication resembling that for
$\mathbbm{C}$, but it fails to capture the addition, which is formal.
However, the homomorphism $h$ just imposes the relations which define
complex addition.  The derived rig homomorphism $h: \URz \rightarrow
\Cz$ behaves similarly.  This homomorphism loses information, just as
the process of taking cardinalities does, so in fact, when $M=U(1)$,
we can define a new cardinality operator:

\begin{definition} If $\catname{X}$ is a $U(1)$-groupoid and
$\Psi:\catname{X} \rightarrow \FSN$ a $U(1)$-stuff type, the
\textbf{complex cardinality} of $\Psi$ is $h$ applied to the usual
$M$-stuff-type cardinality:
\begin{equation}
    |\Psi|_{\mathbbm{C}} = h \left( \sum_{S \in \mathbbm{N}} |\Psi_n| z^n \right) = \sum_{S \in \mathbbm{N}} h|\Psi_n| z^n
\end{equation} where $|\Psi_n|$ is the usual $M$-stuff cardinality,
$h$ is the above homomorphism, and addition is in $\mathbbm{C}$.
\end{definition}

The complex cardinality is a map which takes a $U(1)$-stuff type and
yields a power series in $\Cz$, namely Fock space.  When dealing with
$U(1)$-stuff types, we will write $|\Psi|_{\mathbbm{C}}$ as
$|\Psi|$, unless otherwise noted.

\begin{remark}Note that a type which consists of two states over
$U(1)$-sets of the same set cardinality but opposite phase will have a
cardinality in $\URz$ which contains a formal linear combination which
is in the \textit{kernel} of $h$.  This is the critical fact that when
we represent states in Fock space, there can be \textit{interference}
between states with opposite phases.  In particular, the amplitude for
a (categorified) state containing only those two objects will be zero.
\end{remark}

\subsubsection{Conjugation and The Inner Product}

There is a property of $U(1)$-stuff types which is not generally
shared by $M$-stuff types for arbitrary $M$, resulting from the fact
that it is an Abelian group.  This follows from the fact that there is
a nontrivial monoid isomorphism between $U(1)$ and itself, taking each
element of $U(1)$ to its multiplicative inverse.  There will be such
an isomorphism whenever $M$ is an Abelian group.  Viewing $U(1)$ as
the unit complex numbers, however, allows us to see this as complex
conjugation, which is how we will think of it.  Thus, there is an
operation special to $U(1)$-stuff types:

\begin{definition} If $\catname{X}$ is a groupoid with $U(1)$
labelling $f:\catname{X}\rightarrow U(1)$, its \textbf{conjugate}
groupoid is the $U(1)$ groupoid whose groupoid is labelling is
$\overline{f}$, given by $\overline{f}(x)=\overline{f(x)}$.  When we
write $\catname{X}$ for the $U(1)$-groupoid, we write the conjugate as
$\overline{\catname{X}}$.  If $\Psi:\catname{X} \rightarrow \FSN$ is a
$U(1)$-stuff type, its conjugate $\overline{\Psi}$ is the type which
acts like $\Psi$ on the objects of the underlying groupoid of
$\catname{X}$.
\end{definition}

This allows us to define a variant of the inner product which has the
conjugate-linearity of the usual complex inner product on Fock space.
To distinguish this from the (bilinear) inner product
$\inprod{\Psi,\Phi}$, and call it $\inprod{ \Psi | \Phi }$, also a
more familiar notation to physicists:

\begin{definition}
The \textbf{Fock space inner product} is given by 
\begin{equation}
\inprod{ \Psi | \Phi } = \inprod{ \overline{\Psi} , \Phi }
\end{equation}
\end{definition}

\begin{theorem} The Fock space inner product, for $U(1)$-stuff types
$\Psi$ and $\Phi$ satisfies $|\inprod{\Psi | \Phi}| =
\inprod{|\Psi| | |\Phi|}$, giving the usual conjugate-linear inner
product on $\Cz$.
\begin{proof}
\begin{eqnarray}
|\inprod{\Psi | \Phi}| & = & \sum_{n\in\mathbbm{N}} |\inprod{\Psi | \Phi}_n|\\
   \nonumber & = & \sum_{n\in\mathbbm{N}} |\overline{\Psi_n}| \cdot |\Phi_n|\\
   \nonumber & = & \inprod{\,|\Psi|\, |\, |\Phi|\,}
\end{eqnarray}
\end{proof}
\end{theorem}

So for $U(1)$-stuff types, we might want to define a new inner product
given in terms of the usual $M$-stuff type inner product as
$\inprod{\Psi|\Phi}$.  We will discuss briefly in section
\ref{sec:conjtimereversal} how to interpret this seemingly arbitrary
innovation.

It is worth noting here that the inner product between two states of a
$U(1)$-stuff type may be zero.  Of course, this can happen with stuff
types in any case: for instance, if $\Psi$ is the stuff type for which
there is one object in $\catname{X}$ over every even-cardinality set,
with automorphism group the same as that of the set, and no others;
and $\Phi$ is similar, but has objects over odd-cardinality sets.
These have cardinalities (generating functions) $\cosh(z)$ and
$\sinh(z)$ respectively.  They are examples of what we called
``property types'' in section \ref{sec:forgetful} - namely, these
types can be interpreted as the properties ``being an even set'' and
``being an odd set''.  These two stuff types are orthogonal in the
sense that their inner product is zero.  Here, the interpretation is
that there are \textit{no} sets having both properties, and so the stuff-type
inner product is the empty groupoid.

However, the situation with $U(1)$-types is more subtle: we may have a
nonempty inner product groupoid whose $U(1)$-groupoid cardinality
happens to be zero.  This arises because we may now have negative (and
indeed complex) contributions to the sum giving this cardinality.
This is related to the quantum mechanical phenomenon of ``destructive
interference'' between states.  In our formalism, this interference
occurs when we apply the homomorphism $h: \UR \rightarrow \mathbbm{C}$
and its derived variants.

We interpret the cardinality of the groupoid inner product as the
usual inner product in quantum-mechanics.  This is the amplitude for
finding our system in a given state $\Phi$ after setting it up in a
state $\Psi$, so this says this amplitude (and hence the probability)
is zero.

So the transition amplitudes between some of the ``pure''
(decategorified) states of which $\psi$ and $\phi$ are superpositions
may be nonzero, but the phases with which they appear may allow the
transition between $\psi$ and $\phi$ to have zero amplitude.  Thus,
introducing phases allows destructive interference which makes
otherwise feasible transitions impossible.  We will see in the next
section that this issue of phase is closely related to the concept of
\textit{time evolution} in quantum mechanics, and the
\textit{propagator}.

In particular, in the harmonic oscillator, the phase of a state
changes over time, with a frequency proportional to the energy of that
state.  This is the effect of the \textit{free propagator} for a
system, and it is an operator.  So we must describe $U(1)$-stuff
operators next, and this propagator in particular.

\subsubsection{$U(1)$-Stuff Operators and Time}

We have already noted that $M$-stuff operators act on $M$-stuff types
just like ordinary stuff operators acting on stuff types, except that
the groupoids are now replaced by $M$-groupoids.  The groupoid $T\Psi$
for an $M$-stuff operator $T$ and type $\Psi$ then consists of pairs
of objects $t \in T$ and $x \in \catname{X}$ together with a bijection
of their underlying sets.  As with the product groupoid, this object
is labelled by an element in $M$ given by the product of the labels on
$t$ and $x$.  This is well defined when $M$ is Abelian, as in the case
when $M = U(1)$.

\begin{figure}[h]
\begin{center}
\includegraphics{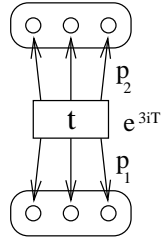}
\end{center}
\caption{\label{fig:timeevol}An Object in The $U(1)$-Stuff Operator $E_T$}
\end{figure}

One class of $U(1)$-stuff operators which is particularly relevant to
quantum mechanics is that of the \textit{time evolution} operators.
These are operators which, when applied to an $M$-stuff type $\Psi$,
produce an $M$-stuff type $\Psi'$ for which the $M$-labels on the
elements of the underlying sets have labels multiplied by a fixed
phase in $U(1)$.  An example of an object in such an operator,
designated $\theta$, is shown in figure \ref{fig:timeevol}.  Here, we
are showing the operator $E_T$, ``time evolution by $T$''.  We show an
object which will evolve a state with three quanta of energy.  Here,
$\mathe^{iT}$ is the change of phase corresponding to time evolution of a
one-energy-quantum state by $T$.  A state
with three energy quanta changes phase by $\mathe^{3iT}$.

\begin{figure}[h]
\begin{center}
\includegraphics{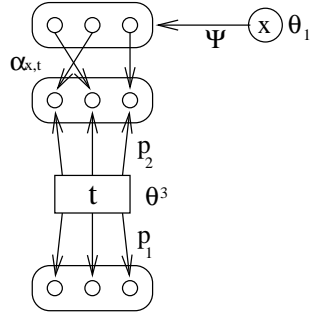}
\end{center}
\caption{\label{fig:timeevolresult}An Object in The $U(1)$-Stuff Type $E_T \Psi$}
\end{figure}

Any object in $E_T$ has two projections - each to an underlying set of
the same size.  This corresponds to the fact that in unperturbed time
evolution no interactions are occurring which would change the energy
level of the system.  An object $t$ in the groupoid of $E_T$, lying
over sets with $k$ quanta of energy.  There will be just one such
object in $E_T$ for each finite set.  It is labelled by the phase by
which a state with $k$ quanta will change in time $T$.  The operator
$E_T$ acts on any $U(1)$-stuff type (categorified state) to give state
to which this has evolved after a time $T$.  An object of the
resulting stuff type is shown in figure \ref{fig:timeevolresult}.

The object in the groupoid of $E_T \Psi$ is the entire ensemble
associated to the finite set $p_2(t)$.  It includes the object $t \in
E_T$ itself, as well as $x \in \catname{X}$, their underlying sets and
the bijection $\alpha_{x,t}$ between them, and also their associated
labels $\mathe^{3i \theta}$ and $\mathe^{i \theta_1}$.  This is an
object in a product $U(1)$-groupoid: $U(1)\catname{Gpd}$ is a weak
2-category with (weak) products, of which this is an example.  This
object in the product groupoid $E_T \times \catname{X}$ is labelled by
the product of the labels on $t$ and $x$, namely $\mathe^{i \theta_1}
\mathe^{3i \theta}$.

Suppose the $U(1)$-stuff type in question happens to just be $Z^k$ - the
categorified state with just $k$ quanta of energy and no phase angle,
or the property type ``being an $k$-element finite set labelled by $1
\in U(1)$''.  Then we get $E_T(Z^k) \cong (\mathe^{i \theta} Z)^k$, and the same
fact holds as an equation for the complex cardinalities.  So the
$U(1)$-cardinality of a $k$-element $U(1)$-set changes by $\mathe^{i k
\theta}$ in time $T$, since each quantum picks up a phase rotation of
$\mathe^{i \theta}$ in time $T$.  In particular, the phase of an
object in a categorified state changes with a frequency proportional
to its energy (the size of the underlying set).

Choosing time units so that $\theta=T$, we get a
phase change of $\mathe^{iTk}$ on a state of energy $k$.  We can write
this as $E_T=\mathe^{iTN}$, where $N$ is the number operator.  To
prove this equality at the categorified level (using a
categorification of the exponential such as we have already discussed)
would require a fully categorified version of the complex numbers.
However, for now we can observe that this will be true at the level of
cardinalities, and take it as a definition.  This arises physically
from the Hamiltonian formulation of quantum mechanics.  We will not
enter into this in detail, but note that the free Hamiltonian is just
$H_0 = N$, the number operator in the exponent of the propagator,
which measures the energy of a state.

Since $\mathe^{iTN}$ has exactly one object for each cardinality, the
product groupoid $\mathe^{iTN} \Psi$ is equivalent to the groupoid
whose objects are the same as those of $\catname{X}$, but whose
$U(1)$-labellings have been multiplied by phases $\mathe^{iTk}$, for
an object with underlying set $\bldsym{k}$.  This is the groupoid of
the state $\Psi$, time-evolved by $T$.  Note that time-evolution by
$-T$ will be given by a similar operator $E_{-T} = \mathe^{-iTN}$: all
the object-labels are the inverses of those of $E_T$.  We will return
to this point in section \ref{sec:conjtimereversal}.

In any case, using these propagators, and the $U(1)$ version of inner
product groupoids whose objects resemble Feynman diagram for
interactions, we recover a combinatorial interpretation for many of
the standard features of the quantum mechanics of the harmonic
oscillator.

\subsubsection{Feynman Diagrams And Perturbation}

Having found a categorification of the quantum harmonic oscillator, we
know that transition amplitudes such as
$\inprod{\psi,p(\frac{\phi^n}{n!})\psi '}$ as a sum over Feynman
diagrams.  In physically realistic settings, this sort of amplitude
often arises when we consider the time-evolution of an oscillator
which is perturbed.  The \textit{free} oscillator evolves in time
according to the operator $E_T=\mathe^{-iTN}$ described above.  The
\textit{perturbed} oscillator, on the other hand, describes a
situation where the energy of the oscillator is modified by another
term: it only approximately matches the description of the free
oscillator we have been using.  Physically, this represents a
potential in which the oscillator is moving.  This means that the
energy of the oscillator is changed by the addition of an extra term,
$V$, which is some function of position, which we think of as a
\textit{potential energy}, in addition to the energy in the oscillator
proper.

In this case, time evolution can be calculated using the new energy,
$H = H_0 + V$.  If $V$ is a function of position, then since the
position is proportional to $\opname{a} + \opname{a}^{\ast} = \phi$,
we have $V = f(\phi)$, for some function $f$.  We will consider the
case where $f=p$ is some polynomial (though naturally any analytic
function can be approximated this way to some degree, so we can obtain
successive approximations by taking a sequence of $p_k$ converging to
$f$).  Since the energy for the free oscillator is already quadratic
in position, we assume that $f$ has minimum degree at least 3.  In
this case, at the decategorified level, the amplitudes for time
evolution by time $t$ associated to the Hamiltonian $H = N + V$, are:
\begin{equation}\label{eq:perturbfeynman}
\inprod{Z^k|\mathe^{-iTH}Z^l} = \sum_{n=0}^{\infty} \int_{0 \leq t_1 \leq \dots \leq t_n \leq T} \inprod{Z^k|\mathe^{-i(T-t_n)N}V\mathe^{-i(t_n-t_{n-1})N}V \dots V\mathe^{-it_1N} Z^l}
\end{equation}
To avoid considering questions of convergence, we think of this purely
as a statement about power series in $T$.  It would take us too far
afield to derive this standard quantum-mechanical fact in full detail,
though background can be found in \cite{bjordrel}, and one derivation
of this equation in our setting can be found in \cite{catquan}.
However, we will point out here that it follows from the fact that the
evolution of a state is governed by the Schr\"odinger equation, which
amounts to:
\begin{equation}
\partial_t \psi = -i (\mathe^{itH_0}V\mathe^{-itH_0})
\end{equation}
Integrating this equation over time, we get
\begin{equation}
\psi(t) = -i \int_0^t (\mathe^{itH_0}V\mathe^{-itH_0})\psi(t_0) \mathd t_0 + \psi(0)
\end{equation}
and by repeated substitution of this expression for $\psi(t)$ into the
integral, we get the sum of integrals which appear in the expression
above.  Taking the inner product with this operator, we finally get
the whole expression.

Ideally, we would like to derive this equation entirely at the
categorified level.  However this would require a more complete
understanding of the categorified version of the complex numbers than
we have constructed here.  To recover time evolution by a phase from
an expression of the form $\mathe^{itH_0}$, we would need to see that
the result is indeed a phase in $U(1)$.

However, knowing that the equation holds at the decategorified level
allows us to give a simple interpretation for the formula.

\begin{theorem} The transition amplitude
$\inprod{z^k|\mathe^{-iTH}z^l}$ for the perturbed harmonic oscillator
with potential $V=f(\phi)$ is given by a sum over all Feynman diagrams
given as composites of those associated with $V$, from a state with
$k$ quanta to one with $l$ quanta.  The sum is of an integral over all
labellings of the edges of the diagrams such that the total phase
along all paths is $\mathe^{-iT}$.
\begin{proof}
This transition amplitude is the $U(1)$-groupoid cardinality of the
inner product $\inprod{Z^k|E^{iTH}Z^l}$, and given by
\ref{eq:perturbfeynman}.  Consider the operator in that equation,
\[
O=\mathe^{-i(T-t_n)N}V\mathe^{-i(t_n-t_{n-1})N}V \dots V\mathe^{-it_1N}
\]
We know that the terms $\mathe^{-i(t_i-t_{i-1})N}$ are just free
propagators, which contribute a phase of $\mathe^{-i()}$ for each
quantum of energy.  We can think of these operators as having objects
given by any number of ``strands'', one for each quantum, and each
strand labelled by a phase $\mathe^{-i(t_i-t_{i-1})}$, the total
number giving the total phase change associated to that energy.

Now consider the $U(1)$-stuff operators $V$.  Each of these operators
has a groupoid whose objects are naturally identified with Feynman
diagrams of the sort associated with $V$.  These do not affect phases.

Composing the operators together, we get all possible composites of
Feynman diagrams of the type associated to $V$, connected by diagrams
whose effect is to label strands by phases associated to the time
intervals between interactions.  To find the total phase associated to
such a composite, we multiply all phases.  This is clearly equivalent
to multiplying the phases on any labelled edges which are joined by
the composition to get a phase on the resulting edge, then multiplying
the product of all edges thus produced.

The transition amplitude we want to recover is:
\[
\sum_{n=0}^{\infty} \int_{0 \leq t_1 \leq \dots \leq t_n \leq T} \inprod{Z^k|O Z^l}
\]
The sum taken over all $n$ simply means that we are taking a groupoid
containing all possible $n$-fold composites of this form.  The
integral over all $n$-part partitions of the interval $[0,T]$ means
each such diagram contributes a phase found by integrating over all
possible ways of dividing the interval into free and interaction
parts.  This contribution is weighted by the size of the symmetry
group, since the inner product inside the integral is just a
$U(1)$-groupoid cardinality.

This proves the statement.
\end{proof}
\end{theorem}

This recovers the usual Feynman rules for calculating transition
amplitudes in the oscillator.

\section{Conclusions}

Here we summarize what we have shown, and suggest directions in which
this work could be carried further.

\subsection{Categorified Quantum Mechanics}

We began by describing the quantum harmonic oscillator and the Weyl
algebra, the algebra of linear operators on its space of states which
correspond to observables and interactions of the oscillator system.
We saw how this could be related - by the Fock representation of the
Weyl algebra - to formal power series with complex coefficients with
exponents counting quanta of energy.

Our aim at the outset was to categorify this aspect of quantum
mechanics.  Categorification of concepts such as ``group'',
``vector-space'', and indeed ``category'' itself have proved
interesting within mathematics, and the resulting 2-groups,
2-vector-spaces, and 2-categories arise naturally in surprising ways.
The idea here was that categorification could be applied in a
physically relevant setting, and could reveal something useful about
the mathematical structures involved.  Here, we began with the Hilbert
space of states of a quantum mechanical system, and the relevant
algebra of operators acting on it.  We have produced
category-theoretic equivalents of these: the 2-categories of stuff
types and of stuff operators can be seen as a categorified Hilbert
space and a categorified algebra.

These are connected to the original setting by concepts of
decategorification which go by the name ``cardinality''.  We have
shown that when we take the cardinalities of all our entities
involving ``stuff'', we recover much of the structure of the Weyl
algebra.  By introducing the idea of $M$-sets, and attendant ideas of
entities labelled with ``phases'' from some monoid, we have improved
this resemblance to quantum mechanics.

Stuff types - groupoids over $\FSN$ - have creation and annihilation
operators which give a purely combinatorial construction which
categorified many features of the Weyl algebra.  They also have a
natural inner product which, in conjunction with these creation and
annihilation operators, allows us to interpret transition amplitudes
as sums over Feynman diagrams.

However, what these categorify is not Fock space, since it only has
scalar multiplication over $\mathbbm{R}^+$, rather than $\mathbbm{C}$,
and cardinalities in $\mathbbm{R}^+[\![z]\!]$, rather than $\Cz$.  Our
$U(1)$-stuff types are a better categorification of Fock space, and
these have cardinalities in $\URz$, which we can map to $\Cz$.  This
map $h$ is not one-to-one, and this fact is responsible for the
phenomenon of interference of states with different phases.

The problem of categorifying quantum mechanics is much more general
than the simple case of a harmonic oscillator we have discussed.
Another approach to bringing category theory to quantum mechanics is
{\cite{coecke}}.  That paper provides good description of a simple
``picture calculus'' for quite general quantum mechanics which uses a
background of category theory.  This is not a categorification in our
sense, but together with some of the structure described in appendix
\ref{sec:slicecategory} may suggest a broader framework for dealing
with the question.

Although we have confined ourselves to the harmonic oscillator in this
paper, we can suggest various directions in which these ideas could be
taken further.  One is to look at the inner product through a more
category-theoretic lens.

\subsubsection{Conjugate-Linearity and the Inner Product}\label{sec:conjtimereversal}

Recall that the inner product for $M$-stuff types had to be modified
somewhat in order to agree with the usual inner product on the Hilbert
space $\Cz$ in the case when $M=U(1)$.  The nontrivial isomorphism of
$U(1)$ with itself provides a notion of complex conjugation.  But how
should we interpret the inner product for $U(1)$-stuff types?  

In fact, it makes more sense when we adopt the interpretation of the
inner product $\inprod{\phi|\psi}$ as pairing a \textit{state vector}
with a \textit{costate covector}.  So the (conjugate-linear) inner
product $\inprod{\psi|T\phi}$ gives the amplitude to find a system
\textit{set up} in state $\phi$ and evolving according to the operator
$T$ to be \textit{measured} in state $\psi$.  This suggests we should
think of observing a system in a certain state as a time-reversed
version of setting the system up in that state.

But if time evolution by $T$ is given by an operator $E_T$,
time-evolution by $-T$ is described by an operator $E_{-T}$.  This has
groupoid and projections to $\FSN$ the same as those for $E_T$, but
the groupoid has objects labels by inverses of the labels on the
objects of $E_T$.  In $U(1)$, this inverse is the same as the complex
conjugate, so that $E_{-T}=\overline{E_T}$.  This suggests an
interpretation of the complex conjugate $\overline{\Phi}$ as a
time-reversal of the original stuff type, consistent with our
interpretation of a measurement process.

\subsubsection{Categorifying $M$}\label{sec:categorifyingM}

The operation $hom$ takes two objects in a category and yields the set
of morphisms between them.  In an enriched category, this can be
replaced by some other kind of collection of morphisms - a vector
space, for instance.  In the case that this collection is always an
object of the same category as the original objects, we have a
``$hom$-object''.  In any case $hom(-,-)$ becomes a functor into the
category in which $hom$-objects are found.

Moreover, the functor $hom(-,B)$ is a covariant functor, while
$hom(A,-)$ is contravariant - that is, a (covariant) functor to
$B^{op}$.  This seems closely analogous to the conjugate-linearity in
the complex inner product on a Hilbert space.  We may ask whether the
inner product on a Hilbert space comes from some $hom$?  In
particular, in order to make a category where morphisms are spans (see
appendix \ref{sec:slicecategory}), we need to look at the opposite
category of a $U(1)$-groupoid.  A groupoid is indistinguishable from
its opposite category after taking cardinality - but what about an
$M$-groupoid?

To make sense of this idea, we could replace a \textit{monoid} $M$
with a \textit{monoidal category}, $\mathcal{M}$. A groupoid with
objects labelled in the monoid - that is, with a function from its set
of objects to $M$, would be replaced by a groupoid $\catname{X}$ with
a functor into the monoidal category $\mathcal{M}$ - so in particular,
we would have labellings of morphisms of $\catname{X}$ with morphisms
of $\mathcal{M}$.  This combination of groupoid and functor can be
interpreted as an object in the category of ``groupoids over
$\mathcal{M}$''.

Given a categorification of $U(1)$, we could ask whether this new
setting more naturally produces the inner product we want for quantum
mechanics.  So far, though, we have not considered how to categorify
the group of phases in order to accomplish this most naturally.

\subsubsection{Non-Counting Measures and $M$-Groupoid Cardinality}

We saw in equation (\ref{eq:perturbfeynman}) and what followed that
the transition amplitudes for the perturbed harmonic oscillator are
given in terms of a sum and integral over all Feynman diagrams with
edges weighted by phases.  To do this, we had to accept the equation
at the equational level and then give it an interpretation in terms of
$U(1)$-stuff types, since we have, to date, not given any categorified
meaning for the integrals.  We can observe, however, that our notion
of cardinality for $M$-groupoids, and by extension $M$-stuff types,
used only the groupoid cardinality derived from counting measure on
sets, weights from the monoid $M$.

To give a categorified interpretation of the integral directly, we
might wish to use the fact that when $M=U(1)$, there is a measure
other than counting measure on $M$ itself.  In this case, the natural
choice is the Haar measure on the Lie group $U(1)$ - though for other
choices of $M$ there may be other natural choices.  Then a cardinality
operator for an $M$-groupoid would involve an integral involving both
$M$ and the groupoid structure.  In the case where the measure on $M$
is just set cardinality, this should reduce to the more combinatorial
definition given here.  We could hope that such a notion of
cardinality would let us give a direct categorified interpretation of
equations such as (\ref{eq:perturbfeynman}).

\subsection{Other Generalizations}

\subsubsection{Higher-Valence Stuff Operators}

We have described stuff types and operators in terms of quantum
mechanics, but it should be clear that they also have an independent
interest as algebraic objects in their own right.  Stuff types form a
categorified Hilbert space, but also a categorified algebra, since
they have a concept of multiplication in the space.

The key fact behind our approach has been that stuff types and stuff
operators form 2-categories of groupoids over either one or two copies
of $\FSN$.  Stuff operators have an action which are the equivalent of
linear operators on this 2-Hilbert space.  This action arises because
of the fact that taking a pullback over one copy of $\FSN$ under both
a stuff type and an operator removes two of the maps to this
underlying $\FSN$, and gives an object with one such map.

In fact, there is no reason why we must restrict ourselves to
groupoids with either one or two maps to $\FSN$ - the categorified
versions of vectors or matrices.  We have done so because these are
the most directly relevant to quantum mechanics, but for categorified
algebra, it makes sense to generalize to look at the equivalent of
``$p$-forms'', or ``$p$-index tensors'' (the natural inner product
obscures the difference between vectors and covectors in this
setting).  These would be groupoids with projections into $p$ copies
of $\FSN$.  Contraction of two tensors over some pair of indices would
amount to identifying the corresponding copies of $\FSN$ and taking a
pullback of the two projections into this copy.

The categorified $p$-forms could be seen as $p$-sort types: types of
structure (or stuff) which could be put on $p$ underlying sets of
different ``sorts'' of objects.

The notion of a \textit{club} described by Max Kelly in {\cite{kelly}}
can be seen as a significant generalization of this setup, where the
categories involved need not be groupoids.  There is a body of results
about these which may bear on the ideas above, and turn out to be
relevant to other physical situations.

\subsubsection{Multisort Species and QFT}

In appendix \ref{generalspecies} we refer to a description of
generalized species as functor categories between $! \mathbbm{G}$ and
$\widehat{\mathbbm{H}}$ for groupoids $\mathbbm{G}$ and $\mathbbm{H}$.
These ``species'' are the ``structure types'' of our terminology,
which correspond to the case where both groupoids are $\bldsym{1}$,
the one-element groupoid, in which case $!\bldsym{1}= \FSN$ and
$\hat{\bldsym{1}} = \Set$.  These correspond to functors from finite
sets to sets of structures which can be put on them.  We also
mentioned the 2-rig $\widehat{!\bldsym{n}}$, a.k.a $\catname{Set [
Z_1, \ldots, Z_n ]}$, the 2-rig of $n$-sort structure types: these
correspond to functors from collections finite sets of $n$ ``sorts''
(i.e. $n$ copies of $\FSN$) to sets of structures which can be put on
these.

We could reverse this point of view in the case of structure types, to
view them as faithful functors from groupoids of structures into the
groupoid of finite sets (giving the ``underlying'' set of a structure)
and then weaken the requirement that the functor be faithful to get
the more general ``stuff types'', so too we can reverse our point of
view of multisort structure types, to view them as faithful functors
from a groupoid of ``$n$-sort structures'' down to $\FSN^n$, giving
the $n$ underlying sets of each sort.  These functors will be faithful
for the same reason as in the case $n = 1$.  Weakening this
requirement would give us a notion of stuff type corresponding to
functions of more than one variable.  Defining creation and
annihilation operators on each sort of element would let us define a
Weyl algebra for $n$ sorts of particles - that is, the algebra of
operators for $n$ quantum harmonic oscillators.  This is interesting,
since a quantum field theory may be represented as a collection of
harmonic oscillators.

Replacing the various sorts of finite sets with finite sets over
monoids - in particular, over $U(1)$, as in our discussion of
$U(1)$-stuff types - we may find an elementary categorical description
of a simple QFT.  Further research in this direction may prove
fruitful.

\subsubsection{Beyond $\FSN$}

We have described the 2-Hilbert space of categorified states (stuff
types) and 2-algebra of operators (stuff operators) for the
categorified quantum harmonic oscillator in terms of some over
categories.  The 2-category $\catname{StuffTypes}$ is the slice
category of $\catname{Gpd}$ over the groupoid $\FSN$, while the
2-category $\catname{StuffOps}$ is the slice category of
$\catname{Gpd}$ over $\FSN^2$.

The groupoid $\FSN$ appears in both of these cases, in the fact that
it does allows stuff operators to act on stuff types by means of
pullbacks.  But the 2-Hilbert space structure of
$\catname{StuffTypes}$ does not depend on the fact that the groupoid
it lies over is $\FSN$: the linear structure is inherited entirely
from the direct sums of groupoids, and the inner product depends only
on the fact that any two stuff types $\Psi$ and $\Phi$ are groupoids
over the \textit{same} groupoid $\catname{G}$, and thus that we can
find a groupoid $\inprod{\Psi,\Phi}$ by taking a weak pullback over
$\catname{G}$.

Similarly, the algebraic properties of $\catname{StuffOps}$ - its
linear structure, composition, and action on $\catname{StuffTypes}$ -
derive from $\catname{Gpd}$ and the possibility of forming pullbacks.
We could derive the same structures for the 2-category of groupoids
over $\catname{G}^2$ for any groupoid $\catname{G}$.  What's more,
just as not all matrices need to be square, and linear transformations
needn't be endofunctions on some single vectorspace, we could take two
different groupoids $\catname{G}$ and $\catname{G'}$, and take the
2-category of groupoids over $\catname{G} \times \catname{G'}$.  We
could compose these in the obvious way, treating them as ``spans''
between $\catname{G}$ and $\catname{G'}$ - and also as functors
between categories of groupoids over $\catname{G}$ and $\catname{G'}$
just as stuff operators are endofunctors of groupoids over $\FSN$.

Why did we choose the specific groupoid $\FSN$ (or its $M$-coloured
counterpart) for the constructions we actually studied?  Because its
decategorification is $\mathbbm{N}$, which is the spectrum of the
number operator for the quantum harmonic oscillator.  This may suggest
how to find other groupoids $\catname{G}$ for which these
constructions have some particular physical interest.  Indeed, as
remarked before, Kelly's theory of ``clubs'' ({\cite{kelly}})
generalizes our framework by, among other things, allowing categories
which are not groupoids.  Perhaps phenomena related to groupoids can
be found which can be given a treatment like the one we have given for
the oscillator.

\section{Acknowledgements}

This work grew out of the regular Quantum Gravity seminar taught by
John Baez at UCR, notes for which are available online as
\cite{catquan}.  I would like to acknowledge his work on this subject
(some published as \cite{finfeyn}), excellent teaching, and helpful
advice and discussions in preparing this paper.  Other students in the
seminar, especially Toby Bartels, Miguel Carrion-Alvarez, Alissa
Crans, and Derek Wise also provided many useful discussions.

{\pagebreak}
\appendix\section{A Little Higher-Dimensional Algebra}\label{sec:catalgebra}

\begin{definition}
  A {\textbf{rig}}, or {\textbf{semiring}}, is a set $R$ with two
  operations, which we customarily denote by $+$ and $\cdot$, referred to as
  addition and multiplication respectively, such that $( R, + )$ is a
  commutative monoid with identity $0$, $( R, \cdot )$ is a monoid with
  identitiy $1$.  We also require that multiplication distributes over
  addition on the left and right, and $0$ is fixed under multiplication by any
  $a \in R$.
\end{definition}

\begin{definition}
  \label{2rigdef}A {\textbf{monoidal category}} $\mathcal{M} $is a
  category equipped with a functor $\otimes : \mathcal{M \times
  \mathcal{M} \rightarrow \mathcal{M}}$, a unit object $1 \in
  \mathcal{M}$, and natural isomorphisms $\alpha, \lambda, \rho$ with
  components $\alpha_{A, B, C} : ( A \otimes B ) \otimes C \rightarrow
  A \otimes ( B \otimes C )$ (the {\textbf{associator}}), $\lambda_A :
  1 \otimes A \rightarrow A$ (the {\textbf{left unit}}), and $\rho_A :
  A \otimes 1 \rightarrow A$ (the {\textbf{right unit}}), satisfying
  coherence conditions{\footnote{See, for instance, MacLane
  {\cite{maclane}}.  This definition includes slightly more than the
  definition of a rig because we here explain the generalization of
  associativity for the monoidal operation.  Also, we extend the
  ``addition'' operation to general colimits - of which coproducts,
  the equivalent of binary sums, are an example.  Otherwise, the two
  defenitions have the same form.}}.  A {\textbf{2-rig}} is a monoidal
  cocomplete category: a monoidal category $C$ which has all
  colimits, such that the functors $X \otimes - : C \rightarrow C$ and
  $- \otimes \catname{X} : C \rightarrow C$ preserve colimits for all objects $X
  \in C$.
\end{definition}

\begin{theorem}
  The category $\SetZ$ is a 2-rig whose monoidal
  operation $\otimes$ is the product $\cdot$ of structure types.
  
  \begin{proof}
    To see that $\cdot$ is a monoidal operation, note that the unit object is
    $1$, which can only be put on the empty set, in exactly one way.  Putting
    a $( 1 \cdot F )$-structure or $( F \cdot 1 )$-structure on a set $S$
    means taking $S = S \uplus \{ \}$, putting an $F$-structure on $S$ and a
    $1$-structure on $\{ \}$.  This is equivalent to putting an $F$-structure
    on $S$, so we have left and right units.  The associator $\alpha_{F, F',
    F''} : ( F \cdot F' ) \cdot F'' \rightarrow F \cdot ( F' \cdot F'' )$ is
    the natural isomorphism induced by the set isomorphism between $( A + B )
    + C$and $A + ( B + C )$.  Splitting $S$ into $A + B + C$ in these two
    ways, and putting an $F$-structure on $A$, $F'$-structure on $B$ and
    $F''$-structure on $C$ can been seen as a way of putting an $( F + F'
    )$-structure on $( A + B )$ and an $F''$-structure on $C$, but also as
    putting an $F$-structure on $A$ and an $( F' + F'' )$-structure on $( B +
    C )$.  In fact this $\alpha$ is a natural isomorphism, so $\cdot$ is
    indeed a monoidal operation.
    
    To see that $\SetZ$ is cocomplete - contains all colimits - note
    that $\Set$ is cocomplete.  Moreover, by taking colimits of the
    sets of structures on each $n$, we can find arbitrary colimits of
    objects of $\SetZ$.
    
    To see that $\SetZ$ is monoidal cocomplete, we now
    only have to have that the multiplication functors $F \cdot -$ and $-
    \cdot F$ preserve colimits for all structure types $F$ (i.e. the
    ``product'' distributes over ``sum'').
    \end{proof}
\end{theorem}

{\pagebreak}

\section{Slice Categories and 2-Categories}\label{sec:slicecategory}

  We saw that in the special case where the groupoid $Z_0$ of
  colourings was a set, the category of groupoid-coloured sets had
  objects which were maps from sets $S$ into $Z_0$, and morphisms were
  commuting diagrams like this:
  \begin{equation}
    \xymatrix{    
      S\ar_{c}[d]\ar^{\sigma}[r] & S'\ar^{c'}[dl] \\
      Z_0 & \\
    }
  \end{equation}
  commutes.

  This is an example of an ``over category''.  These are categories of
  objects ``over'' some given object.  In general, an over category
  can be constructed from any categories $\mathcal{C}$ and object $c
  \in \mathcal{C}$ by taking objects to be maps $f : a \rightarrow c$
  (for $a \in \mathcal{C}$) and morphisms to be commuting triangles.
  We will encounter this sort of construction again when we define
  $M$-sets and $M$-stuff types for a monoid $M$.  This sort of
  construction is often called a ``slice category''.  We will prefer
  the slightly more illustrative terminology ``over category''.  For
  more details, see e.g. \cite{maclane} or \cite{macmoer}.

  In general, however, groupoid-coloured sets have colours taken from
  a groupoid, which is not an object in $\catname{Set}$, so we have
  something somewhat weaker.  One way to say it is that we only have a
  forgetful functor from $Z_0$-$\catname{Set}$ to $\catname{Set}$
  where $Z_0$-sets are taken to the underlying set, and morphisms are
  taken to their underlying set bijections.  It is worth pointing out
  the relationship between this and the change of perspective between
  our original way of defining structure types as functions from
  underlying sets to the ``bundle'' viewpoint, with structured sets
  lying ``over'' their underlying sets.

  Another way to say what we have with $Z_0$-$\Set$ is that it is a
  \textit{weak} over 2-category, when we think of the sets $S$ and
  $S'$ as trivial groupoids, hence $\sigma$, $c$ and $c'$ functors.
  These are defined like over categories, but instead of morphisms
  amounting to commuting triangles, morphisms consist of natural
  isomorphisms $\alpha$ with: \begin{equation}
    \xymatrix{     
      S\ar[d] _{c}="c" \ar^{\sigma}[r] & {S'} \ar^{c'}[dl] \\ 
      Z_0 & \\ 
      \ar@{<=}^{\alpha} "c";"1,2" 
    } 
  \end{equation}

  This is exactly the definition we gave above, where the morphisms
  coming from $\alpha$ give the labels on the strands of $\sigma$.

  The formulation of over 2-categories is particularly relevant to
  stuff types, as we shall see.  First, however, we need to fill in
  some more infrastructure.

\subsection{2-Categories of Stuff Types and Stuff Operators}\label{sec:2cat}

The second - the notion of 2-categories - is just a preliminary suggestion of
a still unfinished subject of higher-dimensional categories.  However, it
turns out to be a crucial idea when we want to describe the connection between
stuff types and the quantum harmonic oscillator.  Groupoids and stuff types
naturally form a 2-category, and in section \ref{stuff-quantum-sec}, we use
the structure of this 2-category to show how the inner product on the space of
categorified states of the oscillator arises naturally.

Here we want to state and prove the important result that stuff types
and groupoids both naturally form a 2-category.  Structure types $F :
\FSN \rightarrow \Set$ formed a functor category whose morphisms were
natural transformations.  This was the ``coefficient'' viewpoint, but
for stuff types we had to take the ``bundle'' viewpoint, and defined
them as functors $F : \catname{X} \rightarrow \FSN$, for some groupoid
$\catname{X}$.  Since a groupoid is already a category, we will see
that all such objects can naturally be formed into something more than
a category.  In particular, what we will get is a 2-category:

\begin{definition}
  A {\textbf{2-category}} $\mathcal{C}$ consists of the following: a
  collection of {\textbf{objects}}, and for every pair of objects $x$
  and $y$, a category $\hom ( x, y )$ whose objects are called
  {\textbf{morphisms}} of $\mathcal{C}$ and whose morphisms are called
  {\textbf{2-morphisms}} of $\mathcal{C}$.  There must be functors
  $\hom ( x, y ) \times \hom ( y, z ) \longrightarrowlim^{m_{x, y, z}}
  \hom ( x, z )$ giving composition $( f, g ) \mapsto g \circ f$.
  There are {\textbf{identity}} morphisms $1_x \in \hom ( x, x )$ with
  {\textbf{unit laws}} which are 2-isomorphisms $\lambda$, $\rho$ from
  $1_y \circ f$ and $f \circ 1_x$ to $f$ for $f \in \hom ( x, y )$.
  There is an {\textbf{associator}}, a 2-isomorphism $\alpha_{f, g, h}
  : ( h \circ g ) \circ f \rightarrow h \circ ( g \circ f )$.  These
  satisfy coherence conditions.
\end{definition}

We are omitting here any discussion the coherence conditions.  Readers
wanting these details can consult {\cite{maclane}} for more details.

A terminological note: what we are calling a 2-category some authors,
such as {\cite{figamhy}} call a {\textit{bicategory}}, and what they call
a {\textit{2-category}} we would call a {\textit{strict 2-category}}, where
the associator and unit laws are identities.  We adopt this convention
because the non-strict case seems to be the more generally useful one,
and deserves a nomenclature which generalizes naturally.

Now, we can make the following observation:

\begin{theorem}
  The collection of all categories, $\catname{Cat}$, naturally forms a
  2-category whose morphisms are functors between categories, and whose
  2-morphisms are natural transformations between functors.  The collection of
  groupoids, $\catname{Gpd}$, is a full sub-2-category of
  $\catname{Cat}$.  In fact, these are strict 2-categories.  
\end{theorem}

This sets up a helpful way of looking at stuff types: we have
described them as ``groupoids over $\FSN$'' - that is, functors from
groupoids $\catname{X}$ to the groupoid $\FSN$.  In particular, since
groupoids from a 2-category, of which $\FSN$ is an object, we can
describe a 2-category of stuff types, namely that of {\textit{groupoids
over $\FSN$}}.  This is a 2-categorical version of an ``over
category''.  This sort of structure is explained in further detail in
appendix \ref{sec:slicecategory}.

We have just described stuff types as functors into $\FSN$ in
$\catname{Gpd}$, so it is natural to ask whether other functors in
$\catname{Gpd}$ are also of interest as further generalizations of
structure types.  In section \ref{generalspecies}, we briefly describe
some work in this direction.  

\begin{definition}\label{stufftype2cat}
  The weak 2-category $\catname{StuffTypes}$ has as objects diagrams
  in $\catname{Gpd}$ of the form
  $\catname{X}\rightarrowlim^{\Psi}\FSN$ (Denoted $(\catname{X},\Psi)$, or
  just $\catname{X}$ or $\Psi$ for short whenever the meaning is
  clear).  Given two objects $(\catname{X_1},\Psi_1)$ and
  $(\catname{X_1},\Psi_1)$, $hom(\Psi_1,\Psi_2)$ has as morphisms
  functors $F:\catname{X_1}\rightarrow \catname{X_2}$ together with a
  natural isomorphism $\alpha$ such that the diagram
\begin{equation}
  \xymatrix{    
    \catname{X_1}\ar[d] _{\Psi_1}="t" \ar^{F}[r] & {\catname{X_2}}\ar^{\Psi_2}[dl] \\
    \FSN & \\
      \ar@{<=}^{\alpha} "t";"1,2"
  }
\end{equation}
  commutes up to $\alpha$.  Given a pair $F$ and $G$ of such morphisms between $X_1$ and $X_2$, the 2-morphisms between them are the natural transformations $\nu$
  between the functors $F$ and $G$ for which the resulting diagram commutes.
\end{definition}

\begin{theorem} The construction given for $\catname{StuffTypes}$
gives a well-defined 2-category.
  \begin{proof} The collections $\opname{hom}(\Psi_1,\Psi_2)$ involve
  functors from $\catname{X_1}$ to $\catname{X_2}$ in $\catname{Gpd}$.
  These are closed under composition.  A morphism in
  $\catname{StuffTypes}$ also includes a natural transformation
  $\alpha$, and these are again closed under composition.  If two
  composable functors $F_1$ and $F_2$ between groupoids make the
  triangles over $\FSN$ commute up to natural isomorphisms $\alpha_1$
  and $\alpha_2$, then $F_2 \circ F_1$ does the same, up to $\alpha_1
  \circ \alpha_2$.  So in particular, the obvious notion of
  composition is well defined, and in fact the
  $\opname{hom}(\Psi_1,\Psi_2)$ are categories.

  Identity morphisms are inherited from $\catname{Gpd}$, and obviously
  make the corresponding triangles commute.  The unit laws and
  associator are just these identity 2-morphisms, so we have a strict
  2-category.
  \end{proof}
\end{theorem}

The construction for stuff operators is similar:

\begin{definition}\label{stuffops2cat}
  The 2-category $\catname{StuffOps}$ has as objects diagrams in
  $\catname{Gpd}$ of the form
  $\FSN\leftarrowlim^{p_1}T\rightarrowlim^{p_2}\FSN)$ (Denoted
  $(T,p_1,p_2)$, or just $T$, for short).  Given two objects
  $T$ and $T'$, $hom(T,T')$ has as morphisms functors
  $F:T\rightarrow T'$ making the diagram
\begin{equation}
  \xymatrix{    
    T\ar_{p_1}[d]\ar_{p_2}[dr]\ar^{F}[r] & {T'}\ar^{p'_1}[dl]\ar^{p'_2}[d] \\
    \FSN & \FSN \\
  }
\end{equation}
  commute up to two natural isomorphisms.  The 2-morphisms are the
  natural transformations $\nu$ between such functors $F$ and $G$
  which make the resulting diagram commute.

\end{definition}
Where we have omitted the detailed diagram for the naturality squares.
It is substantially similar to that for stuff types, in section
\ref{stufftypes}.  We also get a result similar to that for stuff types:

\begin{theorem} The construction given for $\catname{StuffOps}$
gives a well-defined strict 2-category.
  \begin{proof} 
  The proof that this is a 2-category is similar to that for stuff types.
  \end{proof}
\end{theorem}

The algebraic structure of $\catname{StuffOps}$ is of interest.  It
is, in fact, the equivalent of an algebra - having addition, scalar
multiplication (by groupoids) and internal multiplication in the form
of composites.  But since groupoids do not have cardanilities in a
field, we will just point out that if we ignore 2-morphisms, it is a
category, and in fact:

\begin{theorem} The category $\catname{StuffOps}$ (disregarding
2-morphisms) is a 2-rig, where the monoidal operation is composition.
\begin{proof}
First, the $\catname{StuffOps}$ is cocomplete because $\catname{Gpd}$
is, and so any colimit in $\catname{Gpd}$ becomes one in
$\catname{StuffOps}$.  The monoidal operation given by composition
gets all the required natural isomorphisms from those in the weak
pullback.
\end{proof}
\end{theorem}

We have only sketched the main ideas of these proofs, of course (in
particular, we have not even stated the necessary coherence
conditions, let alone proved they are satisfied).  We leave these
details for the interested reader.  However, this finally gives us a
clear description of the categorified version of the algebra of
operators on formal power series.

Analogous results hold for $M$-stuff types.

{\pagebreak}

\section{Categorical Approaches to Generalizing Species}

We have chosen in this paper to generalize Joyal's notion of structure types
in a way which makes use of the classification of functors and their levels of
forgetfulness.  By saying that a structure type is a possibly forgetful
functor which forgets at most {\textit{structure}}, and possibly only
properties, or nothing, we find that it is possible to generalize this to a
\textit{stuff type}, described in section \ref{stufftypes}.  This also hints
at further generalizations which will be possible if we allow ourselves to
consider functors between higher dimensional categories, as described in
section \ref{sec:forgetful}, so that there are more possible degrees of
forgetfulness of functors, and therefore a hierarchy of ``types'' given by
functors which forget ``meta-stuff'' of various degrees.  However, this is not
the only possible direction in which to take the notion of structure
type.  We describe here another direction.

\subsection{Generalized Species} \label{generalspecies}

A generalization of structure types is described in Fiore, Gambino,
and Hyland {\cite{figamhy}} sheds some light on the choice of categories we have made in
defining structure types.  In that paper, conventional structure types are
referred to as {\textit{species}}, and a generalization is developed to $(
\mathbbm{G}, \mathbbm{H} )$-species for arbitrary small groupoids
$\mathbbm{G}$ and $\mathbbm{H}$, which provides, for finite sequences of
{\textit{$\mathbbm{G}$-objects}}, an {\textit{$\mathbbm{H}$-variable set of
structures}} over them.  Structure types are then
$(\bldsym{1},\bldsym{1})$-species, where $\bldsym{1}$ is the groupoid
with one object and identity morphism.

To explain this generalization, we need some terminology.

\begin{definition}
  Suppose $\mathbbm{G}$ is a small groupoid.  Then $! \mathbbm{G}$, the
  {\textbf{free symmetric monoidal completion}} of $\mathbbm{G}$, is the
  smallest symmetric monoidal groupoid containing $\mathbbm{G}$.  We define
  $\widehat{\mathbbm{G}}$, the {\textbf{free cocompletion}} of $\mathbbm{G}$
  is the smallest cocomplete category containing $\mathbbm{G}$.
\end{definition}

Now, we can make an analogy between the creation of a 2-rig from a category
using these constructions and the creation of a rig from a set using the
operation of taking a free abelian group on a set, and the operation of taking
the free monoid on a set.  If we start with a set of generators $S$, and then
take the free Abelian group on $S$, $\mathbbm{Z}[S]$, and then take the
free monoid on $\mathbbm{Z}[S]$, we get a rig, and this is isomorphic to
the rig we get if we take these freely generated structures in the reverse
order.  In particular, if $S = \phi$, the rig we get is $\mathbbm{N}$, if $S =
\{x\}$, we get $\mathbbm{N}[x]$, the free rig on one generator, and so
on.  A similar construction is possible for groupoids (and indeed categories).

To see how this applies to species, we first note that the free symmetric
monoidal completion of a groupoid consists of ``families'' of
$\mathbbm{G}$-objects, whose objects are tuples of objects from $\mathbbm{G}$
and whose morphisms are braids between tuples, with strands labelled by
morphisms of $\mathbbm{G}$.  So in particular, if $\mathbbm{G}$ is the
groupoid $\bldsym{1}$, with one object and only the identity morphism, we
find that $! \mathbbm{G} \cong \FSN$. Moreover, the free
cocompletion of $\mathbbm{G}$ is equivalent to the functor category $\hom (
\mathbbm{G}^{\opname{op}}, \Set )$ of presheaves on
$\mathbbm{G}$ (see, for instance, MacLane \& Moerdijk
{\cite{macmoer}}, I.5, Prop 1).  So in particular, since $\mathbbm{G}$
is a groupoid, and equivalent to its opposite, we have that
$\widehat{! \mathbbm{G}} \cong \hom ( !
\mathbbm{G}, \Set )$.

In the case where $\mathbbm{G} =\bldsym{1}$, we have
$\widehat{!\mathbbm{G}} \cong \hom ( \FSN, \Set ) = \SetZ$.  So the
2-rig of structure types can be seen as the freely generated 2-rig on
one generator.  We may think of this generator as being the basic
``object'', or ``one-element set''.  The first natural extension to
consider is when $\mathbbm{G} =\bldsym{n}$, the groupoid with $n$
objects having only the identity morphisms.  The 2-rig
$\widehat{!\bldsym{n}}$ is also called $\catname{Set [ Z_1, \ldots,
Z_n ]}$, and gives what are called {\textit{multisort species}}.  These
can be described as 2-rigs of structures which can be put on sets of
elements of $n$ different sorts.  For other groupoids $\mathbbm{G}$,
we get different notions of species, many of which appear in various
contexts in the literature, for instance ({\cite{BLL}}).

The generalization of species considered by Fiore, Gambino, and Hyland
{\cite{figamhy}} is the 2-rig $\hom ( ! \mathbbm{G},
\widehat{\mathbbm{H}} )$, for $\mathbbm{G}$and $\mathbbm{H}$ some
small groupoids.  The various examples of $\mathbbm{G}$-species
mentioned above are all seen as functors into $\Set =
\hat{\bldsym{1}}$, so $\mathbbm{H} =\bldsym{1}$.  The 2-rig $\hom ( ! 
\mathbbm{G}, \widehat{\mathbbm{H}} )$ is, in particular, the category
of functors from $! \mathbbm{G}$, the category of {\textit{families of
$\mathbbm{G}$-objects}} to $\mathbbm{\hat{H}}$, the category of
{\textit{$\mathbbm{H}$-variable sets}} - presheaves over $\mathbbm{H}$.

Indeed, it is possible to define a 2-category of species between
groupoids.  In this setting, the category of functors from families of
$\mathbbm{G}$-objects into $\mathbbm{H}$-variable sets plays the role
of $\hom ( \mathbbm{G}, \mathbbm{H} )$, and the objects in the
2-category are small groupoids.

{\pagebreak}

\end{document}